\documentclass[11pt,a4paper]{article}

\usepackage{graphicx} % Required for inserting images.
\usepackage[margin=25mm]{geometry}
\usepackage[font=sf]{caption} % Changes font of captions.

\usepackage{amsmath}
\usepackage{amsfonts}
\usepackage{amssymb}
\usepackage{siunitx}
\usepackage{verbatim}
\usepackage{hyperref} % Required for inserting clickable links.
\usepackage{natbib} % Required for APA-style citations.
\usepackage{times}
\usepackage{latexsym}
\usepackage{threeparttable}  % For table notes
\usepackage{booktabs}        % For better tables (optional)
\usepackage{tabularx}
\usepackage{siunitx}
\usepackage{threeparttable} 
\usepackage{booktabs}
\usepackage{afterpage}
\usepackage{array}
\usepackage{siunitx}
\usepackage{graphicx}        % For \resizebox command
\usepackage{textcomp}
\usepackage{siunitx}
\usepackage{float}
\usepackage{caption}
\usepackage{subcaption}
\usepackage{multirow}
\usepackage{pifont}

\usepackage{cleveref}
\crefname{figure}{Fig.}{Figs.}

\title{A unified solution framework for truck-and-drone routing problems}
\author{Ke Xu and John Gunnar Carlsson}
\date{}
\begin{document}
\maketitle

\begin{abstract}
Coordinated truck-and-drone routing integrates the high capacity and range of ground vehicles with the flexible routing and speed of drones, enabling simultaneous service. Increasingly applied in last-mile delivery, this synchronization helps reduce completion time and operational costs. To improve its efficiency, various coordination modes between trucks and drones have been proposed. Each mode accommodates diverse operational constraints tailored to particular delivery requirements. In existing work, a slight change in the structural framework or operational characteristics could generate a totally different problem variant, which often requires the design of specialized algorithms. Consequently, under the requirement of maintaining structural validation and adapting to multiple operational features, this paper presents a unified three-phase solution framework based on the Lin-Kernighan-Helsgaun algorithm to  solve a wide family of truck-and-drone routing problems. To validate its flexibility and effectiveness, we carry out numerical experiments on three problem variants: the Flying Sidekick Traveling Salesman Problem (FSTSP), the Traveling Salesman Problem with Multiple Drones (TSP-mD), and the Vehicle Routing Problem with Drones (VRP-D), benchmarking each against an effective algorithm. Computational results show that the framework can closely match optimal solutions on small-size instances and even improve the best-known solutions for several medium-size instances. Moreover, additional extensions are discussed to further highlight its versatility.

\vspace{1em}

\end{abstract}

% \begin{keyword}
% vehicle routing problem \sep drone \sep last-mile delivery
% \end{keyword}

\section{Introduction} \label{sec:introduction}
% ``----------Penalty function and resulting figures can be found in \href{https://drive.google.com/drive/folders/13UUiV9mOAK3c969MtVokVwVW2V2fpC7Q}{LKHexperiments1}." vehicle routing and drone scheduling...

% The field of logistics and transportation has undergone transformative advancements over the past decade, with the integration of unmanned aerial vehicles (UAVs), commonly known as drones, into traditional delivery systems emerging as a groundbreaking development. 
%Besides, inefficient last-mile delivery gives rise to a multitude of negative impacts, which can be categorized into three dimensions: economy (vehicle miles traveled VMT and energy consumption), environment (carbon emissions and air pollution), and social (traffic congestion and customer satisfaction) \cite{ranieri2018review}.

During the past decade, exploring coordinated operations of trucks and drones or ground robots has emerged as a significant advancement in logistics research (see, e.g., in \cite{chung2020optimization}, \cite{macrina2020drone}, \cite{li2021ground}, \cite{moshref2021applications}, \cite{dukkanci2024facility}). This coordination has been increasingly applied in last-mile delivery, which is the least efficient phase of distribution networks, contributing up to 75\% of total supply chain costs \cite{o2025optimising}. To optimize last-mile delivery paradigms, \cite{murray2015flying}, \cite{agatz2018optimization}, \cite{campbell2017strategic}, and \cite{carlsson2018coordinated} laid the groundwork. Prominent industry players, such as Amazon \cite{Palmer2020}, Google \cite{Vallance2023}, and DHL \cite{abc2014} are also actively developing cooperative truck-drone delivery systems \cite{Alferez2023}.

However, this integration of drones into conventional Vehicle Routing Problems (\textbf{VRPs}) not only introduces effective logistics but also extends an already diverse problem space into a more intricate domain, involving different variants between trucks and drones. Each variant implements various real-world features (e.g., objective functions of minimizing time versus cost, drone launch and recovery restrictions) on a coordination mode (e.g., one truck + one drone, truck-drone groups). To effectively solve different variants, multiple solution methods have been developed, ranging from metaheuristics such as genetic algorithms \cite{chiang2019impact} and adaptive large neighborhood search (ALNS) \cite{moshref2020design}, to exact methods like branch-and-bound \cite{dell2021algorithms}, as well as continuous approximation techniques \cite{carlsson2018coordinated} and hybrid strategies \cite{luo2021hybrid}.

This complexity of algorithm development, and the fragmented solution space for its variants, motivates us to construct a unified heuristic framework capable of flexibly solving \textbf{truck-and-drone routing} problems. In pursuit of this goal, we introduce \textbf{ penalty functions} developed in an advanced heuristic solver, Lin-Kernighan-Helsgaun (\textbf{LKH-3}), which is well designed for a wide range of constrained VRPs.

%Note that the Lin-Kernighan-Helsgaun algorithm was extended by Keld Helsgaun in %\cite{cook2024constrained}
%to solve the last-mile problem, but they only considered it at the level of a %single driver without truck-drone synchronization. 

The efficiency of our framework is validated through numerical experiments in three variants, involving coordinated structures such as (1) one truck + one drone, (2) one truck + multiple drones, (3) truck groups where each group refers to coordination between a truck and its designated drone, and operational constraints in reality. Each variant is addressed by a specific penalty function, which can not only derive structural routing decisions such as synchronization times and locations of truck(s) and drone(s), but also address real-world operational constraints, including time windows and capacity limitations, by adding the corresponding penalty terms. Therefore, this framework has the potential to solve a wide range truck-and-drone routing problems, encompassing various combinations of structural and operational characteristics.

Before proceeding to the proposed solution framework and the examined cooperation variants, a brief outline of truck-and-drone routing problems is provided. In general, they are characterized by three fundamental elements: (1) a depot where vehicles start and end their routes, (2) a set of customers, each requiring a delivery, and (3) synchronization between truck and drone. Pioneering work was introduced by \cite{murray2015flying}, who proposed two distinct variants. The first is the Flying Sidekick Traveling Salesman Problem (FSTSP), in which a truck and a drone cooperate to deliver packages. In this model, the truck must continue its route after launching the drone and cannot wait at the launch point for retrieval, thereby prohibiting \textbf{loops} in drone trajectory. An alternative version, often referred to as TSP-D in the subsequent literature (see \cite{agatz2018optimization}, \cite{bouman2018dynamic}, \cite{poikonen2019branch}, \cite{roberti2021exact}), relaxes this \textbf{no-loop restriction} by allowing the truck to remain at the launch site for the drone's return, and enhances flexibility by revisiting nodes and disregarding operation time, etc. The explicit distinction between FSTSP and TSP-D was provided in \cite{mahmoudinazlou2024hybrid}.
Meanwhile, the second variant named Parallel Drone Scheduling TSP (PDSTSP) in \cite{murray2015flying} did not involve synchronization: a truck is responsible for distant customers, while drones are independently in charge of customers near the distribution center, resulting in separate routes for drones and the truck. Afterwards, truck-and-drone routing problems have been increasingly enriched by incorporating more than one drone/truck and introducing operational constraints for realistic scenarios. Representative constraints that our adaptable solution framework can readily resolve are as follows.

\begin{enumerate}
\item Drone flight range \citep{ha2018min}: Drones have limited battery capacity, restricting the maximum distance they can travel before recharging.
\item Drone operation time \citep{ha2020hybrid}: Handling times for drone launch and retrieval operations.
\item Customer service time \citep{murray2020multiple}: Serving time for each customer.
\item Drone payload \citep{sacramento2019adaptive}: Drones have lower carrying capacity than trucks, implying that some heavy-demand customers can only be serviced by truck.
\item Non-customer points \citep{carlsson2018coordinated}: Trucks and drones can meet at locations other than customer sites, increasing routing flexibility.
\item Customer time windows \citep{kuo2022vehicle}: Deliveries must be made within time intervals.
\end{enumerate}

% This work aims to address this gap by proposing a novel framework that optimally integrates truck-and-drone operations for efficient last-mile delivery. Our approach not only considers conventional optimization objectives such as makespan and routing costs but also explores the impact of synchronization policies and multi-drone configurations. Through comprehensive computational experiments, we evaluate the performance of our model against existing benchmarks, demonstrating its effectiveness in reducing delivery times and operational costs. 

\begin{figure}[H]
    \centering    \includegraphics[width=1.0\textwidth]{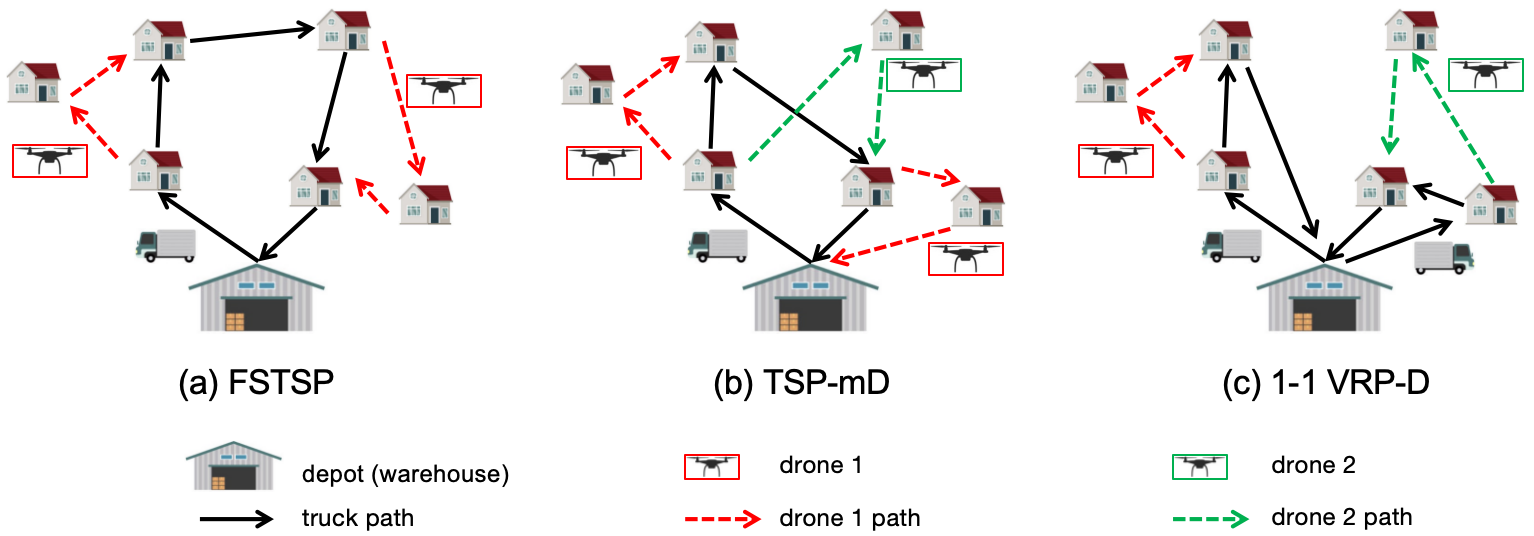}
    \caption{Comparison of delivery schedules of three problem variants.}
    \label{fig:schedules}
\end{figure}

Considering different structures combined with various constraints, a versatile approach to address the full range of truck-and-drone routing problems is highly valuable. Therefore, our flexible solution framework is proposed based on \textbf{LKH-3} and validated through three problem variants: the Flying Sidekick Traveling Salesman Problem (FSTSP, one truck + one drone), the Traveling Salesman Problem with Multiple Drones (TSP-mD, one truck + homogeneous drones), and the Vehicle Routing Problem with Drones (1-1 VRP-D, each group implies coordination between one truck and one drone), as exhibited in Fig.~\ref{fig:schedules}. 

The remainder is organized as follows. Section 2 describes the three-phase solution framework. Section 3 introduces benchmark instances. Section 4 presents literature review and problem description for the three truck-and-drone routing variants. Section 5 reports numerical results to demonstrate the flexibility and efficiency of the proposed framework. Finally, some extensions and future research directions are discussed in Section 6.

\section{Methodology} \label{sec:methodology}

\begin{figure}[H]
    \centering    \includegraphics[width=0.8\textwidth]{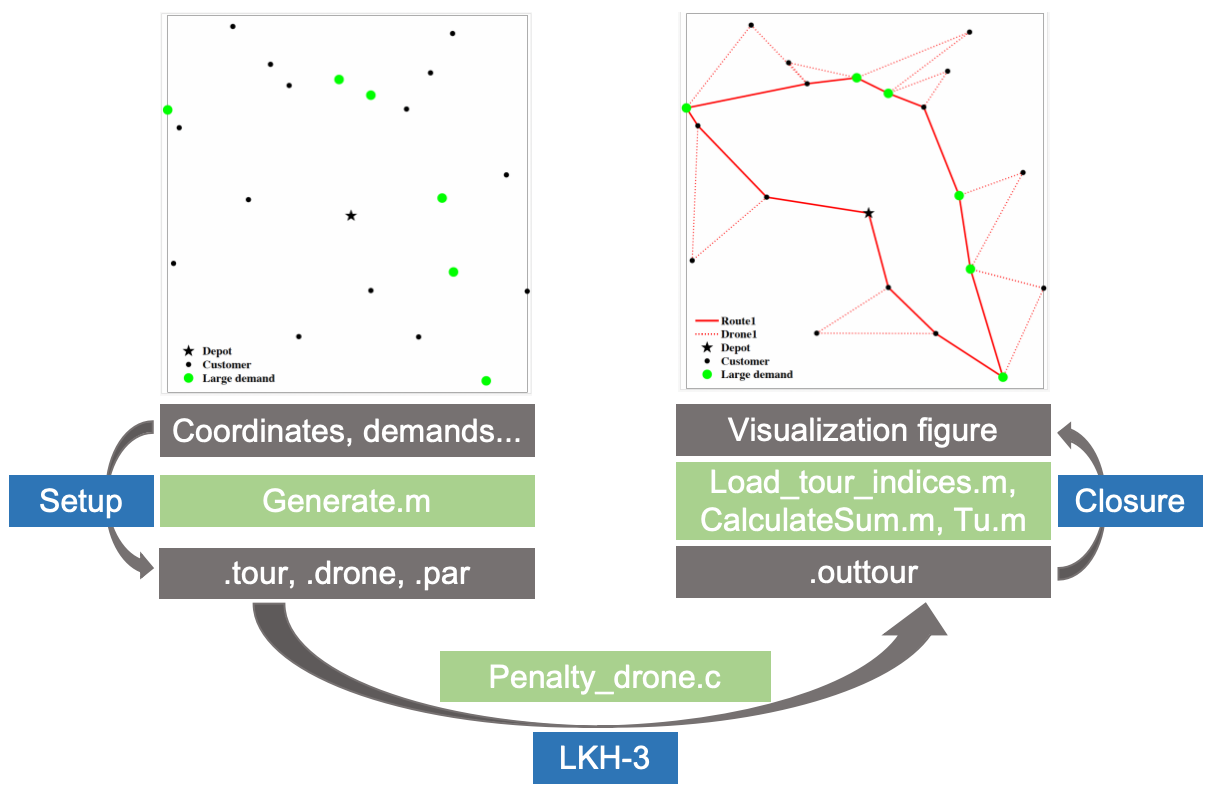}
    \caption{High-level flowchart of the solution framework.}
    \label{fig:method}
\end{figure}

As shown in Fig.~\ref{fig:method}, the proposed solution framework relies on three phases: (1) Setup in Matlab, (2) Execution in LKH-3, and (3) Closure in Matlab. The first two phases are detailed in the following two subsections. Note that we utilize the implementation architecture of the LKH-3 solver. Its execution instructions can be found in the report \cite{helsgaun2017extension}. 

\begin{enumerate}
\item Setup: ``Generate.m" reads basic information such as the coordinates for depot and customers, and prepares three files including ``\texttt{.tour}, \texttt{.drone}, \texttt{.par}" for LKH-3. Note that \texttt{.tour} represents an initial feasible solution that is not absolutely necessary to use, \texttt{.drone} presents one readable truck-and-drone routing problem, \texttt{.par} is a parameter file for LKH-3 where the names of reading and output files are specified. 
\item Execution: By executing the \texttt{.par} file, the penalty function ``Penalty\_drone.c" designed for the routing problem is employed to handle its \textbf{characteristic structure} and \textbf{operational constraints}. Thus, all problem restrictions are satisfied in the output solution \texttt{.outtour}, which achieves the minimum of the problem objective function as well.
\item Closure: ``Load\_tour\_indices.m" and ``CalculateSum.m" reads the solution \texttt{.outtour} and validates its correctness by calculating the objective again, respectively. Moreover, ``Tu.m" is called to draw an illustrative figure. Fig.~\ref{fig:method} (right) shows a loop-free drone trajectory, resulting from the restriction that trucks cannot wait at launch points for retrieval.
\end{enumerate}

\subsection{Details on setup in Matlab}
\begin{figure}[H]
    \centering
    \setlength{\tabcolsep}{1pt} % Space between columns
    \begin{tabular}{c}        \includegraphics[width=0.66\textwidth]{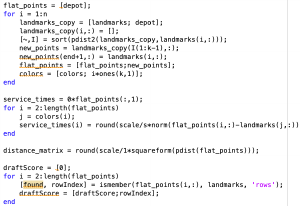} \\
        % \multicolumn{1}{c}{(a) ALNS.}
    \end{tabular}
    \caption{Input files generation in Matlab.}
    \label{fig:fstsp_exe3}
\end{figure}

This subsection describes the setup phase for preparing the input files, i.e., \texttt{.drone} and \texttt{.par} for LKH-3. These files are generated by ``Generate.m" in Matlab, as shown in Fig.~\ref{fig:fstsp_exe3}. The primary resulting data, \texttt{flat\_points}, is generated by the first loop of the seven-line code. Each customer $i$ generates a set of \texttt{new\_points} containing $k$ elements, where $k$ is an input parameter. The first two lines construct \texttt{landmarks\_copy}, including all customer locations except customer $i$. The next three lines determine \texttt{new\_points} by selecting the $k-1$ nearest points from \texttt{landmarks\_copy}, and adding the customer $i$. Then, the \texttt{new\_points} of $k$ elements is added to \texttt{flat\_points}. Meanwhile, these $k$ points are assigned with the same \texttt{color} $i$, which is added to \texttt{colors}. Thus, there are $n$ nonzero colors, with \texttt{color} $0$ assigned to the depot.

Examining those $k-1$ distinct points from \texttt{landmarks\_copy}, each can serve as a launch/retrieval \textbf{(L/R)} site for customer $i$ in LKH-3. Thus, two points among them will be designated: one for launch and the other for retrieval, thus preventing the \textbf{loop} option of waiting at a launch point for retrieval. Furthermore, by adjusting this code block, we can flexibly allow loops (\cite{agatz2018optimization}, \cite{bouman2018dynamic}, \cite{poikonen2019branch}, \cite{roberti2021exact}) and incorporate non-customer points (\cite{carlsson2018coordinated}) for drone operations, which will be discussed in Section 6.1 \hyperref[sec:extensions]{Extensions}.

Subsequently, Fig.~\ref{fig:fstsp_exe3} exhibits \texttt{service\_times} and \texttt{distance\_matrix}, which correspond to the travel time from an L/R point to its target customer, and the travel time between \texttt{flat\_points} locations, respectively. They may be used for objective calculation in LKH-3. Moreover, a location may serve as an L/R candidate for multiple customers, thus one physical point may appear in the \texttt{new\_points} sets of multiple customers. Therefore, we identify the physical locations of \texttt{flat\_points} by assigning distinct scores using \texttt{draftScore}. Thus, $n$ nonzero scores are given to $n$ customers, along with \texttt{draftScore} $0$ assigned to the depot. These scores enable LKH-3 to restrict or allow the truck to revisit a point once it departs, which will be investigated in Section 5 \hyperref[sec:experiments]{Experimental results}.

% relaxes line-breaking rules
\begin{sloppypar}
As a result, the problem information includes \texttt{distance\_matrix}, \texttt{colors}, \texttt{service\_times}, and \texttt{draftScore}. They are documented as EDGE\_WEIGHT\_SECTION, CTSP\_SET\_SECTION, SERVICE\_TIME\_SECTION, and DRAFT\_LIMIT\_SECTION in the \texttt{.drone} file. Further details are provided in \hyperref[app:Inputfiles]{Appendix A}.

\end{sloppypar}

\subsection{Details on execution in LKH-3}
In the second phase, each variant of the truck-and-drone routing problem can require a unique penalty function ``Penalty\_drone.c". Within the penalty function, various \textbf{hard} penalty terms (e.g., larger than 1,000,000) are added if structural settings and operational constraints are violated. In addition, a relatively small penalty term (e.g., smaller than 1,000,000) is added as a \textbf{soft} term corresponding to the objective function value. In the generated solution file \texttt{.outtour}, therefore, the feasible tour having the minimum objective is documented in the format ``A\_B", where A reflects its objective.

Significantly, ``Generate.m" in the first phase and ``Penalty\_drone.c" in LKH-3 can be modified to solve various truck-and-drone routing problems, which will be discussed in Section 6.1 \hyperref[sec:extensions]{Extensions}. To verify the flexibility and effectiveness of the proposed three-phase framework, we evaluate three variants and benchmark each against an effective algorithm. 

\begin{enumerate}
\item Flying Sidekick Traveling Salesman Problem (FSTSP), compared to the Hybrid Genetic Algorithm (HGA) developed in \cite{ha2020hybrid}. Additionally, we assess performance on small instances against the BKN exact results reported in \cite{boccia2023new}.
\item Traveling Salesman Problem with Multiple Drones (TSP-mD), compared to the Adaptive Large Neighborhood Search (ALNS) metaheuristic algorithm proposed by \cite{tu2018traveling}.
\item Vehicle Routing Problem with Drones (VRP-D), compared to the results reported in \cite{sacramento2019adaptive}, who also developed an ALNS approach.

\end{enumerate}
We present the problem details and the corresponding numerical results in Section 4 and 5, respectively.

\section{Instances} \label{sec:instances}

This section provides seven instance sets for truck-and-drone routing benchmarks, as listed in Table~\ref{tab:instance_sets}.

\begin{table}[H]
    \centering
    % \tiny
    % \scriptsize
    % \footnotesize  % Slightly smaller to better fit
    % \small
    % \normalsize  % (default size) 
    % \large
    \caption{Seven instance sets}
    \label{tab:instance_sets}
    \begin{tabular*}{\textwidth}{@{\extracolsep{\fill}}l@{\;}c@{\;}c@{\;}c@{\;}c@{\;}c@{}}
        \toprule
        \textbf{Set} & \textbf{Reference Paper} & \textbf{Dataset Link} &
        \textbf{Aim Problem} & \textbf{\#Ins} &
        \textbf{Customer demand}\\
        \midrule
        Set 1 & \citet{murray2015flying} & \href{http://www.or.unimore.it/site/home/online-resources/exact-models-for-the-fstsp.html}{Link1} & FSTSP; PDSTSP & 72; 240 & Too heavy denoted as 1 \\
        Set 2 & \citet{ha2018min} & \href{https://orlab.com.vn/download/instances-for-the-traveling-salesman-problem-with-drone-tsp-d/}{Link2} & FSTSP & 65 & Too heavy denoted as 1 \\
        Set 3 & \citet{agatz2018optimization} & \href{https://doi.org/10.5281/zenodo.1204676}{Link3}, \href{https://github.com/pcbouman-eur/Drones-TSP}{Link4} & TSP-D & 140 & - \\
        Set 4 & \citet{poikonen2019branch} & \href{https://mario.ruthmair.at/?page_id=226}{Link5} & TSP-D & 200 & - \\
        Set 5 & \citet{murray2020multiple} & \href{https://github.com/optimatorlab/mFSTSP}{Link6} & TSP-mD & 100 & Parcel weights in Lbs \\
        Set 6 & \citet{sacramento2019adaptive} & \href{https://zenodo.org/records/2572764}{Link7} & VRP-D & 112 & Parcel weights in kgs \\
        Set 7 & \citet{merchan20242021} & \href{https://registry.opendata.aws/amazon-last-mile-challenges/}{Link8} & real-world & 6112 & -\\
        \bottomrule
    \end{tabular*}
\end{table}

Sets 1-2 focus on FSTSP instances derived from \cite{murray2015flying} and \cite{ha2018min}. They denote drone-ineligible customers as ``1" due to heavy packages. Note that the 72/240 FSTSP/PDSTSP instances are generated by applying two endurance parameters to 36/120 unique location configurations. The PDSTSP instances can be adapted to FSTSP \cite{boccia2021column}. Of the 65 Set2 instances in \cite{ha2018min}, 60 were tested by the same authors in \cite{ha2020hybrid}, showing improvement. 

Sets 3-4 represent TSP-D instances from \cite{agatz2018optimization} and \cite{poikonen2019branch}. As mentioned in \hyperref[sec:introduction]{Introduction}, this generalized coordination between one truck and one drone removes many restrictions of the FSTSP, e.g., by allowing drones to end their flights at their respective launch nodes. Note that these two sets provide only coordinates for the depot and customers, lacking customer demands. 
% source code for set3 reference paper: https://github.com/pcbouman-eur/Drones-TSP, instances: https://github.com/pcbouman-eur/TSP-D-Instances

Set 5 from \cite{murray2020multiple} considers heterogeneous drones with parcel weights in pounds. Set 6, generated by \cite{sacramento2019adaptive}, provides VRP-D instances with parcel weights in kilograms. Set 7 from \cite{merchan20242021} offers real-world instances inspired by Amazon’s last-mile delivery challenges, providing practical test cases for actual applications.

In Section 5, we analyze FSTSP and TSP-mD using Poikonen instances from Set 4 (6, 9, 19, 29, and 39 customers), and VRP-D using Sacramento instances from Set 6 (6, 10, 12, and 20 customers).

\section{Problem variants} \label{sec:variants}

As outlined in \hyperref[sec:methodology]{Methodology}, to validate the effectiveness of our solution framework, we compare its performance against three algorithms, each of which is specifically tailored to one of the variants: FSTSP (one truck + one drone), TSP-mD (one truck + multiple drones), and VRP-D (truck groups). Thus, each variant's related work and representative features are studied in this section. The numerical results will be compared in Section 5.

\subsection{FSTSP} \label{sec:fstsp}
The Flying Sidekick Traveling Salesman Problem (FSTSP) is the most fundamental variant of truck-and-drone routing problems, incorporating one drone to assist a truck. The drone can be launched from the truck at one node, serve a customer, and then rendezvous with the truck at a different node. Meanwhile, during drone flight, the truck continues to serve customers.

\subsubsection{Related work}

\begin{table}[H]
\centering
\caption{Papers related to FSTSP}
\label{tab:fstsp_papers}
\begin{threeparttable}
\small
\begin{tabular*}{\textwidth}{@{\extracolsep{\fill}}l@{\;}c@{\;}c@{\;}c@{\;}c@{\;}c@{\;}l@{}}
\toprule
\textbf{Authors} & \textbf{L/R} & \textbf{Non-cus} & \textbf{Multi-} & \textbf{Model} & \textbf{Objective} & \textbf{Solution Method} \\
 &  & \textbf{points} & \textbf{visit} &  &  & \\
\midrule
\cite{murray2015flying} & D &  &  & MILP & makespan & Heuristic \\
\cite{ha2015heuristic} & D &  &  &  & makespan & \begin{tabular}[t]{@{}l@{}}Cluster first - route second;\\ Route first - cluster second\end{tabular} \\
\cite{ponza2016optimization} & D &  &  & MILP & makespan & Simulated annealing \\
\cite{ha2018min} & D &  &  & MILP & cost & TSP-LS and GRASP \\
\cite{yurek2018decomposition} & D &  &  & MILP & makespan & Iterative algorithm \\
\cite{de2018randomized} & D &  &  & MILP & makespan & \begin{tabular}[t]{@{}l@{}}Randomized variable\\ neighborhood descent\end{tabular} \\
\cite{marinelli2018route} & D & \ding{51} &  & IP & cost & Greedy heuristic \\
\cite{wang2020cooperative} & D &  &  &  & \begin{tabular}[t]{@{}l@{}}makespan \&\\ cost\end{tabular} & Improved NSGA-II \\
\cite{de2020variable} & D &  &  &  & makespan & HGVNS \\
\cite{ha2020hybrid} & D &  &  &  & \begin{tabular}[t]{@{}l@{}}makespan /\\ cost\end{tabular} & Hybrid genetic algorithm \\
\cite{gonzalez2020truck} & D &  & \ding{51} & MILP & makespan & Iterated greedy heuristic \\
\cite{schermer2020b} & D &  &  & MILP & makespan & B\&C \\
\cite{dell2021algorithms} & D &  &  &  & makespan & B\&B \\
\cite{dell2021drone} & D &  &  &  & makespan & B\&C \\
\cite{el2021parcel} & D &  &  & MIP & makespan & \begin{tabular}[t]{@{}l@{}}Cut generation with\\ bound improvement\end{tabular} \\
\cite{boccia2021column} & D &  &  & ILP & makespan & Column-and-row generation \\
\cite{dell2022exact} & D &  &  & IP & makespan & B\&C \\
\cite{freitas2023exact} & D &  &  & MIP & makespan & Hybrid heuristic \\
\cite{boccia2023new} & D &  &  & MILP & makespan & B\&C \\
\cite{mahmoudinazlou2024hybrid} & D &  &  &  & makespan & Hybrid genetic algorithm \\
\bottomrule
\end{tabular*}
\begin{tablenotes}\footnotesize
\item[1] In the column ``L/R'', ``D'' means the drone must be retrieved at a node different from its launch node.
\item[2] In the column ``Non-cus points'', ``\ding{51}'' means L/R operations can be performed at non-customer nodes.
\item[3] In the column ``Multi-visit'', ``\ding{51}'' means the drone is allowed to service multiple customers within one launch.
\end{tablenotes}
\end{threeparttable}
\end{table}

% optimizing synchronized truck-and-drone operations. This work laid the foundation for subsequent research by demonstrating the potential of synchronized drone and truck operations to enhance delivery efficiency.

\cite{murray2015flying} was the first to propose the FSTSP structure and developed a heuristic to solve their constructed Mixed Integer Linear Programming (MILP) model with a makespan minimization objective. This model was later addressed using simulated annealing in \cite{ponza2016optimization}. Departing from model-based approaches, \cite{ha2015heuristic} introduced two routing strategies: cluster-first \& route-second and route-first \& cluster-second to optimize delivery networks. Then, the MILP model was refined in \cite{yurek2018decomposition} and \cite{de2018randomized} and solved by dedicated heuristics. To enhance operational realism, \cite{marinelli2018route} introduced non-customer launch/retrieval points, while \cite{gonzalez2020truck} allowed the drone to serve multiple customers per sortie. \cite{wang2020cooperative} further expanded the scope by jointly optimizing makespan and cost. 

As FSTSP becomes more complicated, exact methods were developed in \cite{dell2021algorithms}, \cite{el2021parcel}, and \cite{boccia2023new}. In parallel, metaheuristics (\cite{de2020variable}) and hybrid methods (\cite{mahmoudinazlou2024hybrid}) were explored to improve scalability for larger instances. However, designing customized algorithms remains arduous and time-consuming, especially when even a small change in problem assumptions may require model and algorithm adjustments. To address this challenge, \textbf{LKH-3} offers a flexible framework that can handle diverse problem settings by adapting the penalty function.

\subsubsection{Problem description and simplification}

We benchmark the penalty-based solution framework against the Hybrid Genetic Algorithm (HGA) presented in \cite{ha2020hybrid} for FSTSP. Although the original code is not available, our primary goal is to validate the framework in this coordination mode between one truck and one drone. Therefore, we preserve their problem structure (see Table~\ref{tab:fstsp_problemcompare}) and implement a Java version based on their detailed description. Before that, we describe the structure and define the key term \textbf{sortie}. 

Taking a graph $G(V,E)$, $V=\{0, 1, \ldots, n+1\}$ represents a set of depot and customer locations, and $E$ is a set of edges connecting pairs of nodes in $V$. Parcels are delivered to a set of customers $N=\{1, \ldots, n\}$ using a truck and a drone together, with the goal of minimizing the total delivery completion time. Note that although node $0$ indicates the physical depot location, node $n+1$ is a duplicate of the depot. These two nodes correspond to the starting and returning points, respectively. Furthermore, we denote the distance and travel time from node $i$ to node $j$ by truck as $d_{ij}$ and $\tau_{ij}$, respectively. Similarly, $d'_{ij}$ and $\tau'_{ij}$ represent the distance and travel time from node $i$ to node $j$ by drone. The arrival times of the truck and drone at node $i$ are denoted as $t_i$ and $t'_i$. Thus, the objective is to minimize the makespan $\max\{t_{n+1}, t'_{n+1}\}$.
% \textit{}
% \textbf{}
% \mathbf{P}
% \mathcal{P} 
% \mathscr{P}
% \mathbb{P}
% \mathfrak{P}
% \mathrm{P}
A \textbf{sortie} is defined as a 3-tuple $\langle i,j,k \rangle$:
\begin{itemize}
    \item $i \in N_0=\{0,1,\ldots,n\}$ is the \textbf{launch node} where the truck releases the drone.
    \item $j \in N=\{1,\ldots,n\}$ is the \textbf{drone node} where the drone delivers to a customer.
    \item $k \in N_+=\{1,\ldots,n+1\}$ is the \textbf{rendezvous node} where the drone rejoins the truck to recharge and prepare for the next flight.
\end{itemize}
We require $i \neq j \neq k$. Thus, the set of all possible sorties is denoted as $\mathbb{P}$:
\begin{equation*}
\mathbb{P} = \{\langle i,j,k \rangle : i \in N_0, j \in N, k \in N_+, i \neq j \neq k\}.
\end{equation*}
Forms of the type $\langle 0,j,n+1 \rangle$, where both launch and rendezvous positions are at the depot, are excluded from $\mathbb{P}$. Note that this sortie set may become more complex when incorporating real-world operations (e.g., see Table~\ref{tab:vrpd_problem} for VRP-D).

Then, Table~\ref{tab:fstsp_problemcompare} provides the simplification that maintains the FSTSP structure in \cite{ha2020hybrid} to enable algorithm comparison.

% Once being retrieved by the truck, the drone can be re-launched at the current position or travel with the truck to be re-launched at a later customer node. The truck is not allowed to wait the drone at the launch point after launching it. Synchronization is required (truck and drone must meet for launch/retrieval).

\begin{table}[H]
    \centering
    % \tiny
    \scriptsize
    % \footnotesize  % Slightly smaller to better fit
    % \small
    % \normalsize  % (default size) 
    % \large
    \renewcommand{\arraystretch}{1.1} 
    \caption{Original vs. Preserved features in FSTSP}
    \label{tab:fstsp_problemcompare}
    \begin{tabularx}{\textwidth}{>{\centering\arraybackslash}p{2.8cm} >{\centering\arraybackslash}p{7.8cm} >{\centering\arraybackslash}X} 
        \toprule
        \textbf{Factor} & \textbf{Description in \cite{ha2020hybrid}} & \textbf{This work} \\
        \midrule
Problem name & TSP-D & FSTSP \\
Objective & min-cost or min-makespan & min-makespan \\
Decision & Truck tour and scheduling of drone's launch and retrieval & same \\
Vehicle configuration & One truck + one drone & same\\
Depot requirements & Both vehicles must start from and return to the depot & same \\
Drone operation & Drone must be launched from and retrieved by the truck & same \\
Drone carriage & Drone is carried by the truck when not involved in a delivery & same \\
Drone eligibility & Customers in subset $N_D \subset N$ can be served by drone & All can be served by drone \\
Launch & Drone can only service a single customer per flight & same \\
Customer service & Customer is visited and served once by either truck or drone & same \\
Customer demand & - & - \\
Service time & - & - \\
Rendezvous position & Vehicle cannot retrieve drone from its launch position & same \\
Rendezvous & Vehicles must wait for each other & same \\
Waiting costs & Added when vehicles wait for each other at rendezvous points ($w_T = \alpha \times \Psi_T, w_D = \beta \times \Psi_D$) & - \\
Drone endurance & $e$ is 20 or 40 minutes & - \\
Launch/Retrieval time & L time ($D_L=1$ min) and R time ($D_R=1$ min) & - \\
Drone delivery structure & $\mathbb{P} = \{\langle i,j,k \rangle : i,k \in V, j \in N_D, i \neq j \neq k, \tau_{ij}' + \tau_{jk}' \leq e\}$ & $\mathbb{P} = \{\langle i,j,k \rangle : i \in N_0, j \in N, k \in N_+, i \neq j \neq k \}$ \\
Truck travel constraint & $t_{i \rightarrow k}$ + recovery time $\leq e$ & - \\
Drone travel constraint & $\tau_{ij}'$ + $\tau_{jk}'$ + recovery time $\leq e$ & - \\
Vehicle speeds & Both set to 40 km/h & Speed ratio of D/T is 50/35 \\
Launch/Retrieval sites & Customer locations & same \\
Network structure & Different distance matrices for truck and drone & Same Euclidean distance matrices \\
Transportation costs & Different per-distance costs for truck ($C_1$) and drone ($C_2$) & - \\
Fuel consumption & Simplified to cost per distance unit & - \\
Vehicle capacity & - & - \\
Solution method & Hybrid Genetic Algorithm & LKH-3 \\
        \bottomrule
    \end{tabularx}
\end{table}

% \subsubsection{Input file generation}

% Penalty function and resulting figures can be found in \href{https://drive.google.com/drive/folders/13UUiV9mOAK3c969MtVokVwVW2V2fpC7Q}{LKHexperiments1}.
% The selected 6 points serve as candidates for launch/retrieval (L/R) operations associated with customer $i$.

% \begin{figure}[H]
%     \centering
%     \setlength{\tabcolsep}{1pt} % Space between columns
%     \begin{tabular}{cc}        \includegraphics[width=0.6\textwidth]{image/fstsp_exe2.pdf} &        \includegraphics[width=0.3\textwidth]{image/fstsp_exe3.pdf} \\
%         \multicolumn{1}{c}{(a) 2.} & \multicolumn{1}{c}{(b) 3.} 
%     \end{tabular}
%     \caption{FSTSP input files.}
%     \label{fig:fstsp_exe2and3}
% \end{figure}

\subsection{TSP-mD} \label{sec:tspmd}

The FSTSP is extended to the Traveling Salesman Problem with Multiple Drones (TSP-mD), where multiple drones collaborate with a truck for deliveries.

\subsubsection{Related work}

\begin{table}[H]
\centering
\caption{Papers related to TSP-mD}
\label{tab:tspmd_papers}
\setlength{\tabcolsep}{4pt} % Adjust spacing between columns
\begin{threeparttable}
\small
\begin{tabular*}{\textwidth}{@{\extracolsep{\fill}}l@{\;}c@{\;}c@{\;}c@{\;}c@{\;}c@{\;}l@{}}
\toprule
\textbf{Authors} & \textbf{L/R} & \textbf{Non-cus} & \textbf{Multi-} & \textbf{Model} & \textbf{Objective} & \textbf{Solution Method} \\
 &  & \textbf{points} & \textbf{visit} &  &  & \\
\midrule
\cite{mourelo2016optimization} & D &  &  &  & makespan & \begin{tabular}[t]{@{}l@{}}Cluster\&Routing\\(GA, K-means)\end{tabular} \\
\cite{chang2018optimal} & S & \ding{51} &  &  & time & Cluster\&Routing \\
\cite{tu2018traveling} & D &  &  &  & cost & ALNS \\
\cite{yoon2018traveling} & D &  &  & MILP & cost & - \\
\cite{murray2020multiple} & D &  &  & MILP & makespan & Heuristic \\
\multicolumn{7}{@{}l}{(heterogeneous drones)} \\
\cite{raj2020multiple} & D &  &  &  & makespan & Heuristic \\
\multicolumn{7}{@{}l}{(variable drone speeds)} \\
\cite{salama2020joint} & S & \ding{51} &  & IP & cost & Cluster\&Routing \\
\cite{moshref2020truck} & D &  &  & MILP & waiting time & ALNS \\
\cite{moshref2020design} & S &  &  & MILP & time & \begin{tabular}[t]{@{}l@{}}Adaptive tabu search-\\simulated annealing\end{tabular} \\
\cite{gu2020vehicle} & S & \ding{51} &  & MIP & time & \begin{tabular}[t]{@{}l@{}}Location stops\&\\Allocation routing\end{tabular} \\
\cite{poikonen2020multi} & D & \ding{51} & \ding{51} & ILP & makespan & Heuristic \\
\cite{dell2021modeling} & D &  &  & MILP & makespan & B\&C \\
\cite{moshref2021comparative} & S/D &  &  & MIP & time & \begin{tabular}[t]{@{}l@{}}Truck and Drone\\Routing Algorithm\end{tabular} \\
\cite{luo2021hybrid} & D TW &  &  & MIP & cost \& satisfaction & \begin{tabular}[t]{@{}l@{}}Hybrid Optimization and\\Monte Carlo (HOMA)\end{tabular} \\
\cite{luo2021multi} & D &  & \ding{51} & MILP & makespan & Multi-start tabu search \\
\cite{kang2021exact} & S &  &  & MIP & time & logic-based Benders \\
\cite{cavani2021exact} & D &  &  & MILP & makespan & B\&C \\
\cite{wu2022coordinated} & S/D &  & \ding{51} & MINLP & makespan & ALNS \\
\cite{leon2022multi} & D &  & \ding{51} &  & makespan & agent-based approach \\
\cite{salama2022collaborative} & D & \ding{51} &  & MILP & makespan & two-phase search \\
\cite{bruni2022logic} & S &  &  & MILP & time & logic-based Benders \\
\cite{xu2023gv} & D &  & \ding{51} & MILP & makespan & Hybrid metahuristic \\
\cite{tinicc2023exact} & S/D &  &  & MILP & cost & B\&C \\
\cite{morandi2023traveling} & S/D &  & \ding{51} &  & makespan & Gurobi \\
\cite{karakose2024new} & D &  &  &  & makespan & Genetic Algorithm \\
\cite{liu2025adaptive} & D &  & \ding{51} & MILP & makespan & ALNS \\
\cite{rave2025two} & S/D &  &  & MILP & makespan & CPLEX \\

\cite{campbell2017strategic} &  &  &  & Continuous &  &  \\
\bottomrule
\end{tabular*}
\begin{tablenotes}\footnotesize
\item[1] In the column ``L/R", ``D" means drones must be retrieved at nodes different from their launch nodes, ``S" means same, ``S/D" means both options are considered, and ``TW" refers to time windows.
\item[2] In the column ``Non-cus points'', ``\ding{51}'' means L/R operations can be performed at non-customer nodes.
\item[3] In the column ``Multi-visit'', ``\ding{51}'' means drones are allowed to service multiple customers within one launch.
\end{tablenotes}
\end{threeparttable}
\end{table}

TSP-mD was originally designed using cluster-and-routing strategies (\cite{mourelo2016optimization}, \cite{chang2018optimal}, \cite{salama2020joint}). Then, mathematical models were formulated and exact solution methods were developed to optimally solve them, such as branch-and-cut approaches in \cite{cavani2021exact} and \cite{tinicc2023exact}, and logic-based Benders decomposition in \cite{kang2021exact} and \cite{bruni2022logic}. Concurrently, tailored heuristics were proposed for larger instances, including Adaptive Large Neighborhood Search (ALNS) by \cite{tu2018traveling}, Adaptive Tabu Search by \cite{moshref2020truck} and \cite{moshref2020design}, and a heuristic-based algorithm by \cite{poikonen2020multi}. Notably, \cite{campbell2017strategic} adopted a continuous approximation method, offering an analytical perspective distinct from the optimization-based approaches above.

Furthermore, TSP-mD networks do not strictly prohibit \textbf{loop} operations. As shown in the ``L/R" column of Table~\ref{tab:tspmd_papers}, some studies labeled ``D" follow the no-wait principle as in \cite{murray2015flying}, while some labeled ``S" enforce the truck to wait at launch points for retrieval (\cite{chang2018optimal}, \cite{gu2020vehicle}, \cite{kang2021exact}). More recently, hybrid configurations (``S/D") have emerged in \cite{moshref2021comparative} and \cite{tinicc2023exact} to enhance operational flexibility. Further extensions include using non-customer points for L/R operations (\cite{chang2018optimal}, \cite{salama2020joint}, \cite{poikonen2020multi}), incorporating time windows and customer satisfaction metrics (\cite{luo2021hybrid}), servicing multi-customer per sortie (\cite{luo2021multi}, \cite{leon2022multi}), and heterogeneous drones (\cite{murray2020multiple}) that relaxes the common assumption of homogeneity. These developments increase structural and algorithmic complexity, reinforcing the need for a unified framework capable of accommodating diverse problem settings.

\subsubsection{Problem description and simplification}
To verify the proposed TSP-mD penalty function in LKH-3, we compare its performance with the Adaptive Large Neighborhood Search (ALNS) developed by \cite{tu2018traveling}, which adopts limited operational constraints. As the original code is unavailable and we aim to validate the framework through the cooperative mode, we preserve their problem structure (see Table~\ref{tab:tspmd_problemcompare}) and develop a Java implementation based on their description. Before detailing this preservation, we describe the TSP-mD problem as follows. 

Beyond FSTSP, TSP-mD employs multiple drones ($m>1$) alongside a single truck to enhance delivery efficiency. The truck initiates its tour with all drones at the depot and returns after serving all customers. Note that the truck can launch any number of drones from a location and carry the remaining drones while traveling between nodes. Multiple drones can fly at the same time for delivery tasks. Consequently, this parallel operation is constrained by the available drone count, significantly increasing the scheduling complexity. 

The simplification applied to retain the TSP-mD characteristics in \cite{tu2018traveling} for comparison is listed in Table~\ref{tab:tspmd_problemcompare}.

% Coordination of multiple drone flights. Potential waiting times for the truck when multiple drones need to be retrieved at the same customer location. The truck carries drones when they are not engaged in delivery tasks. The objective is to minimize the total completion time (makespan) of the entire delivery process.

% Based on the ALNS algorithm detailed in \cite{tu2018traveling}, we develop a JAVA implementation. Then, comparative analysis is conducted between this implementation and the penalty function used in the TSP-mD framework.

\begin{table}[H]
    \centering
    % \tiny
    \scriptsize
    % \footnotesize  % Slightly smaller to better fit
    % \small
    % \normalsize  % (default size) 
    % \large
    \renewcommand{\arraystretch}{1.1} 
    \caption{Original vs. Preserved features in TSP-mD}
    \label{tab:tspmd_problemcompare}
    \begin{tabularx}{\textwidth}{>{\centering\arraybackslash}p{2.8cm} >{\centering\arraybackslash}p{7.8cm} >{\centering\arraybackslash}X} 
        \toprule
        \textbf{Factor} & \textbf{Description in \cite{tu2018traveling}} & \textbf{This work} \\
        \midrule
Problem name & TSP-mD & same \\
Objective & min-cost & min-makespan \\
Decision & Truck tour, \#drone to be used, L\&R scheduling of drones & same \\
Vehicle configuration & One truck + m same drones (m $>$ 1) & same \\
Depot requirements & Vehicles must start from and return to the depot & same \\
Drone operation & Drone must be launched from and retrieved by the truck & same \\
Drone carriage & Drones are carried by the truck when not involved in a delivery & same \\
Drone number & No more than m drones flying simultaneously & same \\
Drone eligibility & Customers in subset $N_D \subset N$ can be served by drone & All can be served by drone \\
Launch & Drone can only service a single customer each time & same \\
Customer service & Customer is visited and serviced once by either truck or drone & same \\
Customer demand & - & - \\
Service time & - & - \\
Rendezvous position & Vehicle cannot retrieve a drone from its launch position & same \\
Rendezvous & Multiple drones can rejoin at same or different nodes in the truck tour & same \\
% Rendezvous position & Position of launch node must be previous to the rendezvous node in the truck tour & same \\
Waiting time & At $k$,
Vehicles cannot wait for each other more than $\Delta$ unit of time & - \\
Drone endurance & $e$ is 20 minutes & - \\
Launch/Retrieval time & - & - \\
Drone delivery structure & $\mathbb{P} = \{\langle i,j,k \rangle : i,k \in V, j \in N_D, i \neq j \neq k, \tau_{ij}' + \tau_{jk}' \leq e,\left| t_{i\rightarrow k} - (\tau_{ij}' + \tau_{jk}') \right| \leq \Delta \}$ & $\mathbb{P} = \{\langle i,j,k \rangle : i \in N_0, j \in N, k \in N_+, i \neq j \neq k \}$ \\
Vehicle speeds & Both set to 40 km/h & Speed ratio of D/T is 50/35 \\
Launch/Retrieval sites & Customer locations & same \\
Network structure & Manhattan $d_{ij}$ for truck and Euclidean $d_{ij}'$ for drone & Same Euclidean distance matrices \\
Transportation costs & Different per-distance costs for truck ($C_1$) and drone ($C_2$) & - \\
Vehicle capacity & - & - \\
Solution method & Adaptive Large Neighborhood Search (ALNS) metaheuristic & LKH-3 \\
        \bottomrule
    \end{tabularx}
\end{table}

% \subsubsection{Input file generation}

% Penalty function and resulting figures can be found in \href{https://drive.google.com/drive/folders/13UUiV9mOAK3c969MtVokVwVW2V2fpC7Q}{LKHexperiments1}.

\subsection{VRP-D} \label{sec:vrpd}
After introducing multiple drones in TSP-mD, the Vehicle Routing Problem with Drone (VRP-D) extends the problem to multiple identical trucks, each equipped with one or more drones. Every truck begins and ends its route at the depot and coordinates with its designated drone(s) to serve customers along the route. 

This variant integrates the classical Vehicle Routing Problem (VRP) with drone operations and primarily involves two key decisions: (1) assigning customers to truck fleet, and (2) optimizing each truck's route and its coordination with drone(s). Moreover, allowing drones to operate across different trucks adds complexity \cite{kitjacharoenchai2019multiple}.

\subsubsection{Related work}

% Third literature review table with categories
\begin{table}[H]
\centering
\caption{Papers related to VRP-D}
\label{tab:vrpd_papers}
\setlength{\tabcolsep}{4pt} % Adjust spacing between columns
\begin{threeparttable}
\small
\begin{tabular*}{\textwidth}{@{\extracolsep{\fill}}l@{\;}c@{\;}c@{\;}c@{\;}c@{\;}c@{\;}l@{}}
\toprule
\textbf{Authors} & \textbf{L/R} & \textbf{Non-cus} & \textbf{Multi-} & \textbf{Model} & \textbf{obj} & \textbf{Solution Method} \\
 &  & \textbf{points} & \textbf{visit} &  &  & \\
\midrule
\multicolumn{7}{c}{Each truck carries one drone from the beginning. 1-1} \\
\cite{sacramento2019adaptive} & D &  &  & MIP & -cost & ALNS \\
\cite{chiang2019impact} & D &  &  & MIP & -cost & Genetic heuristic \\
\cite{liu2020two} & D &  & \ding{51} & MIP & -cost & Heuristic \\
\cite{euchi2021hybrid} & D &  &  & MILP & -time & Genetic algorithm \\
\cite{kuo2022vehicle} & D TW &  &  & MILP & -cost & VNS \\
\cite{gu2022hierarchical} & D &  & \ding{51} & MILP & -cost & ILS-VND \\
\cite{lei2022dynamical} & D &  &  & MILP & -cost & Dynamical artificial bee \\
\cite{zhen2023branch} & D &  &  & MIP & -cost & B\&P\&C \\
\cite{yin2023branch} & D TW &  & \ding{51} & MILP & -cost & B\&P\&C \\
\cite{kuo2023applying} & D TW &  &  & MIP & 
\begin{tabular}[t]{@{}l@{}}makespan \&\\ carbon emission \end{tabular} & NSGA-II \\
\cite{liu2024cooperated} & D TW &  & \ding{51} & MILP & -cost & Heuristic \\
\begin{tabular}[t]{@{}l@{}}\cite{peng2025multi}\\ (time-dependent truck speed) \end{tabular} & D &  &  & MIP & 
\begin{tabular}[t]{@{}l@{}}cost \&\\ makespan \end{tabular} & ALNS \\
\midrule
\multicolumn{7}{c}{Each truck carries multiple drones from the beginning. 1-m} \\
\cite{wang2017vehicle} & D &  &  &  & makespan & Worst-case analysis \\
\cite{poikonen2017vehicle} & S/D &  &  &  & makespan & Extend analysis \\
\cite{di2017last} & D TW &  & \ding{51} & MIP & -cost & CPLEX \\
\cite{schermer2018algorithms} & D &  &  &  & makespan & Two-phase heuristic \\
\cite{schermer2019hybrid} & D & \ding{51} &  & MILP & makespan & Heuristic \\
\cite{schermer2019matheuristic} & S/D &  &  & MILP & makespan & Heuristic \\
\cite{kitjacharoenchai2019vehicle} & D &  &  & MILP & makespan & CPLEX \\
\cite{poikonen2020multi} & S/D &  &  &  & makespan &  \\
\cite{kitjacharoenchai2020two} & D &  & \ding{51} & MIP & makespan & Heuristic \\
\cite{das2020synchronized} & D TW &  &  & MIP & -cost & Heuristic \\
\cite{tamke2021branch} & D &  &  & MILP & makespan & B\&C \\
\cite{chen2021adaptive} & S &  &  & MILP & -time & ALNS \\
\cite{tamke2023vehicle} & D &  &  & MILP & -cost & Gurobi \\
\multicolumn{7}{@{}l}{(variable drone speeds)} \\
\cite{gao2023scheduling} & D TW & \ding{51} &  &  & -cost & C\&G \\
\cite{han2023vehicle} & D TW &  & \ding{51} &  & -cost & chimp optimization \\
\cite{yang2024optimization} (multiple depots) & D & \ding{51} &  & MILP & -time & 3-phase method \\
\cite{schmidt2025exact} & D &  &  & MIP & -cost & B\&P\&C \\

\bottomrule
\end{tabular*}
\begin{tablenotes}\footnotesize
\item[1] In the column ``L/R", ``D" means drones must be retrieved at nodes different from their launch nodes, ``S" means same, ``S/D" means both options are considered, and ``TW" refers to time windows.
\item[2] In the column ``Non-cus points'', ``\ding{51}'' means L/R operations can be performed at non-customer nodes.
\item[3] In the column ``Multi-visit'', ``\ding{51}'' means drones are allowed to service multiple customers within one launch.
\end{tablenotes}
\end{threeparttable}
\end{table}

% The Vehicle Routing Problem with Drones (VRP-D) moves beyond single-truck models by utilizing multiple trucks. This framework necessitates synchronization between trucks and drones. 

As shown in Table~\ref{tab:vrpd_papers}, VRP-D is classified into two categories based on the number of drones assigned to each truck. The first category, 1-1 VRP-D, equips each truck with a single drone and is relatively straightforward. The second category, 1-m VRP-D, allows multiple drones per truck. In particular, when drones must return to their designated trucks, each truck's route resembles an \textbf{FSTSP} structure in the 1-1 case and a \textbf{TSP-mD} structure in the 1-m case.

In the 1-1 VRP-D setting, most studies adopted the ``D" drone retrieval paradigm, where drones must be recollected at locations different from their launch points. Researchers formulated mathematical models and developed various solution approaches, including exact methods such as branch-price-and-cut (\cite{zhen2023branch}), heuristics (\cite{sacramento2019adaptive}, \cite{chiang2019impact}, \cite{euchi2021hybrid}), and metaheuristics (\cite{lei2022dynamical}). Additionally, \cite{kuo2022vehicle} and \cite{liu2024cooperated} incorporated time window constraints. 

For the more complex 1-m VRP-D setting, early efforts examined theoretical performance through worst-case analysis (\cite{wang2017vehicle}) and extended analysis (\cite{poikonen2017vehicle}). They prioritized establishing theoretical bounds over solving specific instances. Subsequent studies introduced mathematical models, proposed exact solution methods (\cite{tamke2021branch}, \cite{gao2023scheduling}, \cite{schmidt2025exact}) and heuristics (\cite{schermer2019matheuristic}, \cite{kitjacharoenchai2020two}, \cite{chen2021adaptive}). Time windows were further integrated into the models by \cite{dickey2019}, \cite{das2020synchronized} and \cite{gao2023scheduling}, while features such as non-customer L/R sites and multiple-customer per sortie were introduced by \cite{gu2022hierarchical}, \cite{schermer2019hybrid} and \cite{gao2023scheduling}, expanding operational flexibility. 

Furthermore, when drones are not dedicated to specific trucks, as explored in \cite{ibrocska2023multiple}, the problem becomes significantly more complicated. Therefore, to address this growing complexity in both problem formulation and algorithm design, we propose an efficient solution framework by developing adjustable penalty functions within \textbf{LKH-3}.

\subsubsection{Problem description}

% To facilitate the following comparison, we assume that each truck is initially assigned one drone (i.e 1-1 VRP-D), and each drone must return to its original truck. Building on these factors, the performance of the penalty function for VRP-D in \textbf{LKH} is compared to the ALNS method presented in \cite{sacramento2019adaptive}, using the source code generously provided by the authors.

To validate the penalty function developed in LKH-3, we focus on the 1-1 VRP-D configuration with minimal additional constraints, as studied in \cite{sacramento2019adaptive}, and benchmark our performance against their algorithm, the Adaptive Large Neighborhood Search (ALNS). The more complex 1-m VRP-D scenario will be discussed in Section 6.1 \hyperref[sec:extensions]{Extensions}.

Since the source code from \cite{sacramento2019adaptive} is available, we can use their complete formulation for comparison without the problem simplifications required in Section 4.1.2 and 4.2.2 for FSTSP and TSP-mD, respectively. Thus, Table~\ref{tab:vrpd_problem} summarizes the characteristics of the VRP-D problem addressed in \cite{sacramento2019adaptive}. The objective is to minimize the total operational cost of the truck fleet. In particular, the problem incorporates constraints 1-4 described in \hyperref[sec:introduction]{Introduction}: drone flight range (endurance $e$), drone operation time (launch $D_L$ and recovery $D_R$), customer service time ($Serv_D$ by drone and $Serv_T$ by truck), and drone payload capacity (subset of drone-eligible customers $N_D$), with the last one restricting high-demand customers to truck-only service. Hence, the basic 3-tuple sortie set described in Section 4.1.2 is adapted as follows:
\begin{equation*}
\mathbb{P} = \{\langle i,j,k \rangle : i \in N_0, j \in N_D, k \in N_+, i \neq j \neq k, D_L + D_R + Serv_D + \tau_{ij}' + \tau_{jk}' \leq e\}.
\end{equation*}

% Each sortie cannot exceed the drone's maximum flight range. Drone Launch and Retrieval Time by the truck. Customer service time by the truck or a drone. Customers are classified based on package weight, where those with larger demands must be served by a truck. and utilization. The objective is to minimize the total operational cost of the fleet of one truck-one drone combinations.

% In particular, after accounting for the subset of drone-eligible customers $N_D$, flight endurance $e$, drone launch time $D_L$, drone recovery time $D_R$, and drone service time $Serv_D$, 

%
\begin{table}[H]
    \centering
    % \tiny
    \scriptsize
    % \footnotesize  % Slightly smaller to better fit
    % \small
    % \normalsize  % (default size) 
    % \large
    \renewcommand{\arraystretch}{1.1} 
    \caption{1-1 VRP-D features}
    \label{tab:vrpd_problem}
    \begin{tabularx}{\textwidth}{>{\centering\arraybackslash}p{3.0cm} >{\centering\arraybackslash}X } 
        \toprule
        \textbf{Factor} & \textbf{Description in \cite{sacramento2019adaptive}} \\
        \midrule
Problem name & Vehicle Routing Problem with Drones (VRP-D) \\
Objective & Minimize overall operational cost of the fleet of vehicles \\
Decision & Routes of trucks and scheduling of drones' launch and retrieval \\
Vehicle configuration & Multiple homogeneous trucks, each coordinated with a single drone \\
Depot requirements & All vehicles must start from and return to the single depot \\
Drone operation & Drone must be launched from and retrieved by the same truck \\
Drone carriage & Drone is carried by its assigned truck when not in a delivery \\
Drone eligibility & Subset $N_D \subset N$ of customers can be served by drone \\
Drone payload & Maximum drone payload capacity of $Q_D$ \\
Launch & Drone can only serve one customer each time due to payload capacity \\
Customer service & Each customer is visited exactly once, either by a truck or by a drone \\
Customer demand & 86\% of customers can be served by the drone given their lower demands \\
Service time & Service a customer by truck ($Serv_T$) or drone ($Serv_D$) \\
Rendezvous position & Vehicle cannot retrieve its assigned drone from its launch position \\
Rendezvous & Vehicle and its drone must wait for each other \\
Drone endurance & Parameter $e$ for maximum flight endurance \\
Launch/Retrieval time & Launch time ($D_L$) and Retrieval time ($D_R$) \\
Delivery structure & $\mathbb{P} = \{\langle i,j,k \rangle : i \in N_0, j \in N_D, k \in N_+, i \neq j \neq k, D_L + D_R + Serv_D + \tau_{ij}' + \tau_{jk}' \leq e\}$ \\
Vehicle speeds & Truck speed ($v_T$) and Drone speed ($v_D$) (ratio of D/T is 50/35) \\
Network structure & Same Euclidean distance matrices for both truck and drone travel \\
Transportation costs & Different costs for truck ($c_{ij}^T$) and drone ($c_{ij}^D = \alpha \cdot c_{ij}^T$ where $\alpha$ is cost factor) \\
Time constraints & Maximum duration $T_{max}$ of each truck's route \\
Vehicle capacity & Maximum truck capacity of $Q_T$ \\
% Multiple drones constraint & Each truck has exactly one drone, new launches cannot occur while a drone is already flying \\
Synchronization & Time synchronization at launch and recovery points \\
Solution method & Adaptive Large Neighborhood Search (ALNS) metaheuristic \\
        \bottomrule
    \end{tabularx}
\end{table}

% \subsubsection{Input file generation}

% Penalty function and resulting figures can be found in \href{https://drive.google.com/drive/folders/13UUiV9mOAK3c969MtVokVwVW2V2fpC7Q}{LKHexperiments1}.

\section{Experimental results} \label{sec:experiments}
% Thus, the experimental results of these three comparisons are presented in this section.

As illustrated in the \hyperref[sec:methodology]{Methodology} and \hyperref[sec:variants]{Problem variants} sections, we verify the flexibility of our solution framework by benchmarking it against HGA (\cite{ha2020hybrid}) on FSTSP, ALNS (\cite{tu2018traveling}) on TSP-mD, and ALNS (\cite{sacramento2019adaptive}) on VRP-D. 

% We conduct all experiments on Intel Xeon E5-2640 v3 processors at 2.60 GHz with 59 GB memory per node

\begin{table}[H]
    \centering
    % \tiny
    \scriptsize
    % \footnotesize  % Slightly smaller to better fit
    % \small
    % \normalsize  % (default size) 
    % \large
    \caption{Comparison experiments}
    \label{tab:comparison}
    \begin{tabular*}{\textwidth}{@{\extracolsep{\fill}}l@{\;}c@{\;}c@{\;}c@{\;}c@{\;}c@{\;}c@{\;}c@{\;}c@{\;}c@{\;}c@{}}
    %%
    % \begin{tabularx}{\textwidth}{l r r r r r r r r r r}
    % \begin{tabularx}{\textwidth}{l c c c c c c c c c c}
    % \begin{tabularx}{\textwidth}{l *{10}{>{\centering\arraybackslash}X}}
        \toprule
        
        \multicolumn{1}{c}{} & \multicolumn{5}{c}{\textbf{Indirect}} & \multicolumn{5}{c}{\textbf{Direct}} \\
        
        \cmidrule(lr){2-6} \cmidrule(lr){7-11}
        \raisebox{1.5ex}[0pt][0pt]{Variant} &  Paper & Algo & Ins & \#Cus & Results & Paper & Algo & Ins & \#Cus & Results \\
        
        \midrule
FSTSP(min-time) & \cite{ha2020hybrid} & HGA & Set 4 & 6/9/19/29/39 & Tables~\ref{tab:fstsp_1LKH_small}--\ref{tab:fstsp_1LKH_longer} & \cite{boccia2023new} & B\&C & Set 1 & 10/20 & Tables~\ref{tab:Boccia10cusE20_1LKH}--\ref{tab:Boccia20cusE40_1LKH} \\
TSP-mD(min-time) & \cite{tu2018traveling} & ALNS & Set 4 & 6/9/19/29/39 & Tables~\ref{tab:tspmd_1LKH_small}--\ref{tab:tspmd_1LKH_longer} & - & - & - & - & - \\
VRP-D(min-cost) & - & - & - & - & - & \cite{sacramento2019adaptive} & ALNS & Set 6 & 6/10/12/20 & Tables~\ref{tab:vrpd_1LKH}--\ref{tab:vrpd_1LKH_long} \\
        \bottomrule
    \end{tabular*}
\end{table}

In this section, for FSTSP and TSP-mD, we implement HGA and ALNS in Java due to the lack of source code, denoted as ``\textbf{Indirect}" in Table~\ref{tab:comparison}. Each algorithm is run five times with a 5min time limit per run, and the best solutions found are presented. For VRP-D, we directly compare against the best solutions in 10 independent 5min runs reported in \cite{sacramento2019adaptive}, indicated as ``\textbf{Direct}" in Table~\ref{tab:comparison}. Meanwhile, based on the fixed initial seeds, LKH is executed \textbf{once} under the \textbf{5min} limit, as specified by ``TIME\_LIMIT: 300" in the input execution file \texttt{.par}. We conduct all experiments on AMD EPYC 7513 processors at 2.60 GHz and 248 GB memory per node, using 1 CPU core and 8 GB memory allocated per job, supported by the USC Center for Advanced Research Computing (CARC) Discovery cluster. 

In addition, according to the input files generation in Section 2.1, some points in \texttt{flat\_points} may share the same physical location. To address this, we assign different scores to identify locations, where points at the same site receive the same score. This enables us to distinguish whether a site has been visited by the truck. Based on this mechanism, two configurations of LKH are defined: the first (\textbf{LKH Config1}) prohibits the truck from revisiting nodes after departure, while the second (\textbf{LKH Config2}) permits such revisitation. As a result, LKH Config2 may yield a coordinated tour with a lower objective value, such as makespan or operational cost.

Overall, we conduct comparison experiments either by implementing our own version of the referenced algorithms or by directly using reported results. As summarized in Table~\ref{tab:comparison}, we first evaluate the performance of LKH Config1 for FSTSP, comparing it with the HGA developed in \cite{ha2020hybrid} using Poikonen instances (Set 4), and with best-known (BKS) exact solutions reported in \cite{boccia2023new} using Murray instances (Set 1). Next, we present the comparison of TSP-mD between LKH Config1 and the ALNS proposed by \cite{tu2018traveling} using Poikonen instances (Set 4). Finally, we compare LKH Config1 against the ALNS results reported in \cite{sacramento2019adaptive}, and further improve some of their BKSes by extending the LKH runtime to 24h. We also compare LKH Config1 and Config2 on each problem variant.

\subsection{FSTSP}
% HGA is executed five times, each with a 5-minute time limit. The best result is recorded for comparison with LKH. 

In this subsection, we first compare LKH Config1 to HGA on ``Set 4": 125 min-makespan instances of 6/9/19/29/39 customers from \citet{poikonen2019branch}. For the HGA, we present the best solution over 5 runs, each with a 5min time limit. Then, two configurations of LKH on ``Set 4" are compared. Moreover, we compare LKH Config1 with the BKS exact solutions provided in \cite{boccia2023new} on ``Set 1": 36/120 min-makespan instances of 10/20 customers from \cite{murray2015flying} with drone endurance parameter set to 20min or 40min, giving 72/240 instances. 

Before analyzing the comparison tables: Tables~\ref{tab:fstsp_1LKH_small}--\ref{tab:fstsp_1LKH_longer} between HGA and LKH Config1, Tables~\ref{tab:fstsp_2LKH_small}--\ref{tab:fstsp_2LKH_medium} between LKH Config1 and Config2, and Tables~\ref{tab:Boccia10cusE20_1LKH}--\ref{tab:Boccia20cusE40_1LKH} between Boccia and LKH Config1, we sketch several characteristic instances as follows.

\subsubsection{Instances} 
\begin{figure}[H]
    \centering
    \setlength{\tabcolsep}{1pt} % Space between columns
    \begin{tabular}{ccc}
        \includegraphics[width=0.33\textwidth]{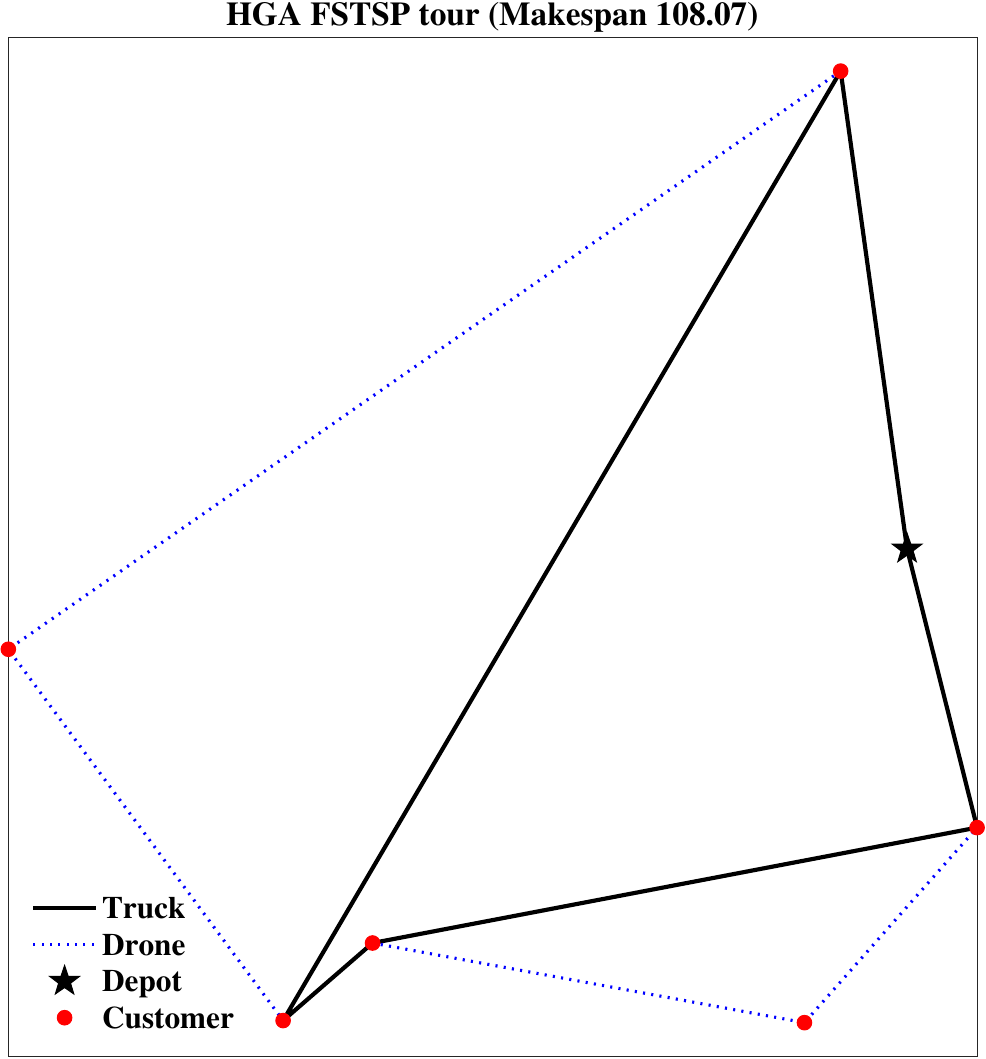} &
        \includegraphics[width=0.33\textwidth]{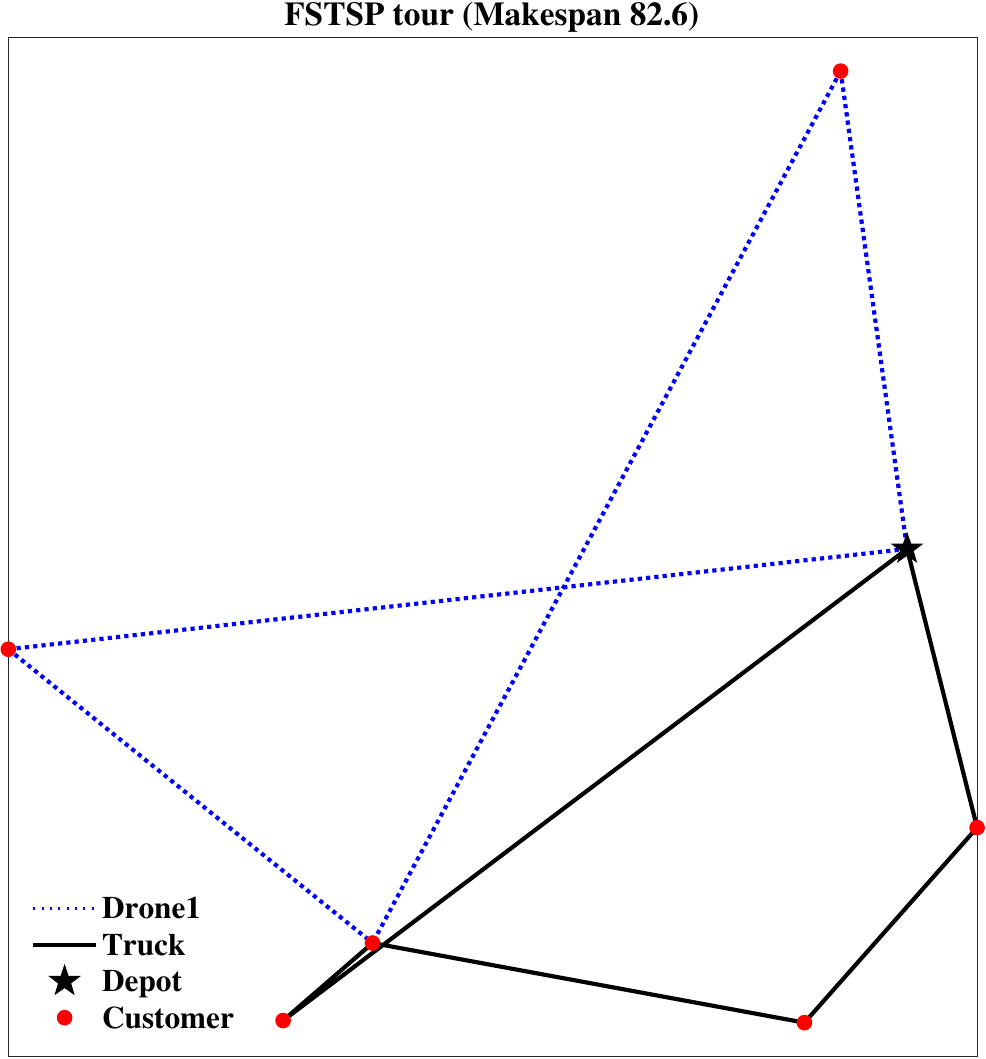} &
        \includegraphics[width=0.33\textwidth]{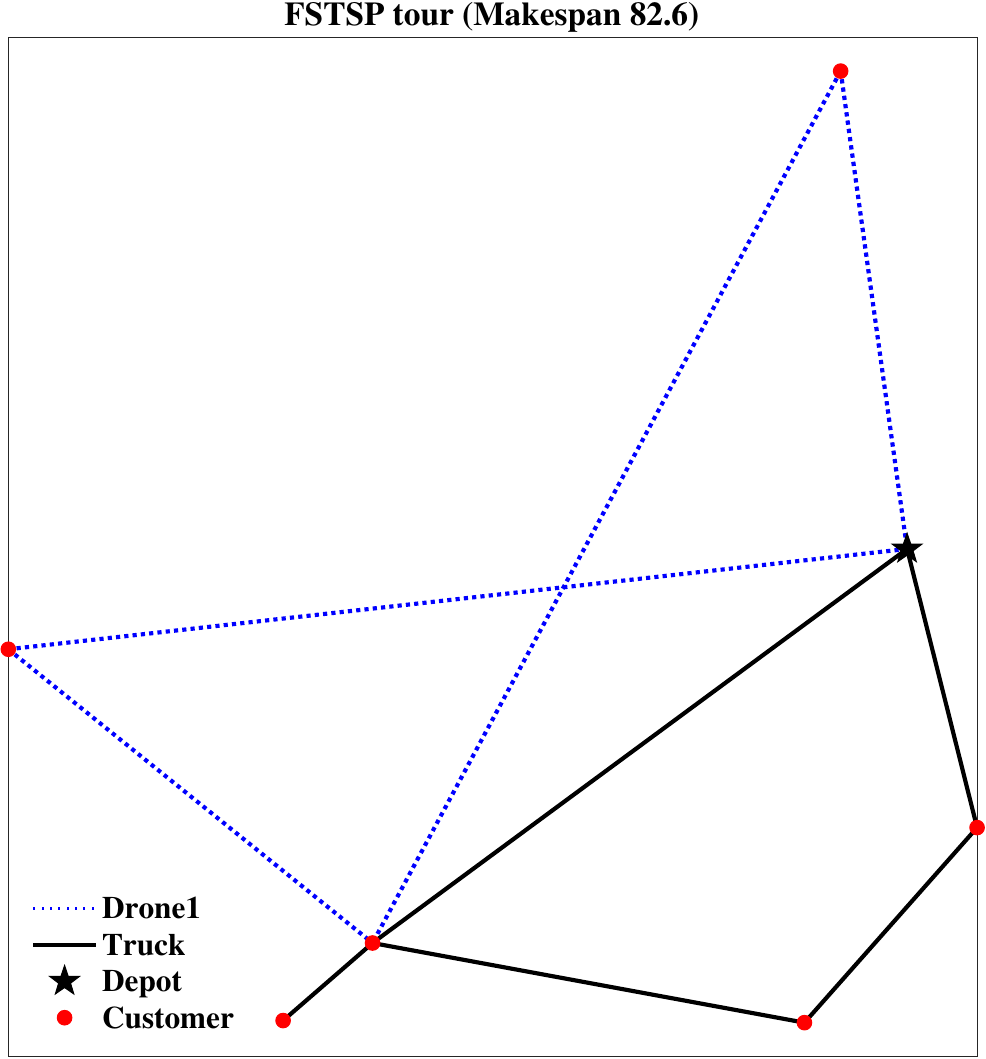} \\
        \multicolumn{1}{c}{(a) HGA.} & \multicolumn{1}{c}{(b) LKH Config1.} & \multicolumn{1}{c}{(c) LKH Config2.}
    \end{tabular}
    \caption{FSTSP tours on instance poi-7-14.}
    \label{fig:FSTSP_figcompare1}
\end{figure}

Firstly, based on the comparison between two configurations in Tables~\ref{tab:fstsp_2LKH_small}--\ref{tab:fstsp_2LKH_medium}, although LKH Config2 allows node revisits, only three instances (poi-7-14, poi-7-23, and poi-40-6) exhibit such behavior. Moreover, revisiting nodes does not necessarily yield improved objective values. For example, in instance poi-7-14 (Fig.~\ref{fig:FSTSP_figcompare1}), LKH Config2 achieves the same makespan of 82.6 despite the truck visiting a customer node in the lower region twice.

\begin{figure}[H]
    \centering
    \setlength{\tabcolsep}{1pt} % Space between columns
    \begin{tabular}{ccc}
        \includegraphics[width=0.33\textwidth]{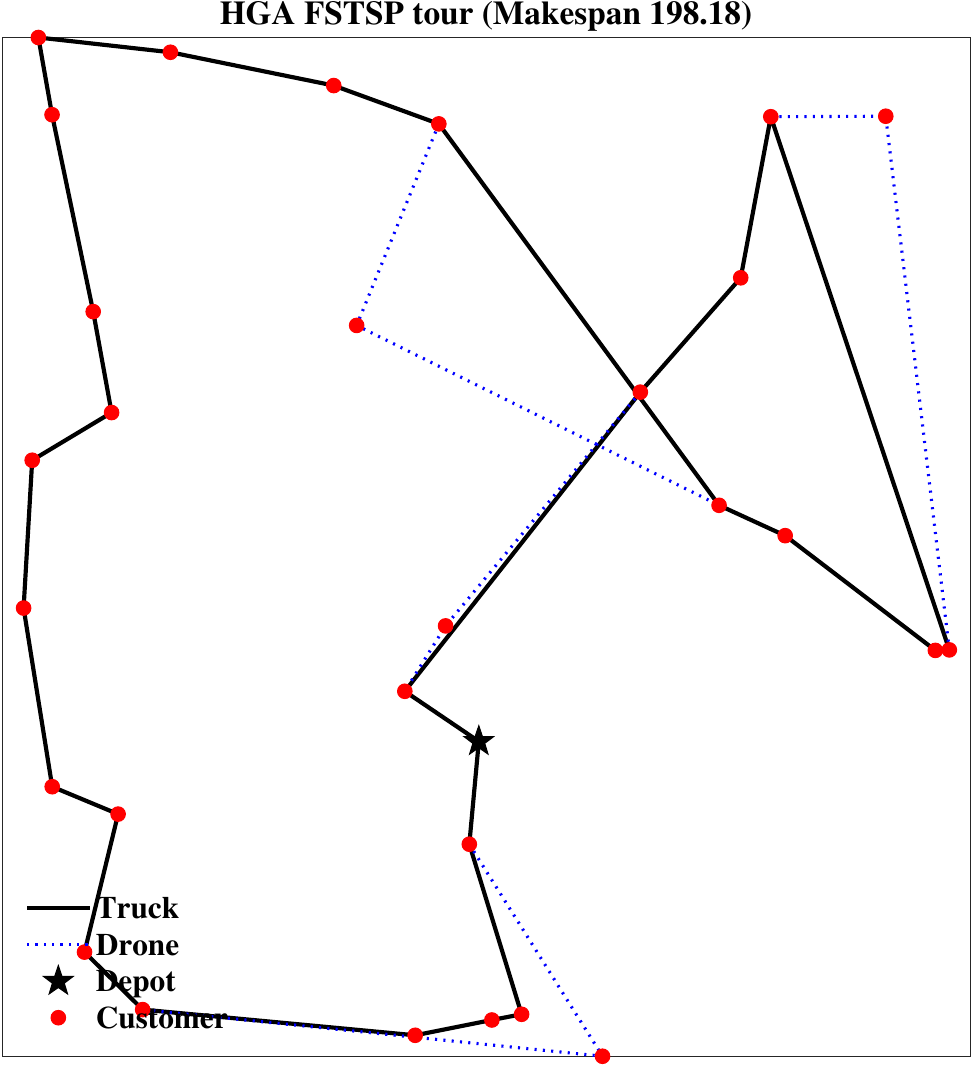} &
        \includegraphics[width=0.33\textwidth]{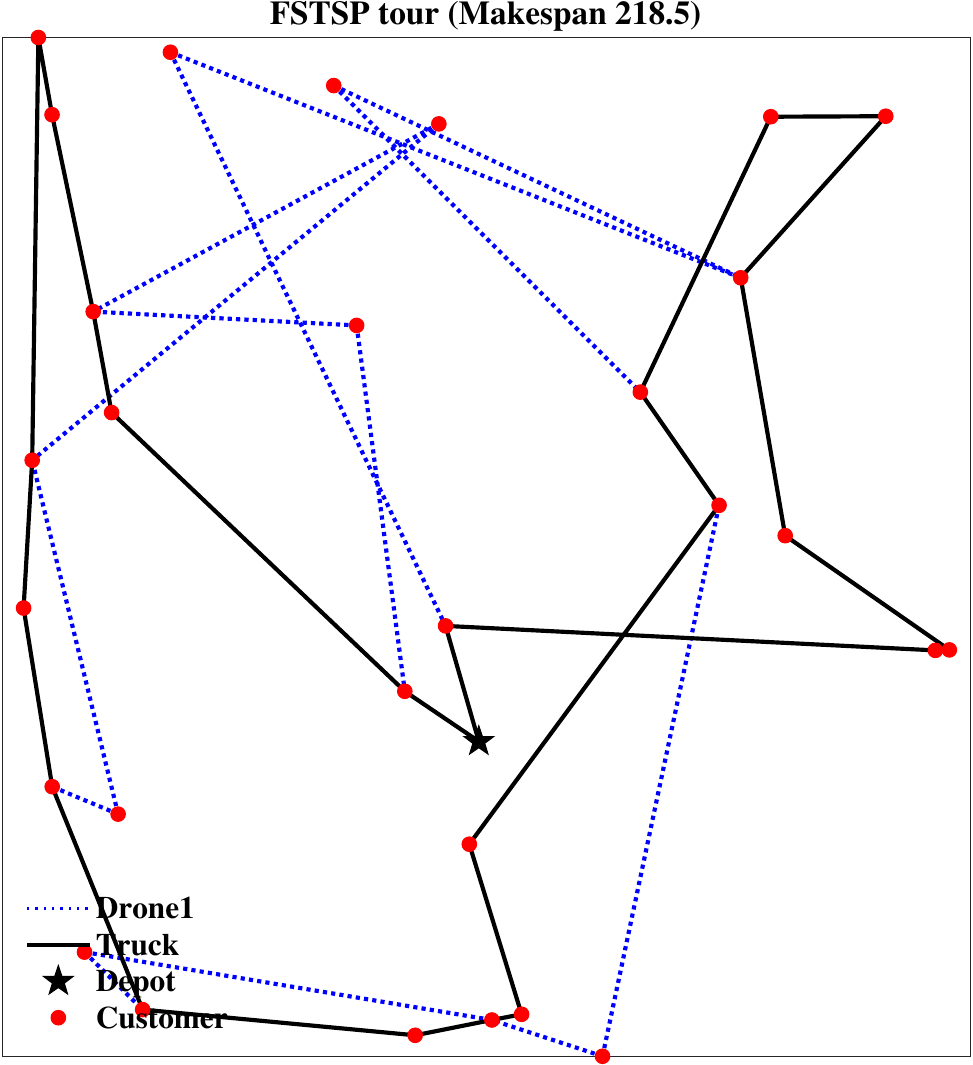} &
        \includegraphics[width=0.33\textwidth]{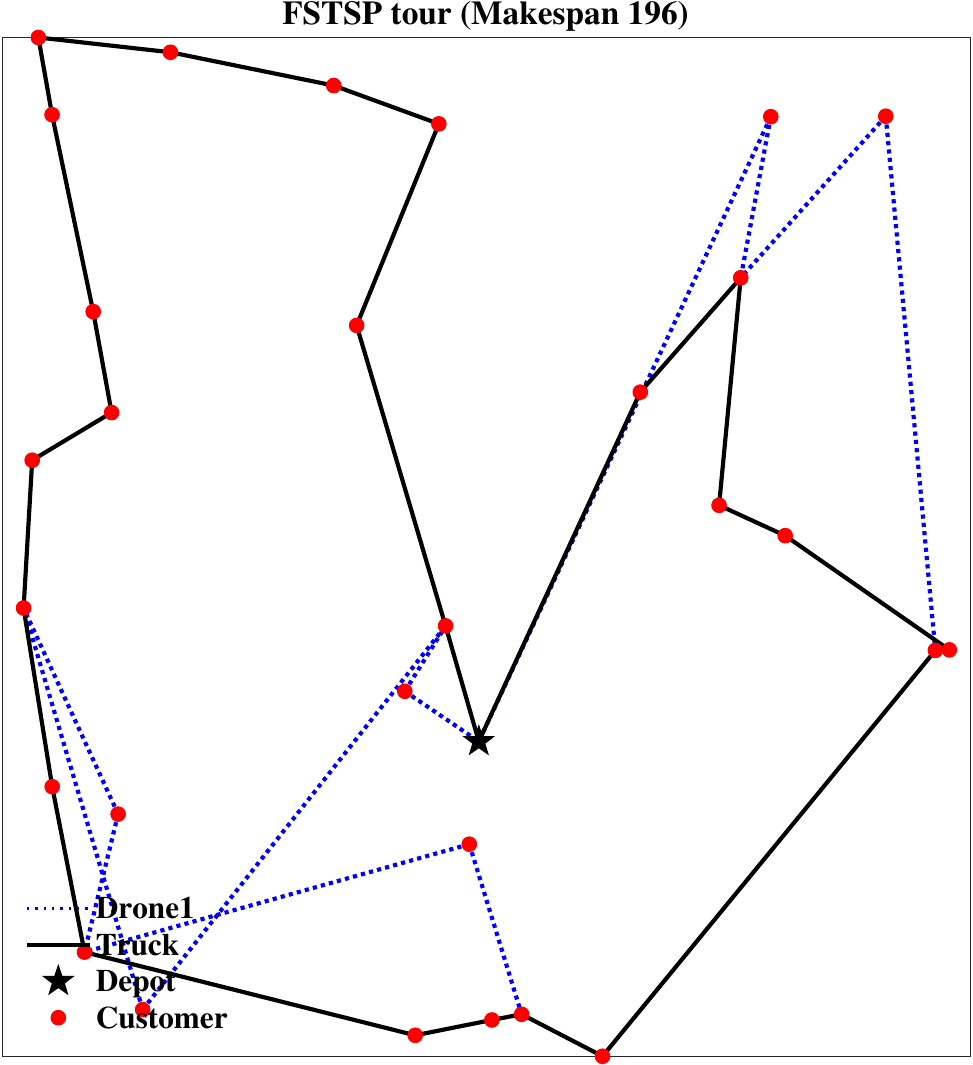} \\
        \multicolumn{1}{c}{(a) HGA.} & \multicolumn{1}{c}{(b) LKH Config1.} & \multicolumn{1}{c}{(c) LKH Config2.}
    \end{tabular}
    \caption{FSTSP tours on instance poi-30-14.}
    \label{fig:FSTSP_figcompare2}
\end{figure}

Secondly, according to the comparison between HGA and LKH Config1 (run once within 5min) in Tables~\ref{tab:fstsp_1LKH_small}--\ref{tab:fstsp_1LKH_medium}, LKH does not consistently outperform HGA. For example, as shown in Fig.~\ref{fig:FSTSP_figcompare2}, HGA achieves a lower objective value of 198.18 on instance poi-30-14. However, it should be emphasized that this result represents the best solution found from five independent runs, each using 5min.

\subsubsection{Computational performance analysis}

\Cref{tab:fstsp_1LKH_small} presents a comprehensive comparison among the exact optimum derived from the Mixed Integer Programming (MIP) model solved via the CPLEX solver, HGA and LKH Config1, on small FSTSP instances. For example, instance poi-7-1 corresponds to the first location configuration among the 25 poi-7- instances, each consisting of one depot and $7-1=6$ customers. Notably, HGA fails to generate valid tours for instances \textbf{poi-7-23} and \textbf{poi-10-23} across all five runs. 

In Table~\ref{tab:fstsp_1LKH_small}, the column $N_D$ indicates the number of drone-eligible customers in each instance, which equals the total number of customers since the drone payload constraint is not considered in this variant. For CPLEX, $z^*$ displays the optimal makespan and Dcus specifies the number of drone-served customers. Note that the required computation time ($t^{\textit{MIP}}$) increases considerably for poi-10- instances due to the larger problem size. For HGA, Table~\ref{tab:fstsp_1LKH_small} reports the best makespan ($M^{\textit{HGA}}$) among five 5min executions, along with the corresponding computational time in seconds. In addition to these items, LKH Config1 records the iteration (Iter) at which the best tour is identified in a 5min execution. 

Furthermore, the column labeled $ \%gap^* = 100 \cdot (z^*-\text{M}^1)/z^* $, quantifies the $\%gap$ of LKH Config1 over CPLEX. As indicated, LKH Config1 matches the optimal makespan for all poi-7- instances, and produces suboptimal results for only 6 out of all 25 poi-10- instances, which are highlighted in bold. On average, LKH takes 6.24s to solve poi-7- instances, exceeding the 2.99s used by CPLEX. However, the case is different in later Section 5.2.2 for a more complex scenario (TSP-mD), where LKH can cost less time than CPLEX when solving the same poi-7- instances to optimum.

Similarly, the column $ \%gap^{\textit{HGA}} = 100 \cdot (M^{\textit{HGA}}-\text{M}^1)/M^{\textit{HGA}} $ captures the improvement over HGA, where positive values indicate better performance by LKH Config1. Considering Table~\ref{tab:fstsp_1LKH_medium} which reports results on medium FSTSP instances, we can conclude that LKH Config1 generally outperforms HGA in instances up to 19 customers, despite being run only once with a 5min runtime limit. 

% FSTSP: cplex 74.85/25=2.994s, lkh 155.97/25=6.2388s

\begin{table}[H]
    \centering
    % \tiny
    % \scriptsize
    \footnotesize  % Slightly smaller to better fit
    % \small
    % \normalsize  % (default size) 
    % \large
    \caption{FSTSP Comparison of HGA and LKH (Config 1) on small instances}
    \label{tab:fstsp_1LKH_small}
    \begin{tabular*}{\textwidth}{@{\extracolsep{\fill}}l@{\;}r@{\;}r@{\;}r@{\;}r@{\;}r@{\;}r@{\;}r@{\;}r@{\;}r@{\;}r@{\;}r@{\;}r@{\;}r@{}}
    %%
    % \begin{tabularx}{\textwidth}{>{\raggedright\arraybackslash}p{1.5cm} >
    % {\raggedleft\arraybackslash}p{0.5cm} >{\raggedleft\arraybackslash}p{0.5cm} >
    % {\raggedleft\arraybackslash}p{0.5cm} >
    % {\raggedleft\arraybackslash}p{1.0cm} >
    % {\raggedleft\arraybackslash}p{0.7cm} >
    % {\raggedleft\arraybackslash}p{0.5cm} >
    % {\raggedleft\arraybackslash}p{0.7cm} >
    % {\raggedleft\arraybackslash}p{1.0cm} >
    % {\raggedleft\arraybackslash}p{0.5cm} >
    % {\raggedleft\arraybackslash}p{0.7cm} >
    % {\raggedleft\arraybackslash}p{0.7cm} >
    % {\raggedleft\arraybackslash}X} 
    % \begin{tabularx}{\textwidth}{l r r r r r r r r r r r r}
    % \begin{tabularx}{\textwidth}{l c c c c c c c c c}
    % \begin{tabularx}{\textwidth}{l *{12}{>{\centering\arraybackslash}X}}
        \toprule
        \multicolumn{2}{c}{} & \multicolumn{3}{c}{\textbf{CPLEX}} & \multicolumn{3}{c}{\textbf{HGA}} & \multicolumn{4}{c}{\textbf{LKH (Config 1)}} & \multicolumn{1}{c}{} & \multicolumn{1}{c}{} \\
        \cmidrule(lr){3-5} \cmidrule(lr){6-8} \cmidrule(lr){9-12}
        \raisebox{1.5ex}[0pt][0pt]{Instance} & \raisebox{1.5ex}[0pt][0pt] {$N_D$} & $z^*$ & Dcus & $t^{\textit{MIP}}$(s) & $M^{\textit{HGA}}$ & Dcus & Time(s) & $\text{M}^1$ & Dcus & Iter & Time(s) & \raisebox{1.5ex}[0pt][0pt] {$\%gap^*$} & \raisebox{1.5ex}[0pt][0pt] {$\%gap^{\textit{HGA}}$} \\
        \midrule
poi-7-1 & 6 & 114.6 & 2 & 2.42 & 128.87 & 1 & 3.3 & 114.6 & 2 & 5904 & 29.2 & 0 & 11.07 \\
poi-7-2 & 6 & 91.6 & 2 & 3.21 & 113.77 & 1 & 2.77 & 91.6 & 2 & 67 & 0.33 & 0 & 19.49 \\
poi-7-3 & 6 & 73.8 & 1 & 2.51 & 96.44 & 1 & 2.83 & 73.8 & 1 & 69 & 0.19 & 0 & 23.48 \\
poi-7-4 & 6 & 100.6 & 3 & 2.53 & 112.34 & 2 & 2.75 & 100.6 & 3 & 2365 & 10.37 & 0 & 10.45 \\
poi-7-5 & 6 & 104.2 & 3 & 2.55 & 124.66 & 1 & 2.56 & 104.2 & 3 & 151 & 0.63 & 0 & 16.41 \\
poi-7-6 & 6 & 93 & 3 & 3.04 & 108.06 & 1 & 2.53 & 93 & 3 & 2327 & 11.25 & 0 & 13.94 \\
poi-7-7 & 6 & 79.4 & 2 & 3.36 & 103.55 & 2 & 4.16 & 79.4 & 2 & 89 & 0.43 & 0 & 23.32 \\
poi-7-8 & 6 & 71 & 2 & 2.46 & 81.42 & 1 & 3.06 & 71 & 2 & 421 & 1.6 & 0 & 12.8 \\
poi-7-9 & 6 & 120.6 & 3 & 2.55 & 133.91 & 1 & 3.23 & 120.6 & 3 & 17866 & 68.46 & 0 & 9.94 \\
poi-7-10 & 6 & 87.4 & 2 & 3.73 & 113.26 & 1 & 2.04 & 87.4 & 2 & 137 & 0.52 & 0 & 22.83 \\
poi-7-11 & 6 & 92.4 & 3 & 3.4 & 109.88 & 1 & 2.16 & 92.4 & 3 & 64 & 0.27 & 0 & 15.91 \\
poi-7-12 & 6 & 97.6 & 2 & 3.44 & 112.87 & 1 & 2.21 & 97.6 & 2 & 818 & 4.04 & 0 & 13.53 \\
poi-7-13 & 6 & 75 & 2 & 3.18 & 78.07 & 1 & 2.31 & 75 & 2 & 994 & 3.5 & 0 & 3.93 \\
poi-7-14 & 6 & 82.6 & 2 & 3.32 & 108.07 & 2 & 2.33 & 82.6 & 2 & 161 & 0.73 & 0 & 23.57 \\
poi-7-15 & 6 & 112.9 & 3 & 3.43 & 126.38 & 2 & 2.24 & 112.9 & 3 & 2356 & 8.85 & 0 & 10.67 \\
poi-7-16 & 6 & 76.5 & 2 & 3.35 & 112.23 & 1 & 2.16 & 76.5 & 2 & 38 & 0.14 & 0 & 31.84 \\
poi-7-17 & 6 & 86.4 & 2 & 3.27 & 107.16 & 1 & 3.48 & 86.4 & 2 & 135 & 0.65 & 0 & 19.37 \\
poi-7-18 & 6 & 76.8 & 2 & 2.53 & 92.51 & 1 & 1.96 & 76.8 & 2 & 964 & 4.22 & 0 & 16.98 \\
poi-7-19 & 6 & 101.4 & 2 & 2.55 & 126.17 & 1 & 2.25 & 101.4 & 2 & 539 & 3.12 & 0 & 19.63 \\
poi-7-20 & 6 & 96.2 & 3 & 2.54 & 107.31 & 1 & 2 & 96.2 & 3 & 62 & 0.12 & 0 & 10.35 \\
poi-7-21 & 6 & 66.1 & 3 & 3.22 & 82.43 & 1 & 2.22 & 66.1 & 3 & 280 & 1.24 & 0 & 19.81 \\
poi-7-22 & 6 & 94.1 & 2 & 2.59 & 113.35 & 1 & 2.09 & 94.1 & 2 & 77 & 0.29 & 0 & 16.98 \\
poi-7-23 & 6 & 129.8 & 2 & 3.48 & - & - & - & 129.8 & 2 & 30 & 0.12 & 0 & - \\
poi-7-24 & 6 & 59.3 & 2 & 3.16 & 83.62 & 1 & 1.91 & 59.3 & 2 & 242 & 1.01 & 0 & 29.08 \\
poi-7-25 & 6 & 92 & 3 & 3.08 & 110.34 & 1 & 2.26 & 92 & 3 & 820 & 4.27 & 0 & 16.62 \\
Average & & & & 2.99 & & & & & & & 6.24 & 0 &  \\
\textbf{poi-10-1} & 9 & \textbf{121.4} & 3 & 1214.33 & \textbf{129.52} & 1 & 4.86 & \textbf{123} & 3 & 282 & 2.1 & \textbf{-1.32} & 5.03 \\
poi-10-2 & 9 & 110 & 3 & 802.15 & 124.16 & 2 & 4.79 & 110 & 3 & 18844 & 164.74 & 0 & 11.4 \\
poi-10-3 & 9 & 116 & 3 & 362 & 150.25 & 1 & 4.95 & 116 & 3 & 5961 & 50.5 & 0 & 22.8 \\
poi-10-4 & 9 & 105.8 & 3 & 453.45 & 125.04 & 1 & 4.55 & 105.8 & 3 & 6969 & 54.64 & 0 & 15.39 \\
poi-10-5 & 9 & 113.9 & 2 & 408.03 & 123.3 & 2 & 5.23 & 113.9 & 2 & 561 & 3.71 & 0 & 7.62 \\
\textbf{poi-10-6} & 9 & \textbf{110.1} & 3 & 487.33 & \textbf{140.12} & 1 & 4.63 & \textbf{117.4} & 3 & 590 & 4.62 & \textbf{-6.63} & 16.21 \\
poi-10-7 & 9 & 85.7 & 3 & 376.66 & 109.76 & 1 & 3.93 & 85.7 & 3 & 24218 & 148.35 & 0 & 21.92 \\
poi-10-8 & 9 & 77.4 & 3 & 1275.99 & 101.07 & 2 & 4.9 & 77.4 & 3 & 822 & 4.39 & 0 & 23.42 \\
\textbf{poi-10-9} & 9 & \textbf{126} & 3 & 315.72 & \textbf{147.59} & 1 & 5.56 & \textbf{127.1} & 4 & 5578 & 48.11 & \textbf{-0.87} & 13.88 \\
\textbf{poi-10-10} & 9 & \textbf{99.7} & 3 & 1507.51 & \textbf{110.97} & 2 & 4.51 & \textbf{100.1} & 3 & 1734 & 14.51 & \textbf{-0.4} & 9.8 \\
poi-10-11 & 9 & 106.8 & 3 & 2377.89 & 131.16 & 1 & 4.25 & 106.8 & 3 & 19502 & 163.84 & 0 & 18.57 \\
poi-10-12 & 9 & 102.5 & 2 & 441.47 & 123.63 & 1 & 4.33 & 102.5 & 2 & 1229 & 6.8 & 0 & 17.09 \\
poi-10-13 & 9 & 89.9 & 2 & 273.7 & 105.38 & 1 & 4.37 & 89.9 & 2 & 733 & 4.65 & 0 & 14.69 \\
poi-10-14 & 9 & 93.3 & 3 & 300.61 & 125.2 & 1 & 5.04 & 93.3 & 3 & 8298 & 62.49 & 0 & 25.48 \\
poi-10-15 & 9 & 116.7 & 4 & 1528.87 & 150.2 & 1 & 3.95 & 116.7 & 4 & 91 & 0.6 & 0 & 22.3 \\
\textbf{poi-10-16} & 9 & \textbf{105.2} & 4 & 383.29 & \textbf{117.69} & 2 & 5.39 & \textbf{106.9} & 4 & 35191 & 270.32 & \textbf{-1.62} & 9.17 \\
poi-10-17 & 9 & 109.5 & 2 & 682.93 & 142.34 & 2 & 6.88 & 109.5 & 2 & 3058 & 24.59 & 0 & 23.07 \\
poi-10-18 & 9 & 98 & 4 & 3506.59 & 110.42 & 2 & 7.14 & 98 & 4 & 2610 & 18.95 & 0 & 11.25 \\
poi-10-19 & 9 & 115.3 & 3 & 501.36 & 125.17 & 2 & 4.95 & 115.3 & 3 & 1623 & 12.98 & 0 & 7.89 \\
poi-10-20 & 9 & 102.1 & 3 & 507.07 & 106.89 & 2 & 4.11 & 102.1 & 3 & 5417 & 46.53 & 0 & 4.48 \\
\textbf{poi-10-21} & 9 & \textbf{75.1} & 3 & 390.05 & \textbf{86.1} & 2 & 6.42 & \textbf{76.3} & 4 & 27957 & 233.1 & \textbf{-1.6} & 11.38 \\
poi-10-22 & 9 & 94.1 & 2 & 512.24 & 103.1 & 2 & 3.9 & 94.1 & 2 & 332 & 2.09 & 0 & 8.73 \\
poi-10-23 & 9 & 124.8 & 2 & 547.76 & - & - & - & 124.8 & 2 & 14966 & 104.42 & 0 & - \\
poi-10-24 & 9 & 95.4 & 3 & 612.32 & 129.24 & 1 & 4.02 & 95.4 & 3 & 9519 & 71.94 & 0 & 26.18 \\
poi-10-25 & 9 & 107.9 & 4 & 848.29 & 115.83 & 2 & 6.71 & 107.9 & 4 & 878 & 6.35 & 0 & 6.85 \\
Average & & & & 824.70 & & & & & & & 53.83 &  &  \\
        \bottomrule
    \end{tabular*}
\end{table}

\begin{table}[H]
    \centering
    % \tiny
    % \scriptsize
    \footnotesize  % Slightly smaller to better fit
    % \small
    % \normalsize  % (default size) 
    % \large
    \caption{FSTSP Comparison of HGA and LKH (Config 1) on medium instances}
    \label{tab:fstsp_1LKH_medium}
    \begin{tabular*}{\textwidth}{@{\extracolsep{\fill}}l@{\;}r@{\;}r@{\;}r@{\;}r@{\;}r@{\;}r@{\;}r@{\;}r@{\;}r@{}}
    %%
    % \begin{tabularx}{\textwidth}{l r r r r r r r r r}
    % \begin{tabularx}{\textwidth}{l c c c c c c c c c}
    % \begin{tabularx}{\textwidth}{l *{9}{>{\centering\arraybackslash}X}}
        \toprule
        \multicolumn{2}{c}{} & \multicolumn{3}{c}{\textbf{HGA}} & \multicolumn{4}{c}{\textbf{LKH (Config 1)}} & \multicolumn{1}{c}{} \\
        \cmidrule(lr){3-5} \cmidrule(lr){6-9}
        \raisebox{1.5ex}[0pt][0pt]{Instance} & \raisebox{1.5ex}[0pt][0pt] {$N_D$} & $M^{\textit{HGA}}$ & Dcus & Time(s) & $\text{M}^1$ & Dcus & Iter & Time(s) & \raisebox{1.5ex}[0pt][0pt] {$\%gap^{\textit{HGA}}$} \\
        \midrule
poi-20-1 & 19 & 177.05 & 3 & 35.11 & 156.5 & 5 & 17588 & 294.25 & 11.61 \\
poi-20-2 & 19 & 167.74 & 3 & 30.81 & 152.3 & 6 & 14680 & 265.29 & 9.2 \\
poi-20-3 & 19 & 215.37 & 2 & 21 & 168.3 & 5 & 17692 & 278.94 & 21.86 \\
poi-20-4 & 19 & 189.63 & 1 & 22.44 & 150 & 6 & 3867 & 68.09 & 20.9 \\
poi-20-5 & 19 & 157.13 & 3 & 17.67 & 147.6 & 6 & 15517 & 254.26 & 6.07 \\
poi-20-6 & 19 & 139.39 & 3 & 23.24 & 132.9 & 6 & 7137 & 124.91 & 4.66 \\
poi-20-7 & 19 & 176.39 & 2 & 25.36 & 158.1 & 7 & 5341 & 92.14 & 10.37 \\
poi-20-8 & 19 & 187.32 & 2 & 22.27 & 162 & 7 & 4067 & 71.44 & 13.52 \\
poi-20-9 & 19 & 163.51 & 2 & 25.39 & 139.7 & 7 & 10518 & 200.23 & 14.56 \\
poi-20-10 & 19 & 156.14 & 3 & 26.94 & 140.1 & 5 & 7732 & 142.31 & 10.27 \\
poi-20-11 & 19 & 163.48 & 2 & 18.28 & 144.7 & 6 & 10807 & 191.07 & 11.49 \\
poi-20-12 & 19 & 170.95 & 2 & 20.27 & 142.8 & 5 & 7062 & 110.5 & 16.47 \\
poi-20-13 & 19 & 174.7 & 2 & 18.75 & 153.4 & 7 & 10752 & 160.27 & 12.19 \\
poi-20-14 & 19 & 168.72 & 2 & 18.99 & 139.6 & 7 & 2316 & 31.79 & 17.26 \\
poi-20-15 & 19 & 194.06 & 2 & 16.91 & 152.1 & 7 & 2340 & 41.95 & 21.62 \\
poi-20-16 & 19 & 163.69 & 2 & 24.94 & 129.8 & 6 & 15053 & 240.9 & 20.7 \\
poi-20-17 & 19 & 166.25 & 2 & 20.94 & 127.1 & 7 & 12900 & 253.4 & 23.55 \\
poi-20-18 & 19 & 179.91 & 2 & 18.43 & 146 & 6 & 4374 & 68.98 & 18.85 \\
poi-20-19 & 19 & 162.57 & 2 & 35.46 & 135 & 6 & 7060 & 106.7 & 16.96 \\
poi-20-20 & 19 & 172.56 & 2 & 19.44 & 150.7 & 7 & 2326 & 35.25 & 12.67 \\
poi-20-21 & 19 & 167.27 & 1 & 17.29 & 137.7 & 6 & 2154 & 35.25 & 17.68 \\
poi-20-22 & 19 & 177.07 & 2 & 29.49 & 151.7 & 5 & 10642 & 182.6 & 14.33 \\
poi-20-23 & 19 & 177.81 & 2 & 18.47 & 153.9 & 7 & 4043 & 60.98 & 13.45 \\
poi-20-24 & 19 & 163.07 & 3 & 19.18 & 140.7 & 6 & 3326 & 61.04 & 13.72 \\
poi-20-25 & 19 & 169.77 & 2 & 19.17 & 148.4 & 6 & 14378 & 239.15 & 12.59 \\
poi-30-1 & 29 & 210.18 & 3 & 102.21 & 183.9 & 7 & 6866 & 233.62 & 12.5 \\
poi-30-2 & 29 & 234 & 3 & 73.12 & 209.3 & 6 & 6609 & 253.53 & 10.56 \\
poi-30-3 & 29 & 212.08 & 2 & 66.34 & 188.1 & 9 & 5252 & 178.68 & 11.31 \\
poi-30-4 & 29 & 207.49 & 4 & 100.99 & 201.7 & 10 & 4863 & 189.66 & 2.79 \\
poi-30-5 & 29 & 189.26 & 3 & 65.89 & 158.6 & 8 & 6215 & 229.86 & 16.2 \\
poi-30-6 & 29 & 193.7 & 4 & 72.82 & 197.9 & 8 & 3310 & 107.18 & -2.17 \\
poi-30-7 & 29 & 180.72 & 5 & 83.14 & 192.1 & 9 & 8704 & 269.11 & -6.3 \\
poi-30-8 & 29 & 190.32 & 4 & 89.47 & 167 & 9 & 8236 & 296.57 & 12.25 \\
poi-30-9 & 29 & 191.31 & 4 & 149.66 & 166.4 & 9 & 6885 & 254.02 & 13.02 \\
poi-30-10 & 29 & 222.74 & 3 & 66.8 & 178.1 & 7 & 3857 & 149.3 & 20.04 \\
poi-30-11 & 29 & 231 & 1 & 65.36 & 185 & 7 & 3746 & 164.28 & 19.91 \\
poi-30-12 & 29 & 238.71 & 2 & 68.52 & 186.1 & 8 & 5821 & 251.45 & 22.04 \\
poi-30-13 & 29 & 225.83 & 2 & 58.63 & 187.6 & 9 & 6926 & 248.86 & 16.93 \\
poi-30-14 & 29 & 198.18 & 4 & 115.16 & 218.5 & 7 & 7128 & 272.72 & -10.25 \\
poi-30-15 & 29 & 175.14 & 4 & 64.45 & 170.1 & 8 & 5592 & 197.46 & 2.88 \\
poi-30-16 & 29 & 189.85 & 5 & 67.12 & 197.5 & 9 & 6414 & 211.29 & -4.03 \\
poi-30-17 & 29 & 224.33 & 1 & 80.35 & 176.4 & 8 & 6393 & 201.34 & 21.37 \\
poi-30-18 & 29 & 227.25 & 3 & 74.81 & 190.2 & 8 & 4488 & 175.54 & 16.3 \\
poi-30-19 & 29 & 184.16 & 4 & 116.48 & 191.9 & 7 & 7328 & 248.2 & -4.2 \\
poi-30-20 & 29 & 183.13 & 3 & 101.17 & 179.6 & 10 & 5059 & 169.78 & 1.93 \\
poi-30-21 & 29 & 194.45 & 3 & 65.91 & 206.5 & 8 & 8840 & 263.8 & -6.2 \\
poi-30-22 & 29 & 225.52 & 2 & 55.36 & 238.4 & 10 & 4940 & 226.82 & -5.71 \\
poi-30-23 & 29 & 213.61 & 1 & 65.94 & 176.4 & 7 & 5516 & 192.62 & 17.42 \\
poi-30-24 & 29 & 188.03 & 2 & 64.19 & 166.4 & 9 & 4357 & 137.14 & 11.5 \\
poi-30-25 & 29 & 227.08 & 2 & 102.09 & 227.3 & 10 & 7431 & 295.88 & -0.1 \\
        \bottomrule
        \multicolumn{10}{r}{Continued on next page}
    \end{tabular*}
\end{table}

\clearpage % Force a page break here

\begin{table}[!t] % Force placement at the absolute top
    \vspace*{0pt} % Remove any vertical space at the top
    \centering
    % \tiny
    % \scriptsize
    \footnotesize  % Slightly smaller to better fit
    % \small
    % \normalsize  % (default size) 
    % \large
    % Use the caption command but suppress the numbering
    \refstepcounter{table}
    \addtocounter{table}{-1} % Keep the table number the same
    \makeatletter
    \def\@captiontype{table}
    {\@makecaption{\tablename~\thetable}{(continued)}}
    \vspace{0.5em}
    \makeatother
    \begin{tabular*}{\textwidth}{@{\extracolsep{\fill}}l@{\;}r@{\;}r@{\;}r@{\;}r@{\;}r@{\;}r@{\;}r@{\;}r@{\;}r@{}}
    %%
    % \begin{tabularx}{\textwidth}{l r r r r r r r r r}
    % \begin{tabularx}{\textwidth}{l c c c c c c c c c}
    % \begin{tabularx}{\textwidth}{l *{9}{>{\centering\arraybackslash}X}}
        \toprule
        \multicolumn{2}{c}{} & \multicolumn{3}{c}{\textbf{HGA}} & \multicolumn{4}{c}{\textbf{LKH (Config 1)}} & \multicolumn{1}{c}{} \\
        \cmidrule(lr){3-5} \cmidrule(lr){6-9}
        \raisebox{1.5ex}[0pt][0pt]{Instance} & \raisebox{1.5ex}[0pt][0pt] {$N_D$} & $M^{\textit{HGA}}$ & Dcus & Time(s) & $\text{M}^1$ & Dcus & Iter & Time(s) & \raisebox{1.5ex}[0pt][0pt] {$\%gap^{\textit{HGA}}$} \\
        \midrule
poi-40-1 & 39 & 234.47 & 3 & 169.27 & 231 & 10 & 3957 & 217.79 & 1.48 \\
poi-40-2 & 39 & 237.08 & 1 & 193.04 & 206.9 & 9 & 4444 & 284.83 & 12.73 \\
poi-40-3 & 39 & 238.83 & 6 & 168.78 & 232.9 & 9 & 5994 & 289.25 & 2.48 \\
poi-40-4 & 39 & 244.33 & 2 & 251.9 & 271.4 & 8 & 4237 & 272.86 & -11.08 \\
poi-40-5 & 39 & 257.48 & 2 & 212.78 & 258.3 & 11 & 5163 & 287.51 & -0.32 \\
poi-40-6 & 39 & 189.56 & 5 & 262.16 & 267.1 & 12 & 6016 & 288.63 & -40.91 \\
poi-40-7 & 39 & 227.86 & 3 & 276.03 & 229.2 & 10 & 4049 & 249.14 & -0.59 \\
poi-40-8 & 39 & 261.6 & 5 & 276.24 & 474.8 & 16 & 4029 & 298.45 & -81.5 \\
poi-40-9 & 39 & 232.27 & 4 & 178.88 & 273.8 & 12 & 3794 & 222.85 & -17.88 \\
poi-40-10 & 39 & 229.02 & 2 & 318.15 & 236.5 & 8 & 4902 & 296.22 & -3.27 \\
poi-40-11 & 39 & 236.14 & 4 & 263.04 & 257.6 & 9 & 5466 & 261.54 & -9.09 \\
poi-40-12 & 39 & 226.88 & 4 & 190.11 & 275.1 & 12 & 3999 & 214.39 & -21.25 \\
poi-40-13 & 39 & 245.74 & 4 & 232.79 & 282.2 & 11 & 3493 & 208.97 & -14.84 \\
poi-40-14 & 39 & 255.15 & 3 & 219.76 & 271.3 & 10 & 4297 & 229.43 & -6.33 \\
poi-40-15 & 39 & 243.37 & 3 & 136.5 & 275.9 & 10 & 5145 & 259.94 & -13.37 \\
poi-40-16 & 39 & 236.97 & 3 & 172.75 & 293.3 & 10 & 3920 & 213.32 & -23.77 \\
poi-40-17 & 39 & 240.74 & 5 & 329.18 & 296.1 & 13 & 4216 & 179.04 & -23 \\
poi-40-18 & 39 & 273.33 & 3 & 147.27 & 329.8 & 13 & 4616 & 221.29 & -20.66 \\
poi-40-19 & 39 & 222.19 & 2 & 165.11 & 214.8 & 11 & 4599 & 254.95 & 3.33 \\
poi-40-20 & 39 & 238.62 & 2 & 230.07 & 244.7 & 9 & 3091 & 181.37 & -2.55 \\
poi-40-21 & 39 & 253.79 & 3 & 169.56 & 229.6 & 6 & 4979 & 276.71 & 9.53 \\
poi-40-22 & 39 & 260.28 & 4 & 147.67 & 255.8 & 11 & 4118 & 260.16 & 1.72 \\
poi-40-23 & 39 & 275.4 & 2 & 267.38 & 330 & 15 & 5311 & 290.45 & -19.83 \\
poi-40-24 & 39 & 213.36 & 4 & 331.6 & 279.8 & 10 & 4211 & 298.27 & -31.14 \\
poi-40-25 & 39 & 228.47 & 5 & 188.2 & 285.1 & 9 & 4157 & 280.93 & -24.79 \\
        \bottomrule
    \end{tabular*}
\end{table}

Based on the comparison above between HGA and LKH Config1 in Tables~\ref{tab:fstsp_1LKH_small}--\ref{tab:fstsp_1LKH_medium}, negative values in $\%gap^{\textit{HGA}}$ imply that LKH underperforms in some instances. To improve LKH outputs, we extend its runtime from 5min to 4h, as shown in Table~\ref{tab:fstsp_1LKH_longer}, where previous Dcus values are listed in parentheses. Note that with the exception of poi-40-17, all gaps become positive, demonstrating better performance. In particular, LKH Config1 reaches optimal solutions for instances poi-10-1, poi-10-6, poi-10-9, and poi-10-16 given 4h, highlighted in bold. Interestingly, the improved solution is not obviously correlated with increasing or decreasing Dcus. These suggest that LKH could improve its makespan to outperform HGA when given a longer runtime, which is not driven solely by the number of drone-service customers. 

\begin{table}[H]
    \centering
    % \tiny
    % \scriptsize
    \footnotesize  % Slightly smaller to better fit
    % \small
    % \normalsize  % (default size) 
    % \large
    \caption{FSTSP Comparison of HGA and LKH (Config 1 4h)}
    \label{tab:fstsp_1LKH_longer}
    \begin{tabular*}{\textwidth}{@{\extracolsep{\fill}}l@{\;}r@{\;}r@{\;}r@{\;}r@{\;}r@{\;}r@{\;}r@{\;}r@{\;}r@{}}
    %%
    % \begin{tabularx}{\textwidth}{l r r r r r r r r r}
    % \begin{tabularx}{\textwidth}{l c c c c c c c c c}
    % \begin{tabularx}{\textwidth}{l *{9}{>{\centering\arraybackslash}X}}
        \toprule
        \multicolumn{2}{c}{} & \multicolumn{3}{c}{\textbf{HGA}} & \multicolumn{4}{c}{\textbf{LKH (Config 1 4h)}} & \multicolumn{1}{c}{} \\
        \cmidrule(lr){3-5} \cmidrule(lr){6-9}
        \raisebox{1.5ex}[0pt][0pt]{Instance} & \raisebox{1.5ex}[0pt][0pt] {$N_D$} & $M^{\textit{HGA}}$ & Dcus & Time(s) & $\text{M}^1$ & Dcus & Iter & Time(s) & \raisebox{1.5ex}[0pt][0pt] {$\%gap^{\textit{HGA}}$} \\
        \midrule
\textbf{poi-10-1} & 9 & 129.52 & 1 & 4.86 & \textbf{121.4} & 3 & 645270 & 5037.54 & 6.27 \\
\textbf{poi-10-6} & 9 & 140.12 & 1 & 4.63 & \textbf{110.1} & 3 & 44011 & 353.99 & 21.42 \\
\textbf{poi-10-9} & 9 & 147.59 & 1 & 5.56 & \textbf{126} & (4) 3 & 78847 & 638.88 & 14.63 \\
poi-10-10 & 9 & 110.97 & 2 & 4.51 & 100.1 & (3) 2 & 55208 & 460.01 & 9.8 \\
\textbf{poi-10-16} & 9 & 117.69 & 2 & 5.39 & \textbf{105.2} & 4 & 44203 & 340.04 & 10.61 \\
poi-10-21 & 9 & 86.1 & 2 & 6.42 & 76.2 & (4) 3 & 126498 & 1075.14 & 11.5 \\
poi-30-6 & 29 & 193.7 & 4 & 72.82 & 161.5 & 8 & 334367 & 12423.06 & 16.62 \\
poi-30-7 & 29 & 180.72 & 5 & 83.14 & 148.7 & (9) 10 & 154150 & 5838.35 & 17.72 \\
poi-30-14 & 29 & 198.18 & 4 & 115.16 & 160.4 & (7) 8 & 85218 & 3509 & 19.06 \\
poi-30-16 & 29 & 189.85 & 5 & 67.12 & 172 & (9) 8 & 26652 & 947.91 & 9.4 \\
poi-30-19 & 29 & 184.16 & 4 & 116.48 & 173.9 & 7 & 81047 & 2786.83 & 5.57 \\
poi-30-21 & 29 & 194.45 & 3 & 65.91 & 178.9 & 8 & 82088 & 2583.13 & 8 \\
poi-30-22 & 29 & 225.52 & 2 & 55.36 & 179.6 & (10) 8 & 305266 & 14309.34 & 20.36 \\
poi-30-25 & 29 & 227.08 & 2 & 102.09 & 221.8 & 10 & 159264 & 6659.3 & 2.33 \\
poi-40-4 & 39 & 244.33 & 2 & 251.9 & 231.1 & (8) 12 & 163357 & 11868.63 & 5.41 \\
poi-40-5 & 39 & 257.48 & 2 & 212.78 & 211 & (11) 9 & 89443 & 5001.01 & 18.05 \\
poi-40-6 & 39 & 189.56 & 5 & 262.16 & 175.7 & (12) 10 & 109075 & 7215.39 & 7.31 \\
poi-40-7 & 39 & 227.86 & 3 & 276.03 & 193.2 & (10) 12 & 108053 & 8161.67 & 15.21 \\
poi-40-8 & 39 & 261.6 & 5 & 276.24 & 215.2 & (16) 11 & 229451 & 14268.98 & 17.74 \\
poi-40-9 & 39 & 232.27 & 4 & 178.88 & 212.4 & (12) 8 & 128096 & 8432.08 & 8.55 \\
poi-40-10 & 39 & 229.02 & 2 & 318.15 & 195 & (8) 12 & 189154 & 11683.48 & 14.85 \\
poi-40-11 & 39 & 236.14 & 4 & 263.04 & 207.2 & (9) 12 & 27291 & 1515.24 & 12.26 \\
poi-40-12 & 39 & 226.88 & 4 & 190.11 & 201.9 & (12) 8 & 44757 & 3271.06 & 11.01 \\
poi-40-13 & 39 & 245.74 & 4 & 232.79 & 200.8 & (11) 10 & 91300 & 5896.61 & 18.29 \\
poi-40-14 & 39 & 255.15 & 3 & 219.76 & 239.6 & (10) 11 & 223884 & 14105.31 & 6.09 \\
poi-40-15 & 39 & 243.37 & 3 & 136.5 & 213.2 & 10 & 49305 & 4239.12 & 12.4 \\
poi-40-16 & 39 & 236.97 & 3 & 172.75 & 200.3 & (10) 8 & 217210 & 14189.27 & 15.47 \\
poi-40-17 & 39 & 240.74 & 5 & 329.18 & 246.9 & (13) 12 & 88549 & 6100.75 & -2.56 \\
poi-40-18 & 39 & 273.33 & 3 & 147.27 & 207.9 & (13) 12 & 255503 & 13714.64 & 23.94 \\
poi-40-20 & 39 & 238.62 & 2 & 230.07 & 200.9 & 9 & 25421 & 1672.51 & 15.81 \\
poi-40-23 & 39 & 275.4 & 2 & 267.38 & 216 & (15) 12 & 211455 & 13243.52 & 21.57 \\
poi-40-24 & 39 & 213.36 & 4 & 331.6 & 201.9 & (10) 12 & 112104 & 7261.89 & 5.37 \\
poi-40-25 & 39 & 228.47 & 5 & 188.2 & 217.7 & (9) 11 & 199460 & 14128.9 & 4.71 \\
        \bottomrule
    \end{tabular*}
\end{table}

After demonstrating the effectiveness of LKH Config1 compared to HGA, we evaluate its two configurations in Tables~\ref{tab:fstsp_2LKH_small}--\ref{tab:fstsp_2LKH_medium}. Recall that the second configuration allows the truck to revisit nodes after it leaves. Thus, we highlight in bold the three instances where LKH Config2 leverages node revisits: poi-7-14, poi-7-23, and poi-40-6. Moreover, Config2 manages to generate optimal solutions for instances poi-10-1, poi-10-6, poi-10-9, poi-10-10, and poi-10-16 in significantly less time than CPLEX (see Table~\ref{tab:fstsp_1LKH_small}), while Config1 fails to achieve optimality in 5min. Then, the positive gap $\%gap^{1,2} = 100 \cdot (\text{M}^1-\text{M}^2)/\text{M}^1$, reflects the makespan reduction achieved by switching from Config1 to Config2. Note that the values $\%gap^{1,2}$ fluctuate across problem instances, indicating that Config2 does not consistently outperform Config1. This variability may be attributed to the removal of the no-revisit restriction, which expands the search space and affects convergence path within the 5min runtime. Consequently, LKH performance gains may depend on specific instance characteristics when turning to Config2.

% second comparison
\begin{table}[H]
    \centering
    % \tiny
    % \scriptsize
    \footnotesize  % Slightly smaller to better fit
    % \small
    % \normalsize  % (default size) 
    % \large
    \caption{FSTSP Comparison of LKH (Config 1) and LKH (Config 2) on small instances}
    \label{tab:fstsp_2LKH_small}
    \begin{tabular*}{\textwidth}{@{\extracolsep{\fill}}l@{\;}r@{\;}r@{\;}r@{\;}r@{\;}r@{\;}r@{\;}r@{\;}r@{\;}r@{\;}r@{}}
    %%
    % \begin{tabularx}{\textwidth}{l r r r r r r r r r r}
    % \begin{tabularx}{\textwidth}{l c c c c c c c c c c}
    % \begin{tabularx}{\textwidth}{l *{10}{>{\centering\arraybackslash}X}}
        \toprule
        \multicolumn{2}{c}{} & \multicolumn{4}{c}{\textbf{LKH (Config 1)}} & \multicolumn{4}{c}{\textbf{LKH (Config 2)}} & \multicolumn{1}{c}{}\\
        \cmidrule(lr){3-6} \cmidrule(lr){7-10}
        \raisebox{1.5ex}[0pt][0pt]{Instance} & \raisebox{1.5ex}[0pt][0pt] {$N_D$} & $\text{M}^1$ & Dcus & Iter & Time(s) & $\text{M}^2$ & Dcus & Iter & Time(s) & \raisebox{1.5ex}[0pt][0pt] {$\%gap^{1,2}$} \\
        \midrule
poi-7-1 & 6 & 114.6 & 2 & 5904 & 29.2 & 114.6 & 2 & 917 & 4.43 & 0 \\
poi-7-2 & 6 & 91.6 & 2 & 67 & 0.33 & 91.6 & 2 & 185 & 0.99 & 0 \\
poi-7-3 & 6 & 73.8 & 1 & 69 & 0.19 & 73.8 & 1 & 11 & 0.04 & 0 \\
poi-7-4 & 6 & 100.6 & 3 & 2365 & 10.37 & 100.6 & 3 & 100 & 0.41 & 0 \\
poi-7-5 & 6 & 104.2 & 3 & 151 & 0.63 & 104.2 & 3 & 299 & 1.54 & 0 \\
poi-7-6 & 6 & 93 & 3 & 2327 & 11.25 & 93 & 3 & 3347 & 14.13 & 0 \\
poi-7-7 & 6 & 79.4 & 2 & 89 & 0.43 & 79.4 & 2 & 105 & 0.34 & 0 \\
poi-7-8 & 6 & 71 & 2 & 421 & 1.6 & 71 & 2 & 303 & 1.36 & 0 \\
poi-7-9 & 6 & 120.6 & 3 & 17866 & 68.46 & 120.6 & 3 & 6792 & 30.23 & 0 \\
poi-7-10 & 6 & 87.4 & 2 & 137 & 0.52 & 87.4 & 2 & 167 & 0.95 & 0 \\
poi-7-11 & 6 & 92.4 & 3 & 64 & 0.27 & 92.4 & 3 & 707 & 3.28 & 0 \\
poi-7-12 & 6 & 97.6 & 2 & 818 & 4.04 & 97.6 & 2 & 47 & 0.25 & 0 \\
poi-7-13 & 6 & 75 & 2 & 994 & 3.5 & 75 & 2 & 1600 & 5.87 & 0 \\
\textbf{poi-7-14} & 6 & \textbf{82.6} & 2 & 161 & 0.73 & \textbf{82.6 revisit} & 2 & 22 & 0.12 & \textbf{0} \\
poi-7-15 & 6 & 112.9 & 3 & 2356 & 8.85 & 112.9 & 3 & 678 & 2.99 & 0 \\
poi-7-16 & 6 & 76.5 & 2 & 38 & 0.14 & 76.5 & 2 & 129 & 0.5 & 0 \\
poi-7-17 & 6 & 86.4 & 2 & 135 & 0.65 & 86.4 & 2 & 20 & 0.1 & 0 \\
poi-7-18 & 6 & 76.8 & 2 & 964 & 4.22 & 76.8 & 2 & 63 & 0.24 & 0 \\
poi-7-19 & 6 & 101.4 & 2 & 539 & 3.12 & 101.4 & 2 & 66 & 0.33 & 0 \\
poi-7-20 & 6 & 96.2 & 3 & 62 & 0.12 & 96.2 & 3 & 545 & 2.44 & 0 \\
poi-7-21 & 6 & 66.1 & 3 & 280 & 1.24 & 66.1 & 3 & 764 & 2.48 & 0 \\
poi-7-22 & 6 & 94.1 & 2 & 77 & 0.29 & 94.1 & 2 & 414 & 1.73 & 0 \\
\textbf{poi-7-23} & 6 & \textbf{129.8} & 2 & 30 & 0.12 & \textbf{129.8 revisit} & 2 & 33 & 0.17 & \textbf{0} \\
poi-7-24 & 6 & 59.3 & 2 & 242 & 1.01 & 59.3 & 2 & 217 & 1.12 & 0 \\
poi-7-25 & 6 & 92 & 3 & 820 & 4.27 & 92 & 3 & 242 & 1.26 & 0 \\
poi-10-1 & 9 & 123 & 3 & 282 & 2.1 & \textbf{121.4} & 3 & 4688 & 40.73 & 1.3 \\
poi-10-2 & 9 & 110 & 3 & 18844 & 164.74 & 110 & 3 & 2835 & 25.86 & 0 \\
poi-10-3 & 9 & 116 & 3 & 5961 & 50.5 & 116 & 3 & 11025 & 91.29 & 0 \\
poi-10-4 & 9 & 105.8 & 3 & 6969 & 54.64 & 105.8 & 3 & 838 & 6.28 & 0 \\
poi-10-5 & 9 & 113.9 & 2 & 561 & 3.71 & 113.9 & 2 & 1776 & 13.93 & 0 \\
poi-10-6 & 9 & 117.4 & 3 & 590 & 4.62 & \textbf{110.1} & 3 & 666 & 4.57 & 6.22 \\
poi-10-7 & 9 & 85.7 & 3 & 24218 & 148.35 & 85.7 & 3 & 539 & 3.08 & 0 \\
poi-10-8 & 9 & 77.4 & 3 & 822 & 4.39 & 77.4 & 3 & 26330 & 255.23 & 0 \\
poi-10-9 & 9 & 127.1 & 4 & 5578 & 48.11 & \textbf{126} & 3 & 1853 & 15.19 & 0.87 \\
poi-10-10 & 9 & 100.1 & 3 & 1734 & 14.51 & \textbf{99.7} & 3 & 487 & 3.8 & 0.4 \\
poi-10-11 & 9 & 106.8 & 3 & 19502 & 163.84 & 110.8 & 3 & 857 & 7.28 & -3.75 \\
poi-10-12 & 9 & 102.5 & 2 & 1229 & 6.8 & 102.5 & 2 & 1887 & 16.15 & 0 \\
poi-10-13 & 9 & 89.9 & 2 & 733 & 4.65 & 89.9 & 2 & 188 & 1.26 & 0 \\
poi-10-14 & 9 & 93.3 & 3 & 8298 & 62.49 & 93.3 & 3 & 10933 & 82.13 & 0 \\
poi-10-15 & 9 & 116.7 & 4 & 91 & 0.6 & 117.8 & 3 & 32937 & 250.08 & -0.94 \\
poi-10-16 & 9 & 106.9 & 4 & 35191 & 270.32 & \textbf{105.2} & 4 & 4144 & 36.46 & 1.59 \\
poi-10-17 & 9 & 109.5 & 2 & 3058 & 24.59 & 109.5 & 2 & 2735 & 22.39 & 0 \\
poi-10-18 & 9 & 98 & 4 & 2610 & 18.95 & 98 & 4 & 5854 & 41.61 & 0 \\
poi-10-19 & 9 & 115.3 & 3 & 1623 & 12.98 & 115.4 & 3 & 721 & 5.2 & -0.09 \\
poi-10-20 & 9 & 102.1 & 3 & 5417 & 46.53 & 102.1 & 3 & 722 & 6.08 & 0 \\
poi-10-21 & 9 & 76.3 & 4 & 27957 & 233.1 & 76.2 & 3 & 3236 & 27.03 & 0.13 \\
poi-10-22 & 9 & 94.1 & 2 & 332 & 2.09 & 94.1 & 2 & 1007 & 9.09 & 0 \\
poi-10-23 & 9 & 124.8 & 2 & 14966 & 104.42 & 124.8 & 2 & 484 & 4.33 & 0 \\
poi-10-24 & 9 & 95.4 & 3 & 9519 & 71.94 & 95.4 & 3 & 9390 & 78.46 & 0 \\
poi-10-25 & 9 & 107.9 & 4 & 878 & 6.35 & 108.8 & 3 & 6731 & 54.04 & -0.83 \\
        \bottomrule
    \end{tabular*}
\end{table}

\begin{table}[H]
    \centering
    % \tiny
    % \scriptsize
    \footnotesize  % Slightly smaller to better fit
    % \small
    % \normalsize  % (default size) 
    % \large
    \caption{FSTSP Comparison of LKH (Config 1) and LKH (Config 2) on medium instances}
    \label{tab:fstsp_2LKH_medium}
    \begin{tabular*}{\textwidth}{@{\extracolsep{\fill}}l@{\;}r@{\;}r@{\;}r@{\;}r@{\;}r@{\;}r@{\;}r@{\;}r@{\;}r@{\;}r@{}}
    %%
    % \begin{tabularx}{\textwidth}{l r r r r r r r r r r}
    % \begin{tabularx}{\textwidth}{l c c c c c c c c c c}
    % \begin{tabularx}{\textwidth}{l *{10}{>{\centering\arraybackslash}X}}
        \toprule
        \multicolumn{2}{c}{} & \multicolumn{4}{c}{\textbf{LKH (Config 1)}} & \multicolumn{4}{c}{\textbf{LKH (Config 2)}} & \multicolumn{1}{c}{}\\
        \cmidrule(lr){3-6} \cmidrule(lr){7-10}
        \raisebox{1.5ex}[0pt][0pt]{Instance} & \raisebox{1.5ex}[0pt][0pt] {$N_D$} & $\text{M}^1$ & Dcus & Iter & Time(s) & $\text{M}^2$ & Dcus & Iter & Time(s) & \raisebox{1.5ex}[0pt][0pt] {$\%gap^{1,2}$} \\
        \midrule
poi-20-1 & 19 & 156.5 & 5 & 17588 & 294.25 & 157.8 & 6 & 14898 & 299.98 & -0.83 \\
poi-20-2 & 19 & 152.3 & 6 & 14680 & 265.29 & 145.5 & 6 & 1778 & 35.65 & 4.46 \\
poi-20-3 & 19 & 168.3 & 5 & 17692 & 278.94 & 162.1 & 7 & 4864 & 88.38 & 3.68 \\
poi-20-4 & 19 & 150 & 6 & 3867 & 68.09 & 139.9 & 5 & 1044 & 20.1 & 6.73 \\
poi-20-5 & 19 & 147.6 & 6 & 15517 & 254.26 & 141.6 & 5 & 8686 & 177.41 & 4.07 \\
poi-20-6 & 19 & 132.9 & 6 & 7137 & 124.91 & 132.5 & 5 & 3160 & 49 & 0.3 \\
poi-20-7 & 19 & 158.1 & 7 & 5341 & 92.14 & 159.3 & 7 & 12739 & 216.84 & -0.76 \\
poi-20-8 & 19 & 162 & 7 & 4067 & 71.44 & 164.3 & 6 & 3938 & 76.35 & -1.42 \\
poi-20-9 & 19 & 139.7 & 7 & 10518 & 200.23 & 139.7 & 7 & 17277 & 297.82 & 0 \\
poi-20-10 & 19 & 140.1 & 5 & 7732 & 142.31 & 140.1 & 5 & 6238 & 125.61 & 0 \\
poi-20-11 & 19 & 144.7 & 6 & 10807 & 191.07 & 144.5 & 6 & 6715 & 142.12 & 0.14 \\
poi-20-12 & 19 & 142.8 & 5 & 7062 & 110.5 & 142.1 & 7 & 14602 & 237.37 & 0.49 \\
poi-20-13 & 19 & 153.4 & 7 & 10752 & 160.27 & 146.5 & 7 & 2558 & 60.02 & 4.5 \\
poi-20-14 & 19 & 139.6 & 7 & 2316 & 31.79 & 136.2 & 8 & 8678 & 185.9 & 2.44 \\
poi-20-15 & 19 & 152.1 & 7 & 2340 & 41.95 & 148.9 & 6 & 7999 & 153.68 & 2.1 \\
poi-20-16 & 19 & 129.8 & 6 & 15053 & 240.9 & 130.8 & 6 & 4901 & 100.39 & -0.77 \\
poi-20-17 & 19 & 127.1 & 7 & 12900 & 253.4 & 124.8 & 5 & 6770 & 138.89 & 1.81 \\
poi-20-18 & 19 & 146 & 6 & 4374 & 68.98 & 144.9 & 7 & 4295 & 94.68 & 0.75 \\
poi-20-19 & 19 & 135 & 6 & 7060 & 106.7 & 134.6 & 5 & 12323 & 233.56 & 0.3 \\
poi-20-20 & 19 & 150.7 & 7 & 2326 & 35.25 & 153.7 & 8 & 14122 & 226.9 & -1.99 \\
poi-20-21 & 19 & 137.7 & 6 & 2154 & 35.25 & 141 & 5 & 5726 & 114.81 & -2.4 \\
poi-20-22 & 19 & 151.7 & 5 & 10642 & 182.6 & 163.1 & 7 & 9708 & 171.84 & -7.51 \\
poi-20-23 & 19 & 153.9 & 7 & 4043 & 60.98 & 139.3 & 3 & 3373 & 62.3 & 9.49 \\
poi-20-24 & 19 & 140.7 & 6 & 3326 & 61.04 & 130.1 & 4 & 1678 & 30.3 & 7.53 \\
poi-20-25 & 19 & 148.4 & 6 & 14378 & 239.15 & 149.4 & 6 & 10167 & 203.22 & -0.67 \\
poi-30-1 & 29 & 183.9 & 7 & 6866 & 233.62 & 186.5 & 8 & 2577 & 118.35 & -1.41 \\
poi-30-2 & 29 & 209.3 & 6 & 6609 & 253.53 & 189.2 & 7 & 418 & 17.3 & 9.6 \\
poi-30-3 & 29 & 188.1 & 9 & 5252 & 178.68 & 181.5 & 8 & 6418 & 261.48 & 3.51 \\
poi-30-4 & 29 & 201.7 & 10 & 4863 & 189.66 & 186.8 & 8 & 3347 & 143.66 & 7.39 \\
poi-30-5 & 29 & 158.6 & 8 & 6215 & 229.86 & 174.1 & 6 & 6822 & 299.46 & -9.77 \\
poi-30-6 & 29 & 197.9 & 8 & 3310 & 107.18 & 187.7 & 9 & 7200 & 294.93 & 5.15 \\
poi-30-7 & 29 & 192.1 & 9 & 8704 & 269.11 & 165.5 & 8 & 3013 & 136.88 & 13.85 \\
poi-30-8 & 29 & 167 & 9 & 8236 & 296.57 & 175.8 & 7 & 3817 & 171.32 & -5.27 \\
poi-30-9 & 29 & 166.4 & 9 & 6885 & 254.02 & 174 & 8 & 4930 & 210.87 & -4.57 \\
poi-30-10 & 29 & 178.1 & 7 & 3857 & 149.3 & 182.8 & 7 & 4331 & 215.64 & -2.64 \\
poi-30-11 & 29 & 185 & 7 & 3746 & 164.28 & 182.9 & 6 & 5418 & 261.62 & 1.14 \\
poi-30-12 & 29 & 186.1 & 8 & 5821 & 251.45 & 187.4 & 7 & 7184 & 255.02 & -0.7 \\
poi-30-13 & 29 & 187.6 & 9 & 6926 & 248.86 & 198.8 & 8 & 3100 & 150.42 & -5.97 \\
poi-30-14 & 29 & 218.5 & 7 & 7128 & 272.72 & 196 & 6 & 7674 & 297.87 & 10.3 \\
poi-30-15 & 29 & 170.1 & 8 & 5592 & 197.46 & 163.6 & 7 & 7133 & 267.94 & 3.82 \\
poi-30-16 & 29 & 197.5 & 9 & 6414 & 211.29 & 207 & 7 & 4445 & 188.85 & -4.81 \\
poi-30-17 & 29 & 176.4 & 8 & 6393 & 201.34 & 199.3 & 8 & 6837 & 253.5 & -12.98 \\
poi-30-18 & 29 & 190.2 & 8 & 4488 & 175.54 & 180.5 & 11 & 5086 & 259.07 & 5.1 \\
poi-30-19 & 29 & 191.9 & 7 & 7328 & 248.2 & 173.8 & 7 & 8199 & 285.31 & 9.43 \\
poi-30-20 & 29 & 179.6 & 10 & 5059 & 169.78 & 165.8 & 7 & 4703 & 204.16 & 7.68 \\
poi-30-21 & 29 & 206.5 & 8 & 8840 & 263.8 & 193.9 & 8 & 8229 & 270.48 & 6.1 \\
poi-30-22 & 29 & 238.4 & 10 & 4940 & 226.82 & 216.1 & 7 & 6200 & 292.34 & 9.35 \\
poi-30-23 & 29 & 176.4 & 7 & 5516 & 192.62 & 178.5 & 6 & 8138 & 284.98 & -1.19 \\
poi-30-24 & 29 & 166.4 & 9 & 4357 & 137.14 & 158.4 & 8 & 5732 & 248 & 4.81 \\
poi-30-25 & 29 & 227.3 & 10 & 7431 & 295.88 & 215.7 & 7 & 5734 & 243.78 & 5.1 \\
        \bottomrule
        \multicolumn{11}{r}{Continued on next page}
    \end{tabular*}
\end{table}

\clearpage % Force a page break here

\begin{table}[!t] % Force placement at the absolute top
    \vspace*{0pt} % Remove any vertical space at the top
    \centering
    % \tiny
    % \scriptsize
    \footnotesize  % Slightly smaller to better fit
    % \small
    % \normalsize  % (default size) 
    % \large
    \refstepcounter{table}
    \addtocounter{table}{-1} % Keep the table number the same
    \makeatletter
    \def\@captiontype{table}
    {\@makecaption{\tablename~\thetable}{(continued)}}
    \vspace{0.5em}
    \makeatother
    \begin{tabular*}{\textwidth}{@{\extracolsep{\fill}}l@{\;}r@{\;}r@{\;}r@{\;}r@{\;}r@{\;}r@{\;}r@{\;}r@{\;}r@{\;}r@{}}
    %%
    % \begin{tabularx}{\textwidth}{l r r r r r r r r r r}
    % \begin{tabularx}{\textwidth}{l c c c c c c c c c c}
    % \begin{tabularx}{\textwidth}{l *{10}{>{\centering\arraybackslash}X}}
        \toprule
        \multicolumn{2}{c}{} & \multicolumn{4}{c}{\textbf{LKH (Config 1)}} & \multicolumn{4}{c}{\textbf{LKH (Config 2)}} & \multicolumn{1}{c}{}\\
        \cmidrule(lr){3-6} \cmidrule(lr){7-10}
        \raisebox{1.5ex}[0pt][0pt]{Instance} & \raisebox{1.5ex}[0pt][0pt] {$N_D$} & $\text{M}^1$ & Dcus & Iter & Time(s) & $\text{M}^2$ & Dcus & Iter & Time(s) & \raisebox{1.5ex}[0pt][0pt] {$\%gap^{1,2}$} \\
        \midrule
poi-40-1 & 39 & 231 & 10 & 3957 & 217.79 & 218.1 & 11 & 4335 & 298.05 & 5.58 \\
poi-40-2 & 39 & 206.9 & 9 & 4444 & 284.83 & 215.2 & 9 & 3427 & 257.29 & -4.01 \\
poi-40-3 & 39 & 232.9 & 9 & 5994 & 289.25 & 261.2 & 10 & 2326 & 191.36 & -12.15 \\
poi-40-4 & 39 & 271.4 & 8 & 4237 & 272.86 & 258.3 & 10 & 3440 & 288.73 & 4.83 \\
poi-40-5 & 39 & 258.3 & 11 & 5163 & 287.51 & 305 & 12 & 3764 & 292.66 & -18.08 \\
\textbf{poi-40-6} & 39 & \textbf{267.1} & 12 & 6016 & 288.63 & \textbf{416.6 revisit} & 21 & 2430 & 291.39 & \textbf{-55.97} \\
poi-40-7 & 39 & 229.2 & 10 & 4049 & 249.14 & 217.7 & 11 & 4041 & 285.63 & 5.02 \\
poi-40-8 & 39 & 474.8 & 16 & 4029 & 298.45 & 261.1 & 9 & 3288 & 274.87 & 45.01 \\
poi-40-9 & 39 & 273.8 & 12 & 3794 & 222.85 & 268.8 & 9 & 2697 & 200.52 & 1.83 \\
poi-40-10 & 39 & 236.5 & 8 & 4902 & 296.22 & 284.3 & 10 & 3412 & 292.5 & -20.21 \\
poi-40-11 & 39 & 257.6 & 9 & 5466 & 261.54 & 250.5 & 9 & 4204 & 283.11 & 2.76 \\
poi-40-12 & 39 & 275.1 & 12 & 3999 & 214.39 & 259.8 & 10 & 3141 & 276.6 & 5.56 \\
poi-40-13 & 39 & 282.2 & 11 & 3493 & 208.97 & 278.3 & 10 & 2830 & 245.08 & 1.38 \\
poi-40-14 & 39 & 271.3 & 10 & 4297 & 229.43 & 262.6 & 11 & 3921 & 288.03 & 3.21 \\
poi-40-15 & 39 & 275.9 & 10 & 5145 & 259.94 & 218.6 & 10 & 4086 & 282.28 & 20.77 \\
poi-40-16 & 39 & 293.3 & 10 & 3920 & 213.32 & 256.2 & 9 & 2955 & 252.22 & 12.65 \\
poi-40-17 & 39 & 296.1 & 13 & 4216 & 179.04 & 281.2 & 12 & 3095 & 244.6 & 5.03 \\
poi-40-18 & 39 & 329.8 & 13 & 4616 & 221.29 & 285 & 11 & 3282 & 230.63 & 13.58 \\
poi-40-19 & 39 & 214.8 & 11 & 4599 & 254.95 & 248.5 & 12 & 2652 & 267.91 & -15.69 \\
poi-40-20 & 39 & 244.7 & 9 & 3091 & 181.37 & 229.7 & 10 & 4253 & 298.26 & 6.13 \\
poi-40-21 & 39 & 229.6 & 6 & 4979 & 276.71 & 221.1 & 10 & 3798 & 298.97 & 3.7 \\
poi-40-22 & 39 & 255.8 & 11 & 4118 & 260.16 & 253.5 & 11 & 3520 & 233.42 & 0.9 \\
poi-40-23 & 39 & 330 & 15 & 5311 & 290.45 & 290.6 & 10 & 3269 & 243.26 & 11.94 \\
poi-40-24 & 39 & 279.8 & 10 & 4211 & 298.27 & 206.5 & 10 & 2900 & 191.85 & 26.2 \\
poi-40-25 & 39 & 285.1 & 9 & 4157 & 280.93 & 232.4 & 9 & 3956 & 272.85 & 18.48 \\
        \bottomrule
    \end{tabular*}
\end{table}

Until now, we have used Poikonen instances (up to 39 customers) from Set 4 to compare LKH Config1 against a heuristic, i.e., HGA. To further demonstrate the effectiveness of our solution framework on FSTSP, we apply LKH Config1 to the Murray instances from Set1 with 10 and 20 customers, for which BKS exact solutions are available from \cite{boccia2023new}. Note that the solution files for 10-cus instances can be obtained from \cite{dell2021benchmark}.
Following their problem setup, we incorporate the constraints of drone endurance and L/R operation times into LKH Config1, and execute it once within the limit of 1h. 

As shown in Tables~\ref{tab:Boccia10cusE20_1LKH}--\ref{tab:Boccia10cusE40_1LKH}, the number of drone-eligible customers is denoted as $N_D$, different from the total number of customers ($10$). The optimal makespan ($opt$), the number of drone-serviced customers (Dcus), and the corresponding computational time are sourced from \cite{boccia2023new}. Note that the solutions with zero Dcus are actually TSP tours. In Table~\ref{tab:Boccia10cusE20_1LKH} and Table~\ref{tab:Boccia10cusE40_1LKH} with 20 and 40 minutes endurance, the gap column shows $ \%gap = 100 \cdot (opt-\text{M}^1)/opt $. Only 3 instances are highlighted in bold since LKH could not find the optimal solutions within the same time limit of 1h.

\begin{table}[H]
    \centering
    % \tiny
    \scriptsize
    % \footnotesize  % Slightly smaller to better fit
    % \small
    % \normalsize  % (default size) 
    % \large
    \caption{FSTSP (10 cus) Comparison of Boccia and LKH (Config 1) with endurance $E=20$}
    \label{tab:Boccia10cusE20_1LKH}
    \begin{tabular*}{\textwidth}{@{\extracolsep{\fill}}l@{\;}r@{\;}r@{\;}r@{\;}r@{\;}r@{\;}r@{\;}r@{\;}r@{\;}r@{\;}r@{\;}r@{\;}r@{\;}r@{\;}r@{\;}r@{\;}r@{\;}r@{}}
    %%
    % \begin{tabularx}{\textwidth}{l r r r r r r r r r r}
    % \begin{tabularx}{\textwidth}{l c c c c c c c c c c}
    % \begin{tabularx}{\textwidth}{l *{10}{>{\centering\arraybackslash}X}}
        \toprule
        \multicolumn{2}{c}{} & \multicolumn{3}{c}{\textbf{Boccia}} & \multicolumn{3}{c}{\textbf{LKH (Config 1)}} & \multicolumn{1}{c}{} & \multicolumn{2}{c}{} & \multicolumn{3}{c}{\textbf{Boccia}} & \multicolumn{3}{c}{\textbf{LKH (Config 1)}} & \multicolumn{1}{c}{} \\
        
        \cmidrule(lr){3-5} \cmidrule(lr){6-8}
        \cmidrule(lr){12-14} \cmidrule(lr){15-17}
        
        \raisebox{1.5ex}[0pt][0pt]{Instance} & \raisebox{1.5ex}[0pt][0pt] {$N_D$} & $opt$ & Dcus & Time(s) & $\text{M}^1$ & Dcus & Time(s) & \raisebox{1.5ex}[0pt][0pt] {$\%gap$} & \raisebox{1.5ex}[0pt][0pt]{Instance} & \raisebox{1.5ex}[0pt][0pt] {$N_D$} & $opt$ & Dcus & Time(s) & $\text{M}^1$ & Dcus & Time(s) & \raisebox{1.5ex}[0pt][0pt] {$\%gap$} \\
        \midrule
37v1 & 9 & 57.45 & 0 & 0.96 & 57.6 & 0 & 0.6 & 0 & 40v7 & 8 & 49.9 & 2 & 1.15 & 50 & 2 & 6.15 & 0 \\
37v2 & 9 & 53.79 & 1 & 0.72 & 53.7 & 1 & 7.56 & 0 & 40v8 & 8 & 62.7 & 2 & 0.93 & 63 & 2 & 5.82 & 0 \\
37v3 & 9 & 54.66 & 0 & 0.41 & 54.7 & 0 & 6.49 & 0 & \textbf{40v9} & 8 & 42.53 & \textbf{4} & 1.25 & 43 & \textbf{3} & 11.84 & \textbf{-1.11} \\
37v4 & 9 & 67.46 & 0 & 0.33 & 67.7 & 0 & 0.35 & 0 & 40v10 & 8 & 43.08 & 4 & 2.05 & 43 & 4 & 279.36 & 0 \\
37v5 & 9 & 51.78 & 2 & 3.42 & 51.8 & 2 & 51.43 & 0 & 40v11 & 8 & 49.2 & 3 & 0.68 & 49.2 & 3 & 0.87 & 0 \\
37v6 & 9 & 48.6 & 2 & 2.65 & 48.8 & 2 & 7.72 & 0 & 40v12 & 8 & 62 & 3 & 0.62 & 62.2 & 3 & 9.7 & 0 \\
37v7 & 9 & 49.58 & 2 & 2.36 & 49.9 & 2 & 21.3 & 0 & 43v1 & 8 & 69.59 & 0 & 0.15 & 68.2 & 0 & 0.24 & 0 \\
37v8 & 9 & 62.38 & 2 & 1.99 & 62.9 & 2 & 3.84 & 0 & 43v2 & 8 & 72.15 & 0 & 0.18 & 70.6 & 0 & 0.35 & 0 \\
37v9 & 9 & 43.48 & 2 & 3.04 & 43.6 & 2 & 138.65 & 0 & 43v3 & 8 & 77.34 & 0 & 0.12 & 76.3 & 0 & 0.55 & 0 \\
37v10 & 9 & 41.91 & 3 & 2.03 & 41.9 & 3 & 8.32 & 0 & 43v4 & 8 & 90.14 & 0 & 0.34 & 89.3 & 0 & 0.53 & 0 \\
37v11 & 9 & 42.9 & 2 & 1.45 & 42.8 & 2 & 686.24 & 0 & 43v5 & 8 & 58.71 & 2 & 0.98 & 58.3 & 2 & 0.89 & 0 \\
37v12 & 9 & 56.85 & 3 & 1.44 & 56.9 & 3 & 1.55 & 0 & 43v6 & 8 & 59.09 & 2 & 1.54 & 58.6 & 2 & 3.34 & 0 \\
40v1 & 8 & 49.43 & 2 & 0.81 & 49.6 & 2 & 111.38 & 0 & 43v7 & 8 & 65.52 & 2 & 1.24 & 64.8 & 2 & 1.57 & 0 \\
40v2 & 8 & 51.71 & 1 & 0.84 & 51.9 & 1 & 0.97 & 0 & 43v8 & 8 & 84.81 & 1 & 1.66 & 84.4 & 1 & 2.59 & 0 \\
40v3 & 8 & 57.1 & 1 & 0.91 & 57.2 & 1 & 2.86 & 0 & 43v9 & 8 & 46.93 & 3 & 2.94 & 46.7 & 3 & 3.68 & 0 \\
40v4 & 8 & 69.9 & 1 & 0.8 & 70.2 & 1 & 6.59 & 0 & 43v10 & 8 & 47.93 & 3 & 2.69 & 47.4 & 3 & 30.32 & 0 \\
40v5 & 8 & 45.46 & 2 & 2.45 & 45.7 & 2 & 79.11 & 0 & 43v11 & 8 & 57.38 & 3 & 0.85 & 57 & 3 & 97.89 & 0 \\
40v6 & 8 & 44.51 & 2 & 1.31 & 44.7 & 2 & 33.72 & 0 & 43v12 & 8 & 69.2 & 3 & 0.78 & 68.8 & 3 & 28.02 & 0 \\
Average & & & & & & & & & & & & & 1.34 & & & 45.90 & \\
        \bottomrule
    \end{tabular*}
\end{table}

\begin{table}[H]
    \centering
    % \tiny
    \scriptsize
    % \footnotesize  % Slightly smaller to better fit
    % \small
    % \normalsize  % (default size) 
    % \large
    \caption{FSTSP (10 cus) Comparison of Boccia and LKH (Config 1) with endurance $E=40$}
    \label{tab:Boccia10cusE40_1LKH}
    \begin{tabular*}{\textwidth}{@{\extracolsep{\fill}}l@{\;}r@{\;}r@{\;}r@{\;}r@{\;}r@{\;}r@{\;}r@{\;}r@{\;}r@{\;}r@{\;}r@{\;}r@{\;}r@{\;}r@{\;}r@{\;}r@{\;}r@{}}
    %%
    % \begin{tabularx}{\textwidth}{l r r r r r r r r r r}
    % \begin{tabularx}{\textwidth}{l c c c c c c c c c c}
    % \begin{tabularx}{\textwidth}{l *{10}{>{\centering\arraybackslash}X}}
        \toprule
        \multicolumn{2}{c}{} & \multicolumn{3}{c}{\textbf{Boccia}} & \multicolumn{3}{c}{\textbf{LKH (Config 1)}} & \multicolumn{1}{c}{} & \multicolumn{2}{c}{} & \multicolumn{3}{c}{\textbf{Boccia}} & \multicolumn{3}{c}{\textbf{LKH (Config 1)}} & \multicolumn{1}{c}{} \\
        
        \cmidrule(lr){3-5} \cmidrule(lr){6-8}
        \cmidrule(lr){12-14} \cmidrule(lr){15-17}
        
        \raisebox{1.5ex}[0pt][0pt]{Instance} & \raisebox{1.5ex}[0pt][0pt] {$N_D$} & $opt$ & Dcus & Time(s) & $\text{M}^1$ & Dcus & Time(s) & \raisebox{1.5ex}[0pt][0pt] {$\%gap$} & \raisebox{1.5ex}[0pt][0pt]{Instance} & \raisebox{1.5ex}[0pt][0pt] {$N_D$} & $opt$ & Dcus & Time(s) & $\text{M}^1$ & Dcus & Time(s) & \raisebox{1.5ex}[0pt][0pt] {$\%gap$} \\
        \midrule
37v1 & 9 & 50.57 & 1 & 2.2 & 50.5 & 1 & 10.12 & 0 & 40v7 & 8 & 49.23 & 3 & 2.35 & 49.3 & 3 & 2798.85 & 0 \\
37v2 & 9 & 47.31 & 1 & 1.43 & 47.1 & 1 & 9.01 & 0 & \textbf{40v8} & 8 & 62.03 & \textbf{3} & 1.83 & 62.6 & \textbf{2} & 13.98 & \textbf{-0.92} \\
37v3 & 9 & 53.69 & 1 & 2.38 & 53.8 & 1 & 1.14 & 0 & 40v9 & 8 & 42.53 & 4 & 1.25 & 42.6 & 4 & 920.73 & 0 \\
37v4 & 9 & 67.46 & 0 & 2.57 & 67.7 & 0 & 0.51 & 0 & 40v10 & 8 & 43.08 & 4 & 2.34 & 43 & 4 & 4.65 & 0 \\
37v5 & 9 & 45.84 & 2 & 2.75 & 45.8 & 2 & 3.44 & 0 & 40v11 & 8 & 49.2 & 3 & 0.8 & 49.2 & 3 & 2.51 & 0 \\
37v6 & 9 & 44.6 & 2 & 3.17 & 44.4 & 2 & 53.78 & 0 & 40v12 & 8 & 62 & 3 & 0.62 & 62.2 & 3 & 3.37 & 0 \\
37v7 & 9 & 47.62 & 2 & 3.3 & 47.6 & 1 & 1.94 & 0 & 43v1 & 8 & 57.01 & 1 & 2.35 & 55.3 & 1 & 3.57 & 0 \\
\textbf{37v8} & 9 & 60.42 & \textbf{2} & 2.36 & 60.6 & \textbf{1} & 26.21 & \textbf{-0.3} & 43v2 & 8 & 58.05 & 1 & 2.24 & 56.7 & 1 & 1.98 & 0 \\
37v9 & 9 & 42.42 & 2 & 2.08 & 42.3 & 2 & 87.51 & 0 & 43v3 & 8 & 69.43 & 2 & 2.33 & 68.6 & 2 & 201.1 & 0 \\
37v10 & 9 & 41.91 & 3 & 2.22 & 41.9 & 3 & 1.55 & 0 & 43v4 & 8 & 83.7 & 1 & 2.82 & 83.6 & 1 & 9.23 & 0 \\
37v11 & 9 & 42.9 & 2 & 2.44 & 42.8 & 2 & 0.21 & 0 & 43v5 & 8 & 52.09 & 2 & 2.43 & 51.9 & 2 & 10.69 & 0 \\
37v12 & 9 & 55.7 & 2 & 2.2 & 55.8 & 2 & 5.02 & 0 & 43v6 & 8 & 52.33 & 2 & 2.78 & 52 & 2 & 1.7 & 0 \\
40v1 & 8 & 46.89 & 1 & 3.16 & 47.1 & 1 & 2.99 & 0 & 43v7 & 8 & 61.88 & 3 & 2.62 & 61.5 & 2 & 0.41 & 0 \\
40v2 & 8 & 46.42 & 1 & 2.31 & 46.1 & 1 & 0.36 & 0 & 43v8 & 8 & 73.73 & 3 & 2.86 & 73.6 & 3 & 163.13 & 0 \\
40v3 & 8 & 53.93 & 1 & 3.2 & 53.9 & 1 & 8.81 & 0 & 43v9 & 8 & 46.93 & 3 & 4.74 & 46.7 & 3 & 3.24 & 0 \\
40v4 & 8 & 68.4 & 1 & 3.72 & 68.8 & 1 & 13.83 & 0 & 43v10 & 8 & 47.93 & 3 & 3.61 & 47.4 & 3 & 23.4 & 0 \\
40v5 & 8 & 43.53 & 2 & 2.62 & 43.9 & 2 & 6.84 & 0 & 43v11 & 8 & 56.4 & 3 & 1.82 & 55.9 & 3 & 859.72 & 0 \\
40v6 & 8 & 44.08 & 2 & 1.92 & 44.3 & 2 & 61.41 & 0 & 43v12 & 8 & 69.2 & 3 & 1.17 & 68.8 & 3 & 40.9 & 0 \\
Average & & & & & & & & & & & & & 2.42 & & & 148.83 & \\
        \bottomrule
    \end{tabular*}
\end{table}

As the problem scale increases to 20 customers (see Tables~\ref{tab:Boccia20cusE20_1LKH}--\ref{tab:Boccia20cusE40_1LKH}), only 6 instances highlighted in bold were not solved to optimality by \cite{boccia2023new}. Notably, for instance 5246 under a 20 minutes endurance constraint in Table~\ref{tab:Boccia20cusE20_1LKH}, the BKS upper bound is improved by LKH from 74.19 to 72 using just 7.57s. Furthermore, the gaps $ \%gap = 100 \cdot (opt-\text{M}^1)/opt $ are nearly zero, indicating that LKH can match the BKS optimal solutions across most cases. 

When increasing the endurance from 20 to 40 minutes in Table~\ref{tab:Boccia20cusE40_1LKH}, the number of unsolved instances (bold) increases to 5, due to solution space expansion. Correspondingly, the average time in \cite{boccia2023new} rises significantly from 77.36s to 535.44s, while LKH only grows from 235.69s to 369.81s, demonstrating its competitiveness under more complex conditions. Note that LKH enhances instance 5252, reducing the BKS upper bound from 85.33 to 84.8 in 276.57s.

\begin{table}[H]
    \centering
    \tiny
    % \scriptsize
    % \footnotesize  % Slightly smaller to better fit
    % \small
    % \normalsize  % (default size) 
    % \large
    \caption{FSTSP (20 cus) Comparison of Boccia and LKH (Config 1) with endurance $E=20$}
    \label{tab:Boccia20cusE20_1LKH}
    \begin{tabular*}{\textwidth}{@{\extracolsep{\fill}}l@{\;}r@{\;}r@{\;}r@{\;}r@{\;}r@{\;}r@{\;}r@{\;}r@{\;}r@{\;}r@{\;}r@{\;}r@{\;}r@{\;}r@{\;}r@{\;}r@{\;}r@{\;}r@{\;}r@{\;}r@{}}
    %%
    % \begin{tabularx}{\textwidth}{l r r r r r r r r r r}
    % \begin{tabularx}{\textwidth}{l c c c c c c c c c c}
    % \begin{tabularx}{\textwidth}{l *{10}{>{\centering\arraybackslash}X}}
        \toprule
        \multicolumn{2}{c}{} & \multicolumn{2}{c}{\textbf{Boccia}} & \multicolumn{2}{c}{\textbf{LKH}} & \multicolumn{1}{c}{} & \multicolumn{2}{c}{} & \multicolumn{2}{c}{\textbf{Boccia}} & \multicolumn{2}{c}{\textbf{LKH}} & \multicolumn{1}{c}{} & \multicolumn{2}{c}{} & \multicolumn{2}{c}{\textbf{Boccia}} & \multicolumn{2}{c}{\textbf{LKH}} & \multicolumn{1}{c}{}\\
        
        \cmidrule(lr){3-4} \cmidrule(lr){5-6} \cmidrule(lr){10-11} \cmidrule(lr){12-13}
        \cmidrule(lr){17-18} \cmidrule(lr){19-20}
        
        \raisebox{1.5ex}[0pt][0pt]{Ins} & \raisebox{1.5ex}[0pt][0pt] {$N_D$} & $UB$ & Time & $\text{M}^1$ & Time & \raisebox{1.5ex}[0pt][0pt] {$\%gap$} & \raisebox{1.5ex}[0pt][0pt]{Ins} & \raisebox{1.5ex}[0pt][0pt] {$N_D$} & $UB$ & Time & $\text{M}^1$ & Time & \raisebox{1.5ex}[0pt][0pt] {$\%gap$} & \raisebox{1.5ex}[0pt][0pt]{Ins} & \raisebox{1.5ex}[0pt][0pt] {$N_D$} & $UB$ & Time & $\text{M}^1$ & Time & \raisebox{1.5ex}[0pt][0pt] {$\%gap$} \\
        \midrule
4847 & 17 & 267.05 & 5.39 & 267.7 & 175.42 & -0.24 & 5025 & 17 & 131.43 & 9.87 & 131.2 & 316.14 & 0 & 5154 & 17 & 121.9 & 49.95 & 120.1 & 91.22 & 0 \\
4849 & 17 & 248.3 & 0.84 & 246.7 & 17.79 & 0 & 5027 & 17 & 114.31 & 14.32 & 113.3 & 37.78 & 0 & 5156 & 16 & 124.46 & 21.18 & 124.6 & 289.74 & -0.11 \\
4853 & 17 & 232.87 & 1.01 & 233.6 & 8.95 & -0.31 & 5030 & 16 & 116.1 & 44.69 & 115.7 & 24.64 & 0 & 5159 & 16 & 145.79 & 18.46 & 145.9 & 16.67 & -0.08 \\
4856 & 18 & 253.33 & 1.47 & 252.9 & 26.28 & 0 & 5032 & 17 & 117.55 & 13.28 & 117.3 & 65.65 & 0 & 5201 & 16 & 148.02 & 9.25 & 147.9 & 33.24 & 0 \\
4858 & 18 & 240.63 & 1.39 & 238.2 & 62.34 & 0 & 5034 & 18 & 105.1 & 42.98 & 103.9 & 968.77 & 0 & 5203 & 18 & 138.59 & 21.24 & 138.2 & 310.57 & 0 \\
4902 & 17 & 242.32 & 1.25 & 242 & 236.89 & 0 & 5036 & 17 & 124.33 & 22.99 & 123.7 & 18.85 & 0 & 5205 & 18 & 134.78 & 49.77 & 133.9 & 70.53 & 0 \\
4907 & 18 & 239.28 & 0.91 & 240.1 & 318.59 & -0.34 & 5039 & 18 & 130.91 & 7.5 & 131.8 & 20.06 & -0.68 & 5207 & 18 & 121.47 & 23.18 & 121.8 & 429.93 & -0.27 \\
4909 & 17 & 222.88 & 1.15 & 222.5 & 50.46 & 0 & 5041 & 17 & 125.32 & 10.83 & 124.8 & 4.21 & 0 & 5209 & 16 & 135.92 & 15.72 & 135.4 & 309.78 & 0 \\
4912 & 18 & 267.62 & 0.87 & 267.1 & 5.08 & 0 & 5044 & 16 & 120.77 & 25.8 & 120.2 & 104.68 & 0 & 5212 & 17 & 137.67 & 18.77 & 137.5 & 596.79 & 0 \\
4915 & 17 & 259.39 & 0.66 & 259.1 & 5.3 & 0 & 5047 & 17 & 112.85 & 9.61 & 111.2 & 32.83 & 0 & 5214 & 18 & 126.25 & 27.93 & 124.9 & 29.51 & 0 \\
4917 & 17 & 173.97 & 3.6 & 172.9 & 29.8 & 0 & 5049 & 17 & 197.76 & 3.31 & 196.7 & 56.47 & 0 & 5216 & 18 & 101.07 & 46.56 & 98.8 & 1177.24 & 0 \\
4920 & 17 & 170.05 & 4.92 & 169.8 & 9.48 & 0 & 5051 & 17 & 180.62 & 9.92 & 180.8 & 27.55 & -0.1 & 5218 & 17 & 115.62 & 27.93 & 114.4 & 11.23 & 0 \\
4922 & 17 & 169.6 & 2.66 & 168.7 & 12.08 & 0 & 5053 & 16 & 176.51 & 2.01 & 175.3 & 54.67 & 0 & 5220 & 16 & 119.04 & 28.49 & 118.9 & 991.79 & 0 \\
4924 & 16 & 159.57 & 2.56 & 158.6 & 86.51 & 0 & 5055 & 16 & 177.29 & 4.23 & 177.3 & 10.44 & -0.01 & 5223 & 18 & 94.59 & 42.48 & 93.4 & 75.39 & 0 \\
4926 & 18 & 155.98 & 13.54 & 155.6 & 6.45 & 0 & 5057 & 16 & 180.7 & 6.02 & 181.1 & 5.79 & -0.22 & 5225 & 17 & 129.72 & 12.36 & 127.5 & 121.17 & 0 \\
4928 & 18 & 166 & 18.64 & 166.7 & 70.65 & -0.42 & 5059 & 16 & 150.82 & 5.64 & 149.7 & 27.48 & 0 & 5227 & 17 & 116.19 & 75.77 & 115.8 & 115.31 & 0 \\
4931 & 17 & 172.49 & 14.91 & 173.1 & 69.88 & -0.35 & 5102 & 16 & 165.49 & 5.24 & 165.7 & 593.29 & -0.13 & 5229 & 17 & 92.87 & 70.17 & 92.1 & 128.57 & 0 \\
4933 & 18 & 159.39 & 13.83 & 159.3 & 265.12 & 0 & 5104 & 16 & 181.61 & 4.93 & 181.5 & 16.27 & 0 & 5231 & 17 & 98.36 & 38.46 & 99.6 & 23.78 & -1.26 \\
4935 & 17 & 176.69 & 18.86 & 176.4 & 75.37 & 0 & 5106 & 16 & 158.49 & 12.78 & 157.1 & 10.69 & 0 & 5233 & 17 & 111.62 & 20.99 & 111.5 & 148.6 & 0 \\
4937 & 17 & 172.5 & 5.51 & 170.5 & 24.29 & 0 & 5108 & 16 & 172.12 & 1.96 & 170 & 17.35 & 0 & 5235 & 17 & 118.89 & 12.2 & 119.1 & 1691.07 & -0.18 \\
4939 & 18 & 201.03 & 5.34 & 201.2 & 1021.09 & -0.08 & 5110 & 18 & 135.43 & 17.93 & 135.2 & 6.03 & 0 & 5238 & 18 & 78.93 & 562.17 & 76.8 & 932.52 & 0 \\
4941 & 16 & 253.08 & 1.27 & 252.6 & 86.5 & 0 & 5112 & 17 & 131.23 & 13.2 & 130.9 & 526.95 & 0 & 5240 & 18 & 84.51 & 183.79 & 84 & 264.95 & 0 \\
4944 & 18 & 247.03 & 2.98 & 247.3 & 1531.87 & -0.11 & 5115 & 18 & 127.39 & 19.19 & 126.7 & 16.86 & 0 & 5242 & 17 & 85.65 & 10.29 & 85 & 494.92 & 0 \\
4946 & 16 & 237.21 & 1.36 & 236.9 & 3.19 & 0 & 5117 & 17 & 130.36 & 46.49 & 129.2 & 252.58 & 0 & 5244 & 18 & 86.81 & 33.92 & 83.7 & 560.73 & 0 \\
4948 & 17 & 258.06 & 1.38 & 257.3 & 4.92 & 0 & 5119 & 17 & 118.97 & 28.65 & 117.9 & 31.07 & 0 & \textbf{5246} & 18 & \textbf{74.19} & \textbf{3600} & \textbf{72} & \textbf{7.57} & \textbf{2.95} \\
4950 & 17 & 240.99 & 1.44 & 239.2 & 266.81 & 0 & 5121 & 17 & 131.19 & 25.25 & 131.8 & 385.22 & -0.46 & 5248 & 17 & 83.1 & 1760.19 & 81 & 14.81 & 0 \\
4952 & 18 & 218.09 & 4.77 & 218.5 & 47.13 & -0.19 & 5123 & 17 & 121.95 & 16.46 & 124.3 & 487.63 & -1.93 & 5250 & 17 & 81.93 & 42.83 & 81.3 & 72.03 & 0 \\
4954 & 16 & 261.06 & 1.4 & 261.4 & 189.58 & -0.13 & 5125 & 18 & 130.96 & 13.52 & 128.9 & 20.13 & 0 & 5252 & 17 & 86.15 & 110.9 & 87.9 & 409.72 & -2.03 \\
4957 & 17 & 252.28 & 1.2 & 253.6 & 3.64 & -0.52 & 5127 & 16 & 132.24 & 7.97 & 131.9 & 179.89 & 0 & 5255 & 17 & 82.5 & 89.69 & 81.3 & 101.91 & 0 \\
4959 & 17 & 249.92 & 1.11 & 248.6 & 291.25 & 0 & 5130 & 18 & 126.49 & 45.19 & 125.6 & 199.14 & 0 & 5257 & 18 & 78.27 & 329.14 & 78.5 & 687.49 & -0.29 \\
5001 & 17 & 115.5 & 17.58 & 115.5 & 381.81 & 0 & 5132 & 16 & 106.36 & 50.9 & 106.1 & 12.06 & 0 & 5306 & 17 & 100.46 & 19.14 & 101.3 & 22.46 & -0.84 \\
5003 & 17 & 172.73 & 8.03 & 170.8 & 991.3 & 0 & 5134 & 17 & 102.57 & 21.86 & 101.5 & 28.87 & 0 & 5310 & 16 & 92.41 & 47.1 & 91.7 & 66.86 & 0 \\
5006 & 17 & 155.39 & 5.9 & 154.8 & 1507.83 & 0 & 5136 & 16 & 98.91 & 32.02 & 98.2 & 46.07 & 0 & 5312 & 17 & 83.59 & 20.19 & 82.2 & 10.88 & 0 \\
5008 & 18 & 159.74 & 5.87 & 160.4 & 46.47 & -0.41 & 5138 & 16 & 91.61 & 73.1 & 90.7 & 338.58 & 0 & 5321 & 17 & 101.49 & 11.89 & 102 & 136.83 & -0.5 \\
5010 & 18 & 145.48 & 15 & 146.1 & 208.32 & -0.43 & 5141 & 17 & 95.53 & 68.49 & 96.1 & 236.62 & -0.6 & 5324 & 18 & 101.55 & 77.96 & 102.8 & 24.95 & -1.23 \\
5012 & 16 & 171.38 & 8.02 & 171.9 & 54.07 & -0.3 & 5143 & 17 & 97.69 & 26.96 & 96.6 & 682.85 & 0 & 5330 & 16 & 118.45 & 7.55 & 118.1 & 31.22 & 0 \\
5015 & 17 & 172.67 & 3.43 & 176.1 & 2.57 & -1.99 & 5145 & 16 & 95.09 & 43.67 & 92.7 & 58.78 & 0 & 5334 & 16 & 102.02 & 26.4 & 103.3 & 730.5 & -1.25 \\
5017 & 16 & 166.47 & 8.86 & 166 & 4.55 & 0 & 5148 & 16 & 90.58 & 52.35 & 88.9 & 102.82 & 0 & 5336 & 18 & 104.46 & 12.53 & 102.8 & 103.39 & 0 \\
5020 & 17 & 156.92 & 3.46 & 156.5 & 40.52 & 0 & 5150 & 18 & 82.19 & 172.37 & 81.2 & 545.24 & 0 & 5345 & 17 & 114.19 & 7.88 & 112.6 & 23.73 & 0 \\
5022 & 17 & 146.21 & 13.71 & 145.1 & 1101.22 & 0 & 5152 & 17 & 90.58 & 384.83 & 89.1 & 31.61 & 0 & 5351 & 17 & 115.1 & 53.96 & 116.2 & 979.35 & -0.96 \\
Average & & & & & & & & & & & & & & & & & 77.36 & & 235.69 & \\
        \bottomrule
    \end{tabular*}
    \begin{tablenotes}\footnotesize
    \item Note: Out of 120 instances, 119 non-bold ones are solved to optimality by \cite{boccia2023new}.
    \end{tablenotes}
\end{table}

\begin{table}[H]
    \centering
    \tiny
    % \scriptsize
    % \footnotesize  % Slightly smaller to better fit
    % \small
    % \normalsize  % (default size) 
    % \large
    \caption{FSTSP (20 cus) Comparison of Boccia and LKH (Config 1) with endurance $E=40$}
    \label{tab:Boccia20cusE40_1LKH}
    \begin{tabular*}{\textwidth}{@{\extracolsep{\fill}}l@{\;}r@{\;}r@{\;}r@{\;}r@{\;}r@{\;}r@{\;}r@{\;}r@{\;}r@{\;}r@{\;}r@{\;}r@{\;}r@{\;}r@{\;}r@{\;}r@{\;}r@{\;}r@{\;}r@{\;}r@{}}
    %%
    % \begin{tabularx}{\textwidth}{l r r r r r r r r r r}
    % \begin{tabularx}{\textwidth}{l c c c c c c c c c c}
    % \begin{tabularx}{\textwidth}{l *{10}{>{\centering\arraybackslash}X}}
        \toprule
        \multicolumn{2}{c}{} & \multicolumn{2}{c}{\textbf{Boccia}} & \multicolumn{2}{c}{\textbf{LKH}} & \multicolumn{1}{c}{} & \multicolumn{2}{c}{} & \multicolumn{2}{c}{\textbf{Boccia}} & \multicolumn{2}{c}{\textbf{LKH}} & \multicolumn{1}{c}{} & \multicolumn{2}{c}{} & \multicolumn{2}{c}{\textbf{Boccia}} & \multicolumn{2}{c}{\textbf{LKH}} & \multicolumn{1}{c}{}\\
        
        \cmidrule(lr){3-4} \cmidrule(lr){5-6} \cmidrule(lr){10-11} \cmidrule(lr){12-13}
        \cmidrule(lr){17-18} \cmidrule(lr){19-20}
        
        \raisebox{1.5ex}[0pt][0pt]{Ins} & \raisebox{1.5ex}[0pt][0pt] {$N_D$} & $UB$ & Time & $\text{M}^1$ & Time & \raisebox{1.5ex}[0pt][0pt] {$\%gap$} & \raisebox{1.5ex}[0pt][0pt]{Ins} & \raisebox{1.5ex}[0pt][0pt] {$N_D$} & $UB$ & Time & $\text{M}^1$ & Time & \raisebox{1.5ex}[0pt][0pt] {$\%gap$} & \raisebox{1.5ex}[0pt][0pt]{Ins} & \raisebox{1.5ex}[0pt][0pt] {$N_D$} & $UB$ & Time & $\text{M}^1$ & Time & \raisebox{1.5ex}[0pt][0pt] {$\%gap$} \\
        \midrule
4847 & 17 & 255.6 & 51.14 & 257.1 & 83.83 & -0.59 & 5025 & 17 & 118.43 & 533.62 & 125.1 & 359.2 & -5.63 & 5154 & 17 & 104.28 & 296.76 & 106.4 & 337.29 & -2.03 \\
4849 & 17 & 225.15 & 108.81 & 228.8 & 62.6 & -1.62 & 5027 & 17 & 111.81 & 333.1 & 111.1 & 525.04 & 0 & 5156 & 16 & 107.12 & 120.46 & 112.2 & 109.49 & -4.74 \\
4853 & 17 & 211.38 & 115.55 & 214.7 & 780.31 & -1.57 & 5030 & 16 & 100.64 & 3490.25 & 102.4 & 1706.78 & -1.75 & 5159 & 16 & 120.02 & 158.89 & 120 & 41.13 & 0 \\
4856 & 18 & 234.63 & 113.71 & 234 & 166.14 & 0 & 5032 & 17 & 103.06 & 204.23 & 107.1 & 25.82 & -3.92 & 5201 & 16 & 140.3 & 813.36 & 141.4 & 54.18 & -0.78 \\
4858 & 18 & 215.63 & 52.83 & 217 & 187.8 & -0.64 & 5034 & 18 & 102.16 & 168.37 & 101.5 & 134.71 & 0 & 5203 & 18 & 124.06 & 359.09 & 122.8 & 206.88 & 0 \\
4902 & 17 & 223.84 & 48.92 & 224.6 & 180.43 & -0.34 & 5036 & 17 & 112.1 & 88.98 & 110.8 & 68.24 & 0 & 5205 & 18 & 117.68 & 329.91 & 120.7 & 701.28 & -2.57 \\
4907 & 18 & 192.87 & 41.95 & 192.6 & 594.78 & 0 & \textbf{5039} & 18 & \textbf{115.85} & \textbf{3600} & \textbf{116} & \textbf{704.11} & \textbf{-0.13} & 5207 & 18 & 113.9 & 318.79 & 114.5 & 498.99 & -0.53 \\
4909 & 17 & 216.04 & 48.41 & 216.9 & 466.16 & -0.4 & 5041 & 17 & 113.88 & 276.09 & 114.2 & 74.94 & -0.28 & \textbf{5209} & 16 & \textbf{121.24} & \textbf{3600} & \textbf{123.8} & \textbf{21.54} & \textbf{-2.11} \\
4912 & 18 & 238.04 & 51.64 & 244.1 & 63.48 & -2.55 & 5044 & 16 & 115.09 & 510.86 & 115.9 & 140.7 & -0.7 & 5212 & 17 & 133.06 & 443.64 & 135.1 & 415.61 & -1.53 \\
4915 & 17 & 229.03 & 39.29 & 229.7 & 249.92 & -0.29 & 5047 & 17 & 103.86 & 354.8 & 103.6 & 244.8 & 0 & 5214 & 18 & 121.85 & 282.34 & 121.9 & 424.44 & -0.04 \\
4917 & 17 & 165.31 & 29.56 & 165.6 & 178.43 & -0.18 & 5049 & 17 & 180.71 & 1065.02 & 181.2 & 482.72 & -0.27 & 5216 & 18 & 91.94 & 146.8 & 91.8 & 14.63 & 0 \\
4920 & 17 & 157.03 & 131.76 & 159.2 & 248.2 & -1.38 & 5051 & 17 & 162.52 & 803.47 & 164.5 & 1005.8 & -1.22 & 5218 & 17 & 96.06 & 161.47 & 95.7 & 1424.18 & 0 \\
4922 & 17 & 164.1 & 133.3 & 162.5 & 15.6 & 0 & 5053 & 16 & 145.42 & 115.44 & 145.3 & 1741.75 & 0 & 5220 & 16 & 114.93 & 171.7 & 118.4 & 158.91 & -3.02 \\
4924 & 16 & 158.98 & 93.41 & 157.6 & 31.82 & 0 & 5055 & 16 & 171.64 & 77.22 & 172.6 & 116.77 & -0.56 & 5223 & 18 & 91.22 & 503.58 & 92.4 & 13.64 & -1.29 \\
4926 & 18 & 151.19 & 305.57 & 152.1 & 62.21 & -0.6 & 5057 & 16 & 172.88 & 479.58 & 173.3 & 299.45 & -0.24 & 5225 & 17 & 123.58 & 84.26 & 121.3 & 919.71 & 0 \\
\textbf{4928} & 18 & \textbf{146.76} & \textbf{3600} & \textbf{149.9} & \textbf{15.77} & \textbf{-2.14} & 5059 & 16 & 134.83 & 261.99 & 142.8 & 42.06 & -5.91 & 5227 & 17 & 101.74 & 149.88 & 100.3 & 892.4 & 0 \\
4931 & 17 & 165.22 & 162.1 & 165.8 & 340.21 & -0.35 & 5102 & 16 & 163.77 & 124.13 & 163.2 & 283.17 & 0 & 5229 & 17 & 88.67 & 70.36 & 91.2 & 631.84 & -2.85 \\
4933 & 18 & 155.13 & 401.8 & 156.6 & 77.44 & -0.95 & 5104 & 16 & 176.7 & 104.39 & 177.5 & 183.65 & -0.45 & 5231 & 17 & 94.02 & 1583.46 & 96.2 & 179.3 & -2.32 \\
4935 & 17 & 155.77 & 1860.91 & 157 & 318.47 & -0.79 & 5106 & 16 & 141.77 & 196.89 & 144.1 & 285.06 & -1.64 & 5233 & 17 & 107.9 & 138.11 & 108.7 & 786.31 & -0.74 \\
4937 & 17 & 159.18 & 490.75 & 157.3 & 11.07 & 0 & 5108 & 16 & 162.21 & 935.15 & 161.8 & 14.07 & 0 & 5235 & 17 & 98.96 & 357.6 & 101.2 & 1392.86 & -2.26 \\
4939 & 18 & 179.1 & 2986.07 & 179.3 & 77.05 & -0.11 & 5110 & 18 & 126.2 & 1322.61 & 127.5 & 22.08 & -1.03 & 5238 & 18 & 77.99 & 281.27 & 76.8 & 43.34 & 0 \\
4941 & 16 & 241.9 & 26.95 & 242.7 & 720.75 & -0.33 & 5112 & 17 & 123.13 & 287.31 & 127 & 1261.33 & -3.14 & 5240 & 18 & 81.02 & 134.01 & 83.7 & 452.76 & -3.31 \\
4944 & 18 & 233.75 & 54.47 & 240.6 & 14.68 & -2.93 & 5115 & 18 & 123.08 & 3115.85 & 126.2 & 8.32 & -2.53 & 5242 & 17 & 84.52 & 81.7 & 84 & 173.24 & 0 \\
4946 & 16 & 220.39 & 75.87 & 220.9 & 300.43 & -0.23 & 5117 & 17 & 118.87 & 505.66 & 123.5 & 158.18 & -3.9 & 5244 & 18 & 85.32 & 33.78 & 81.9 & 618.22 & 0 \\
4948 & 17 & 244.13 & 18.35 & 245.9 & 55.65 & -0.73 & 5119 & 17 & 113 & 223.39 & 111.3 & 1145.97 & 0 & 5246 & 18 & 63.25 & 125.96 & 61 & 1002.11 & 0 \\
4950 & 17 & 225.61 & 28.95 & 225 & 271.82 & 0 & 5121 & 17 & 112.03 & 144.85 & 118.4 & 63.03 & -5.69 & 5248 & 17 & 77.55 & 199.2 & 80.7 & 7.55 & -4.06 \\
4952 & 18 & 210.67 & 346.84 & 211.9 & 5.05 & -0.58 & 5123 & 17 & 113.75 & 229.14 & 114.5 & 105.38 & -0.66 & 5250 & 17 & 81.73 & 2120.91 & 81.8 & 516.24 & -0.09 \\
4954 & 16 & 239.67 & 40.83 & 244.5 & 72.9 & -2.02 & 5125 & 18 & 118.86 & 592.08 & 121.7 & 1318.6 & -2.39 & \textbf{5252} & 17 & \textbf{85.33} & \textbf{3600} & \textbf{84.8} & \textbf{276.57} & \textbf{0.62} \\
4957 & 17 & 225.85 & 88.06 & 230.7 & 508.6 & -2.15 & 5127 & 16 & 124.57 & 219.7 & 124.9 & 128.83 & -0.26 & 5255 & 17 & 71.13 & 345.29 & 72.3 & 370.54 & -1.64 \\
4959 & 17 & 231.95 & 36.42 & 231.9 & 12.69 & 0 & 5130 & 18 & 119.01 & 405.15 & 120 & 341.95 & -0.83 & 5257 & 18 & 75.17 & 753.53 & 77.6 & 2316.02 & -3.23 \\
5001 & 17 & 114.09 & 129.04 & 114.7 & 711.31 & -0.53 & 5132 & 16 & 102.77 & 471.4 & 104.2 & 17.23 & -1.39 & 5306 & 17 & 72.37 & 353.28 & 75.4 & 120.78 & -4.19 \\
5003 & 17 & 162.39 & 40.76 & 168.2 & 14.91 & -3.58 & 5134 & 17 & 101.24 & 117.72 & 100.5 & 133.6 & 0 & 5310 & 16 & 88.33 & 2028.71 & 90 & 525.34 & -1.89 \\
5006 & 17 & 135.92 & 94.63 & 135.6 & 532.44 & 0 & 5136 & 16 & 93.24 & 111.43 & 92.4 & 9.79 & 0 & 5312 & 17 & 82.2 & 90.37 & 81.8 & 19.92 & 0 \\
5008 & 18 & 152.11 & 249.34 & 155.4 & 1806.4 & -2.16 & 5138 & 16 & 87.53 & 645.75 & 90.7 & 104.33 & -3.62 & 5321 & 17 & 79.6 & 92.65 & 83.4 & 1495.17 & -4.77 \\
5010 & 18 & 133.85 & 377.62 & 137.9 & 88.65 & -3.03 & \textbf{5141} & 17 & \textbf{92.44} & \textbf{592.25} & \textbf{93.6} & \textbf{185.63} & \textbf{-1.25} & 5324 & 18 & 79.75 & 1108.16 & 79.2 & 13.06 & 0 \\
5012 & 16 & 154.34 & 232.78 & 155.9 & 274.76 & -1.01 & 5143 & 17 & 93.87 & 463.49 & 96.5 & 13.24 & -2.8 & 5330 & 16 & 96.35 & 632.45 & 98.1 & 31.69 & -1.82 \\
5015 & 17 & 154.38 & 300.89 & 155.8 & 1041.65 & -0.92 & 5145 & 16 & 88.13 & 71.65 & 92.7 & 14.32 & -5.19 & 5334 & 16 & 88.89 & 98.94 & 90.9 & 65.88 & -2.26 \\
5017 & 16 & 156.18 & 189.4 & 157.7 & 92.78 & -0.97 & 5148 & 16 & 865.29 & 128.27 & 86.5 & 454.36 & 0 & 5336 & 18 & 86.39 & 3080.09 & 89 & 127.87 & -3.02 \\
5020 & 17 & 141.92 & 146.23 & 144.2 & 696.24 & -1.61 & 5150 & 18 & 77.42 & 160.82 & 77.8 & 564.99 & -0.49 & 5345 & 17 & 99.63 & 83.02 & 102.6 & 75.05 & -2.98 \\
5022 & 17 & 128.16 & 1657.48 & 129.8 & 896.18 & -1.28 & 5152 & 17 & 79.17 & 163 & 82.2 & 17.26 & -3.83 & 5351 & 17 & 90.1 & 221.93 & 91.9 & 63.78 & -2 \\
Average & & & & & & & & & & & & & & & & & 535.44 & & 369.81 & \\
        \bottomrule
    \end{tabular*}
    \begin{tablenotes}\footnotesize
    \item Note: Out of 120 instances, 115 non-bold ones are solved to optimality by \cite{boccia2023new}.
    \end{tablenotes}
\end{table}

In summary, with a single 5min execution shown in Tables~\ref{tab:fstsp_1LKH_small}--\ref{tab:fstsp_1LKH_medium}, LKH Config1 derives optimal results on small Poikonen instances, outperforms the HGA algorithm (five 5min runs) from \cite{ha2020hybrid} on instances up to 19 customers and most 29-customer scenarios. Under an extended runtime of 4h in Table~\ref{tab:fstsp_1LKH_longer}, LKH surpasses HGA across all instances up to 39 customers except one instance poi-40-17. Moreover, when executed once within 1h for Murray instances Set1 reported in \cite{boccia2023new}, LKH generates the BKS solutions for 69 out of 72 instances with 10 customers and for the majority of 20-customer instances. Importantly, it improves the BKS upper bound for two of the six 20-customer instances that were not solved to optimality.

%85/120, 31/120: 116 out of 240 instances with 20 customers

\subsection{TSP-mD}
In this subsection, we continue to use the 125 min-makespan Poikonen instances of 6/9/19/29/39 customers to compare LKH Config1 against ALNS, a well-designed heuristic for TSP-mD. ALNS is executed five times, and each is constrained to 5min. The best result is chosen for comparison with LKH, which is run \textbf{once} within 5min. Moreover, we assume that there are \textbf{5 drones} available at the depot. This parameter is read from the input file \texttt{.drone}. Additionally, we compare LKH Config1 and Config2 on those 125 Poikonen instances.

Before analyzing the comparison tables: Tables~\ref{tab:tspmd_1LKH_small}--\ref{tab:tspmd_1LKH_longer} (ALNS vs. LKH Config1), Tables~\ref{tab:tspmd_2LKH_small}--\ref{tab:tspmd_2LKH_medium} (LKH Config1 vs. Config2), we first highlight several representative instances below.

\subsubsection{Instances}
\begin{figure}[H]
    \centering
    \setlength{\tabcolsep}{1pt} % Space between columns
    \begin{tabular}{ccc}
        \includegraphics[width=0.33\textwidth]{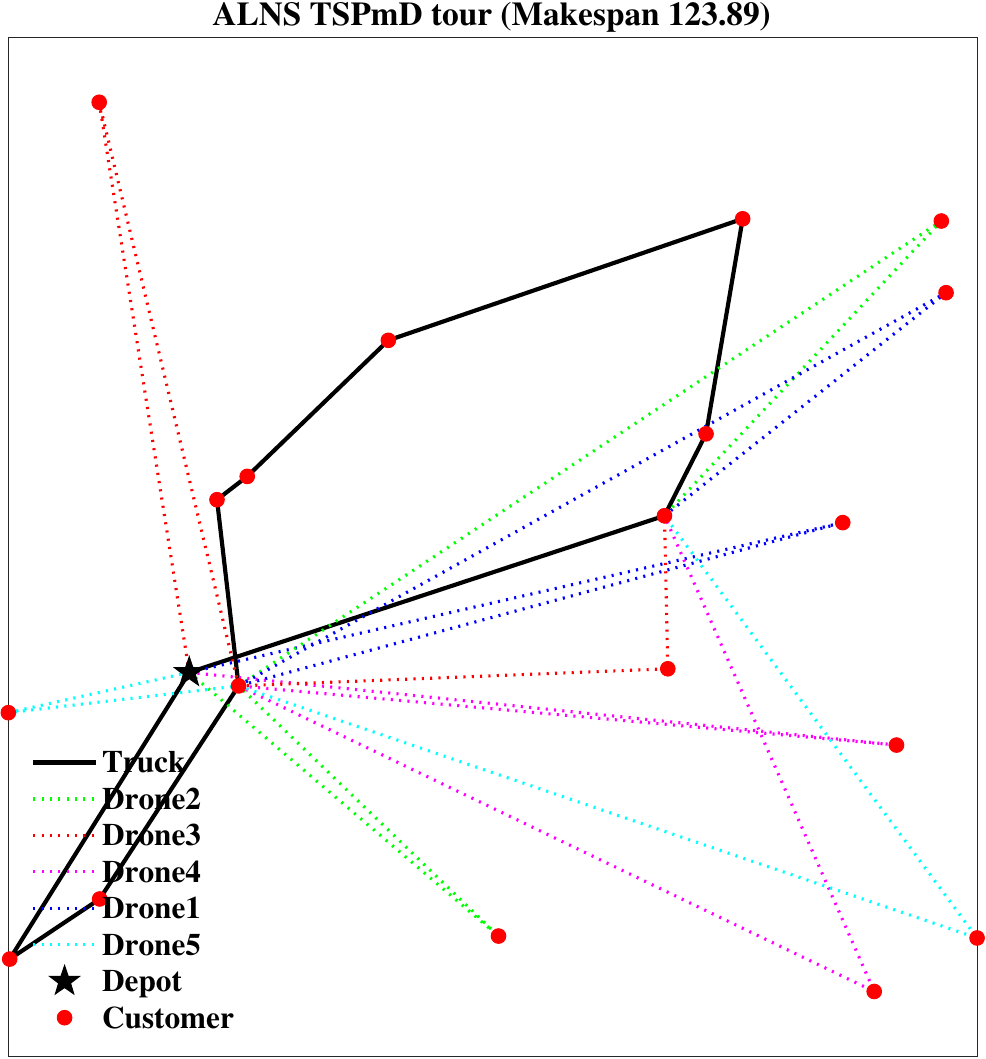} &
        \includegraphics[width=0.33\textwidth]{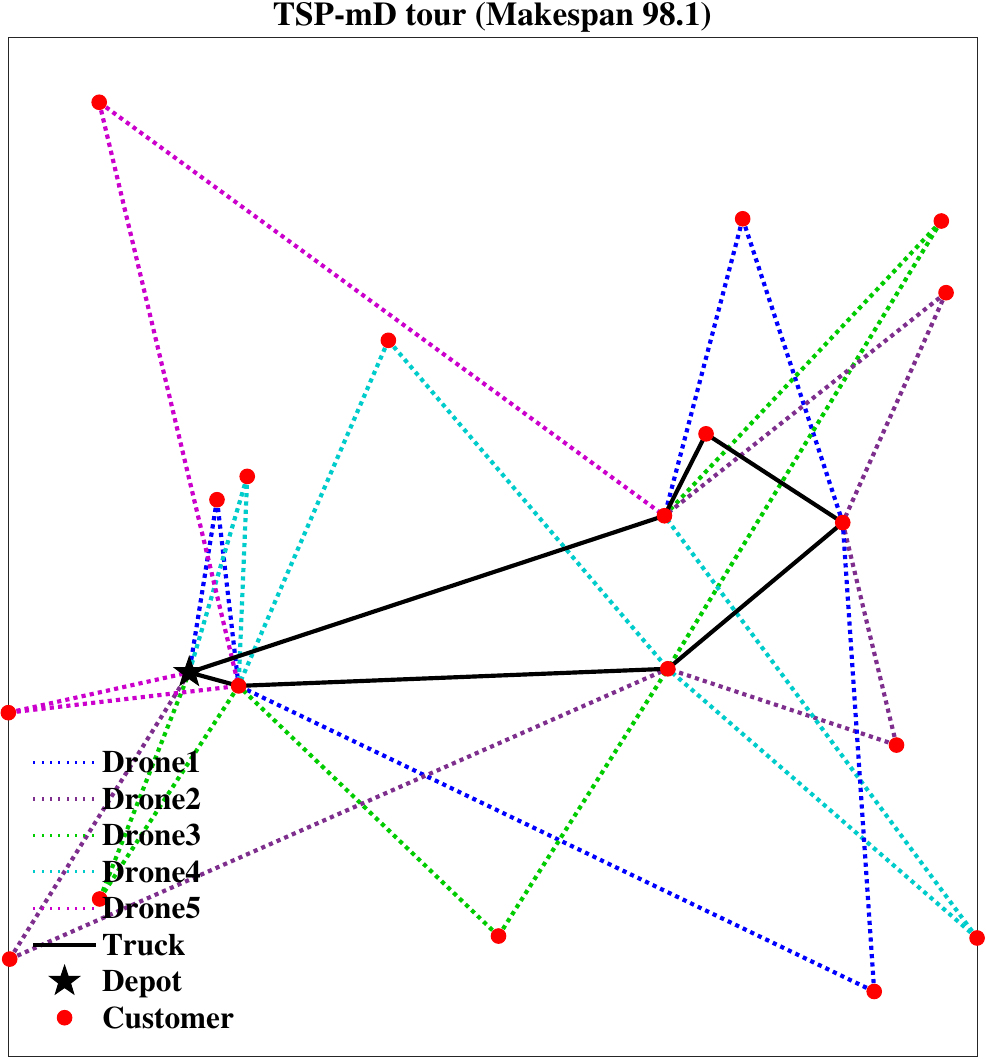} &
        \includegraphics[width=0.33\textwidth]{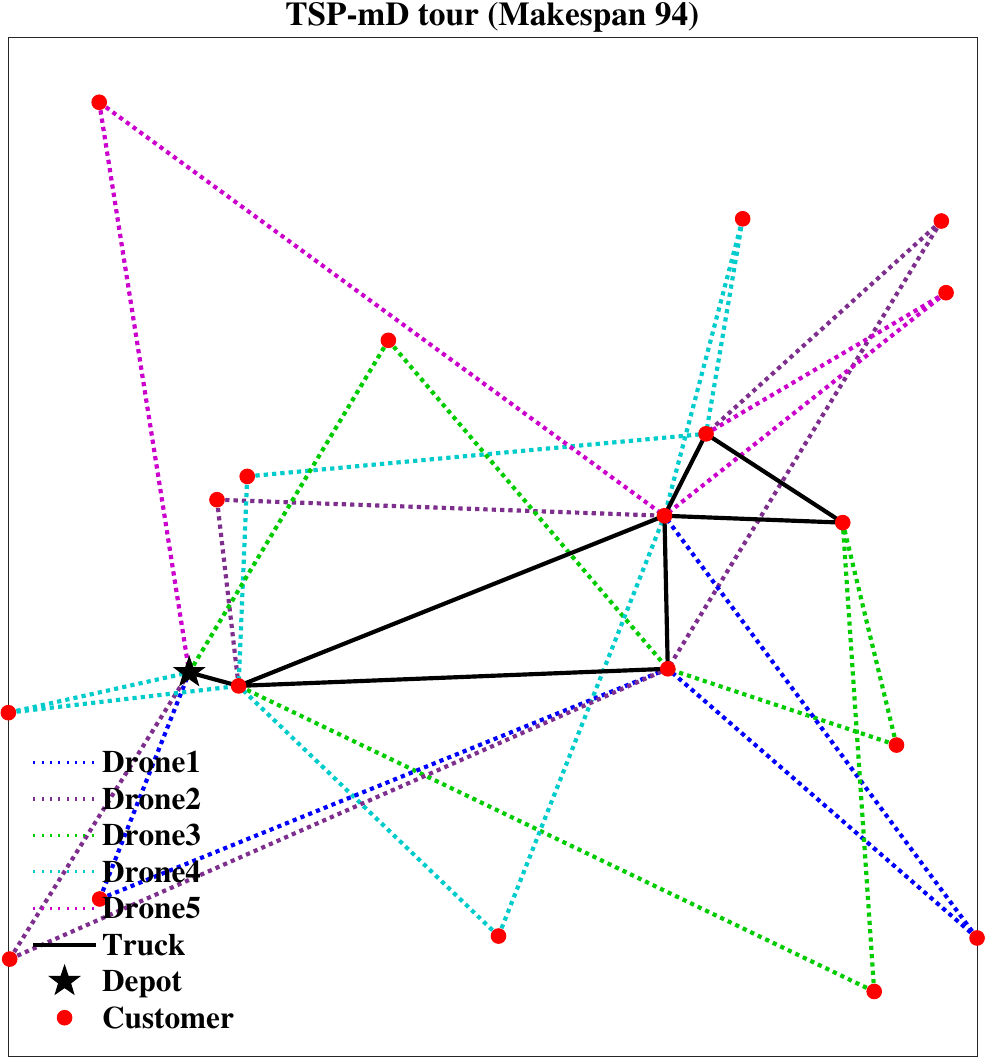} \\
        \multicolumn{1}{c}{(a) ALNS.} & \multicolumn{1}{c}{(b) LKH Config1.} & \multicolumn{1}{c}{(c) LKH Config2.}
    \end{tabular}
    \caption{TSPmD tours on instance poi-20-25.}
    \label{fig:tspmd_figcompare1}
\end{figure}

Firstly, according to the comparison between two configurations in Tables~\ref{tab:tspmd_2LKH_small}--\ref{tab:tspmd_2LKH_medium}, LKH Config2 adopts the node revisit option in bold instances. Taking poi-20-25 for example (Fig.~\ref{fig:tspmd_figcompare1}), LKH Config2 reduces the makespan from 98.1 under Config1 to 94 by revisiting two customer nodes. Both LKH tours outperform ALNS, which produces the highest makespan of 123.89. Note that all five available drones are utilized in these three tours, with each drone's path represented by dashed lines of distinct colors.

\begin{figure}[H]
    \centering
    \setlength{\tabcolsep}{1pt} % Space between columns
    \begin{tabular}{ccc}
        \includegraphics[width=0.33\textwidth]{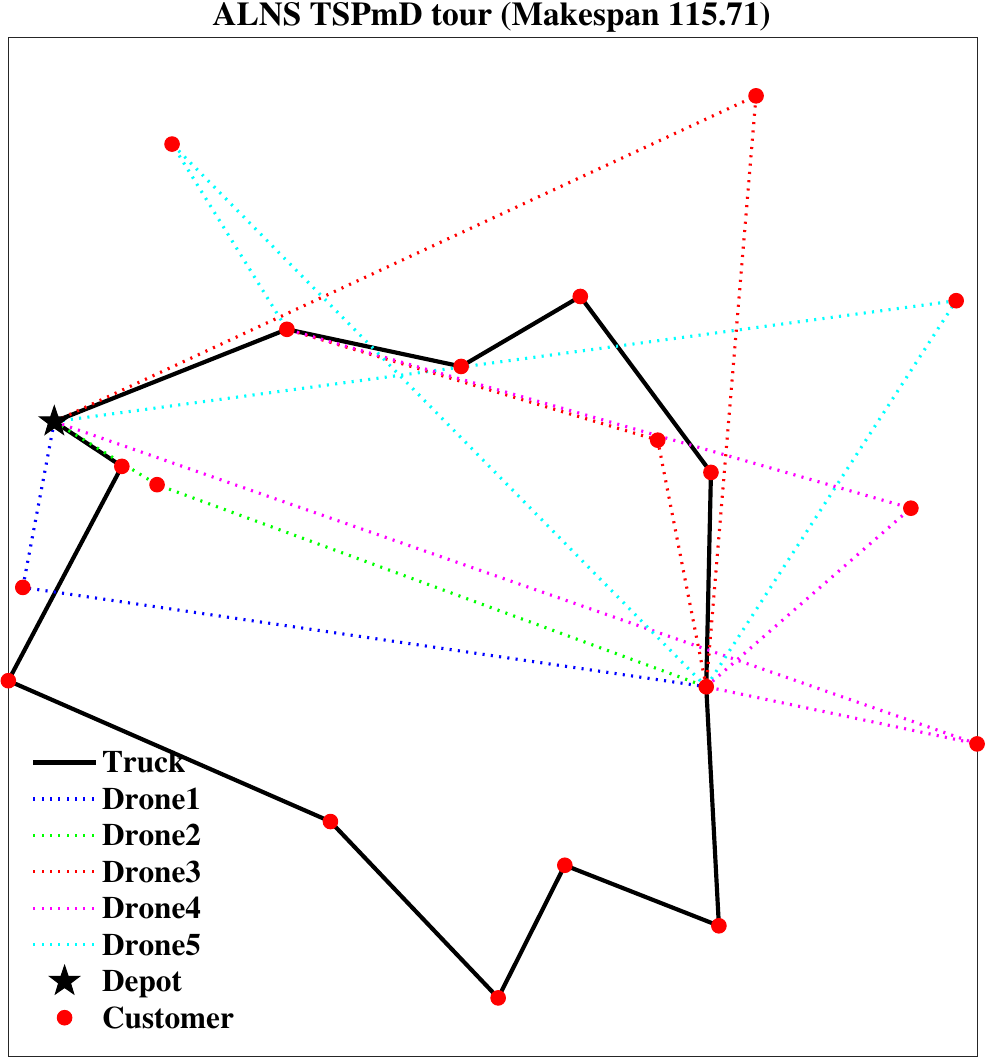} &
        \includegraphics[width=0.33\textwidth]{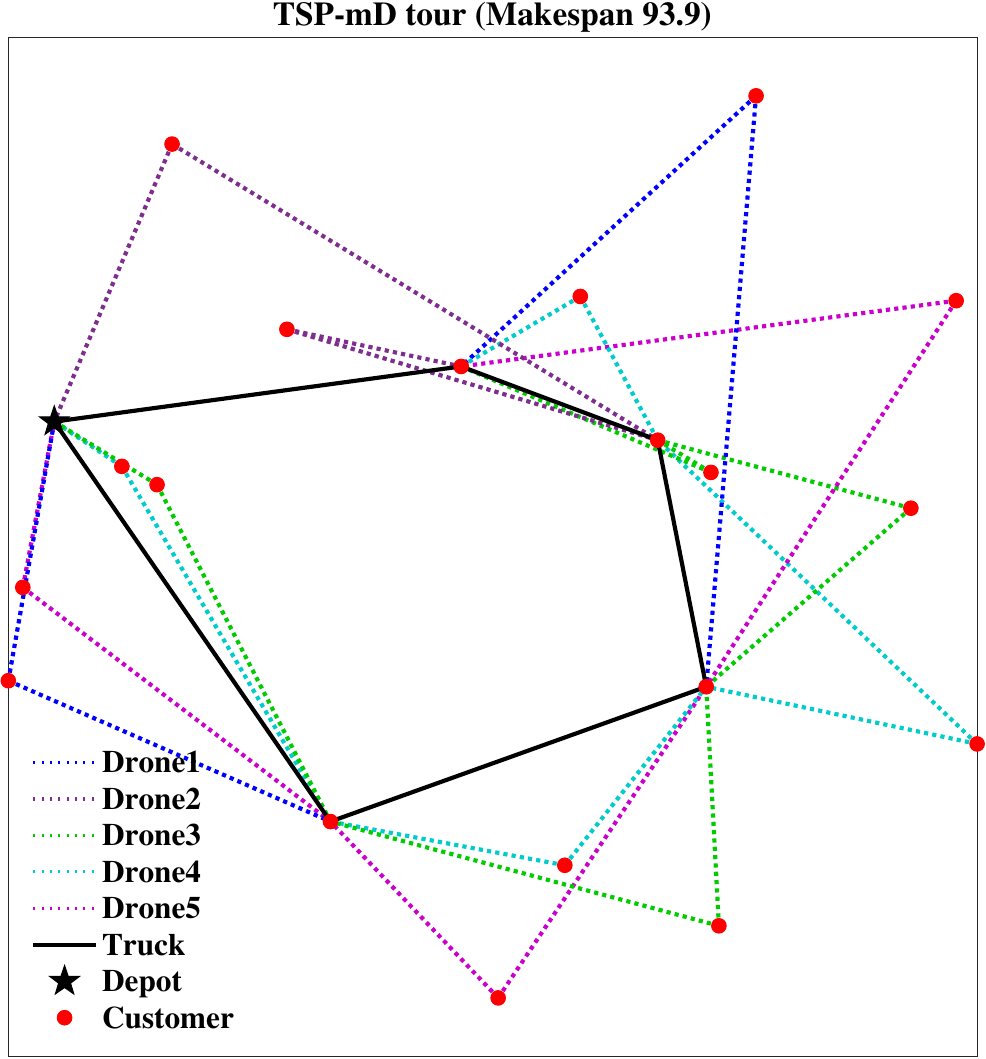} &
        \includegraphics[width=0.33\textwidth]{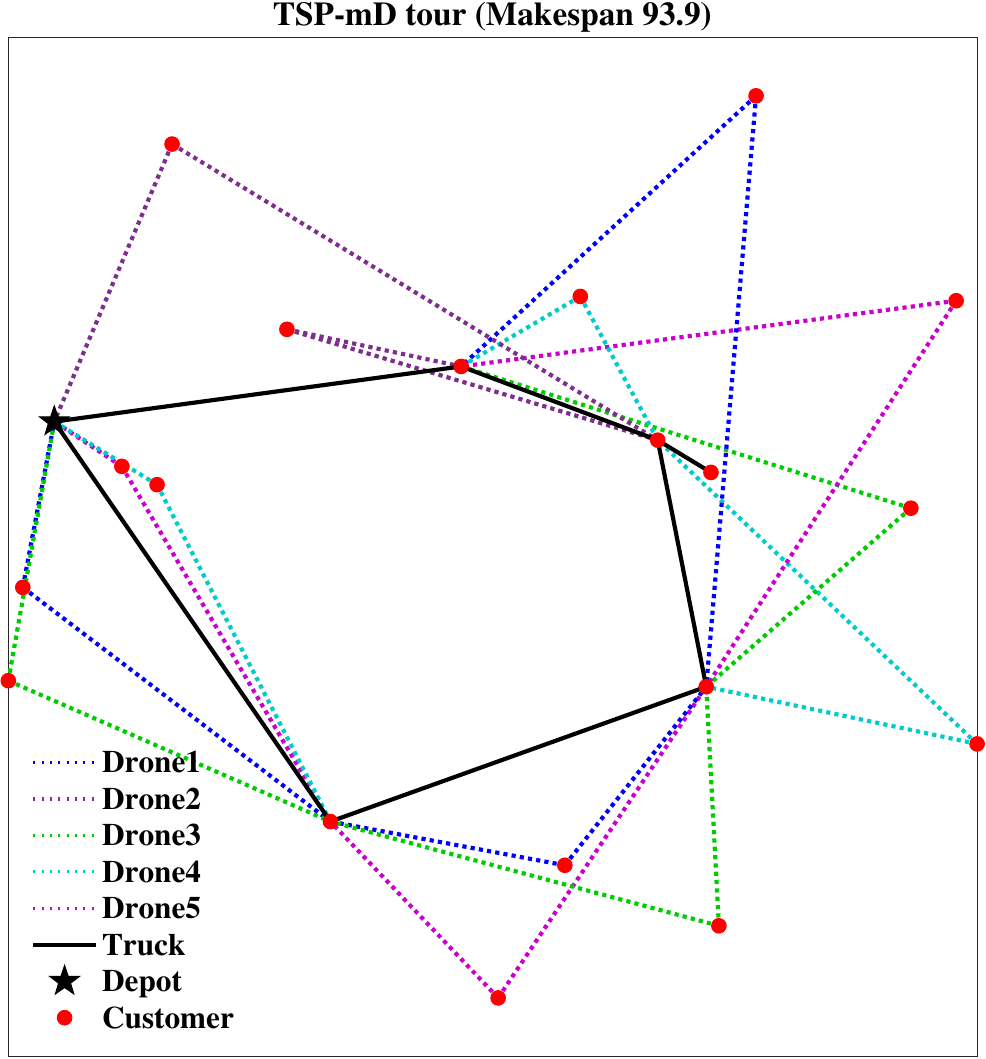} \\
        \multicolumn{1}{c}{(a) ALNS.} & \multicolumn{1}{c}{(b) LKH Config1.} & \multicolumn{1}{c}{(c) LKH Config2.}
    \end{tabular}
    \caption{TSPmD tours on instance poi-20-19.}
    \label{fig:tspmd_figcompare2}
\end{figure}

However, this makespan reduction by revisiting nodes is not consistent. As shown in Fig.~\ref{fig:tspmd_figcompare2} for instance poi-20-19, allowing the truck to revisit nodes does not affect the makespan, which remains at 93.9 under both LKH configurations.

\begin{figure}[H]
    \centering
    \setlength{\tabcolsep}{1pt} % Space between columns
    \begin{tabular}{ccc}
        \includegraphics[width=0.33\textwidth]{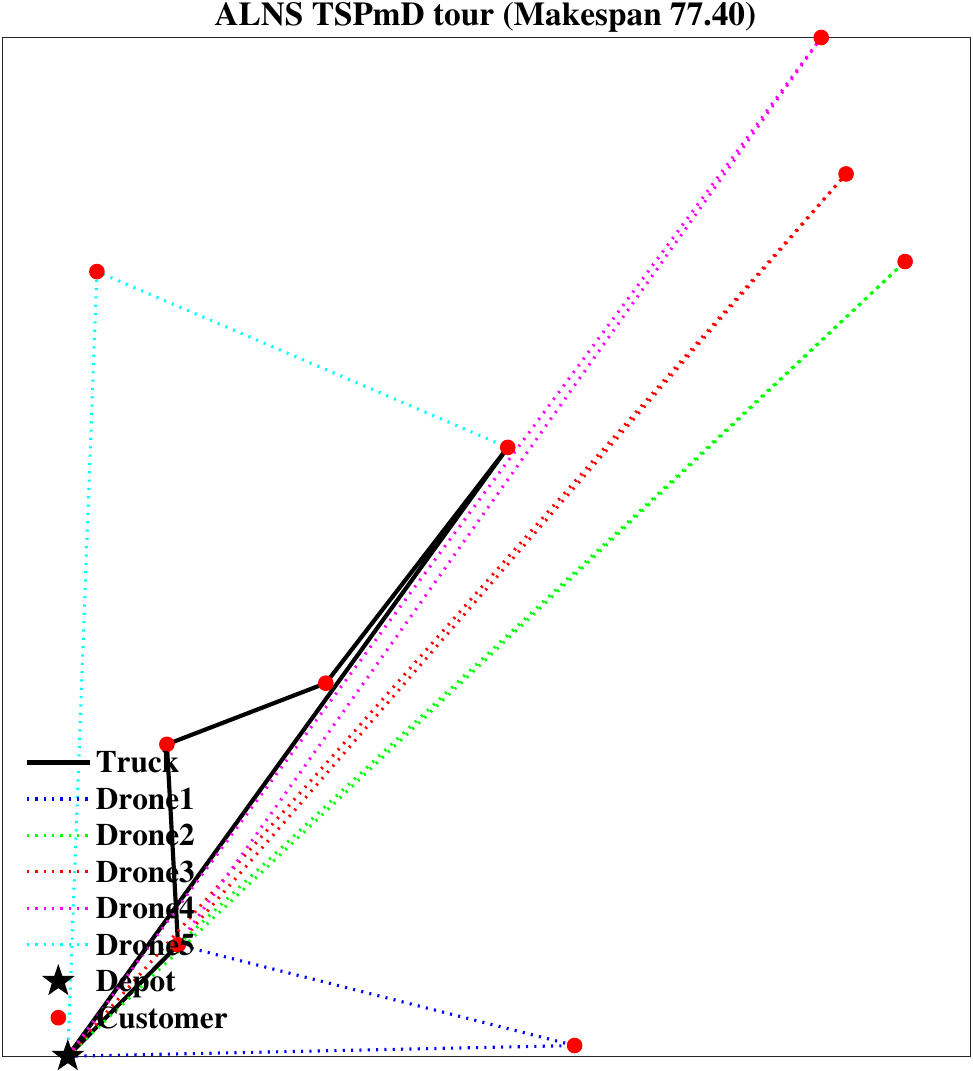} &
        \includegraphics[width=0.33\textwidth]{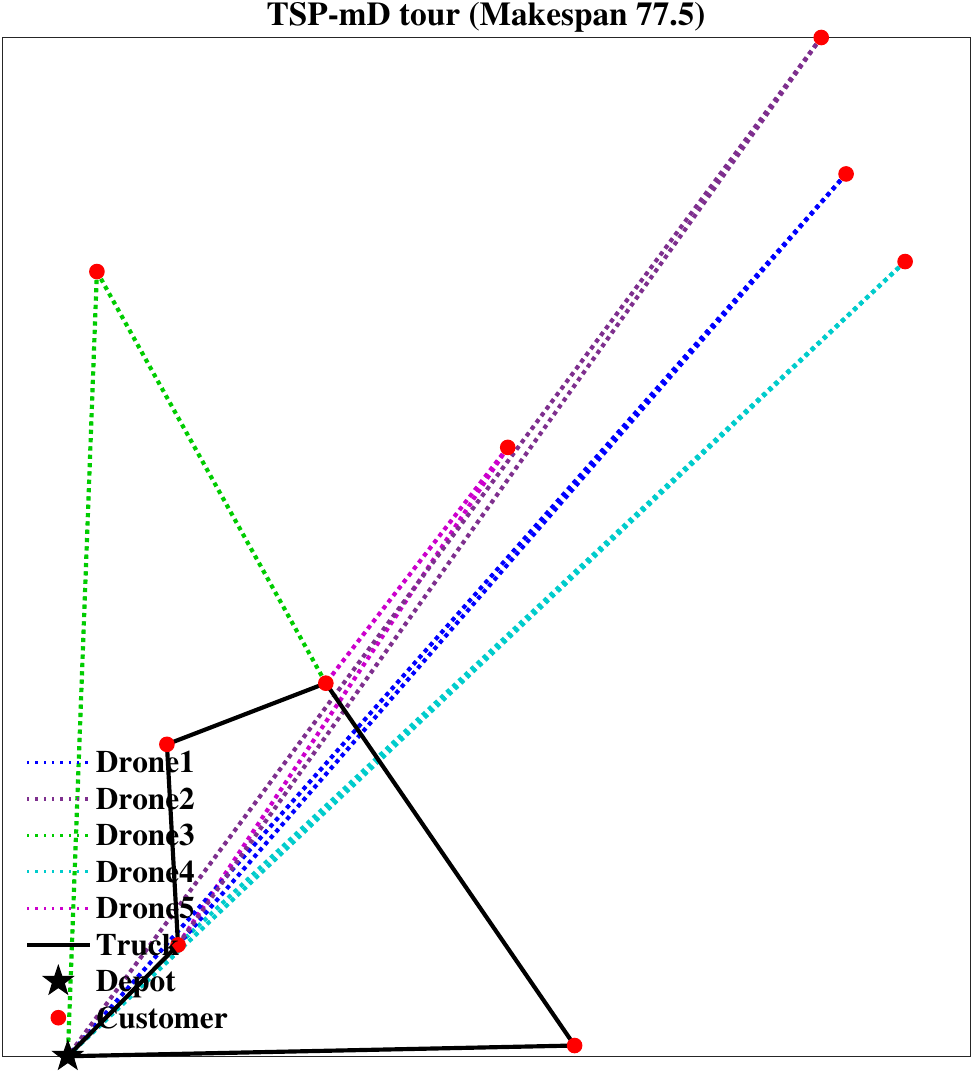} &
        \includegraphics[width=0.33\textwidth]{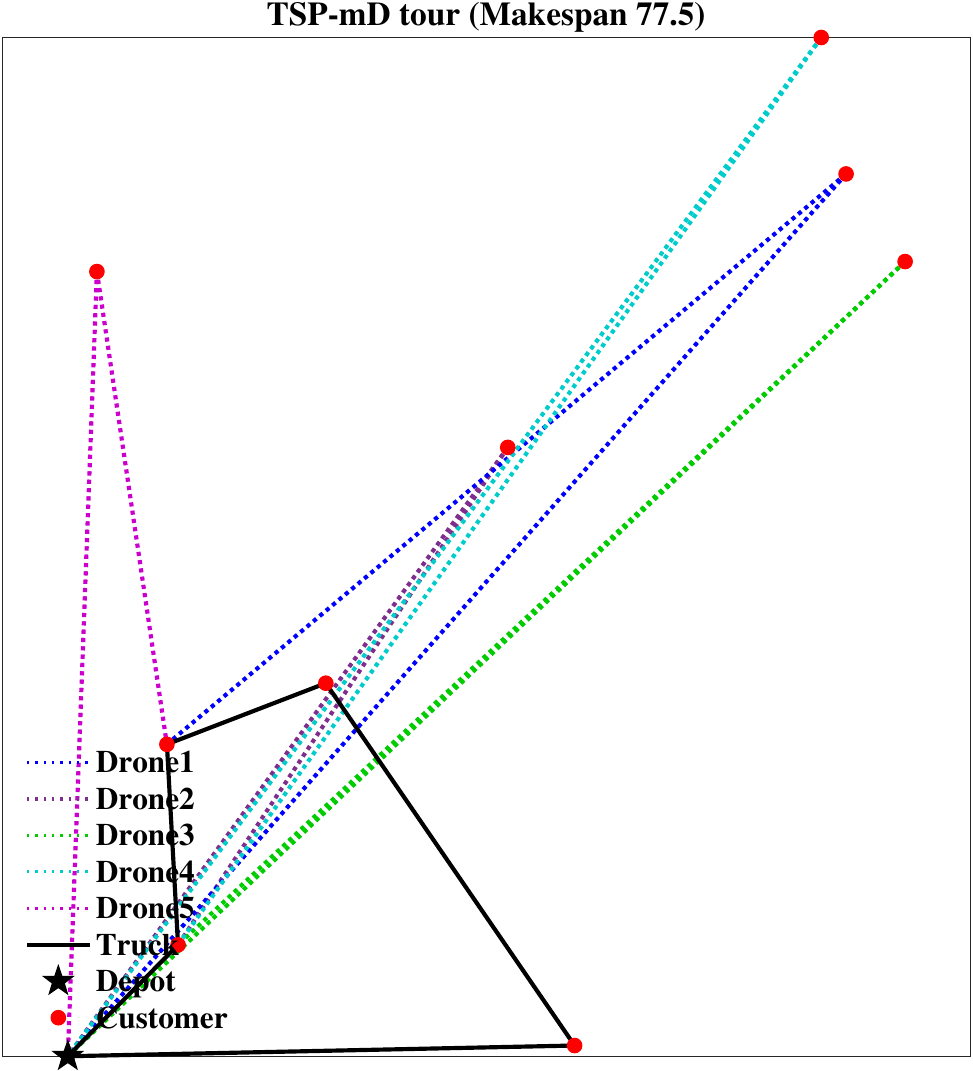} \\
        \multicolumn{1}{c}{(a) ALNS.} & \multicolumn{1}{c}{(b) LKH Config1.} & \multicolumn{1}{c}{(c) LKH Config2.}
    \end{tabular}
    \caption{TSPmD tours on instance poi-10-4.}
    \label{fig:tspmd_figcompare3}
\end{figure}

Secondly, as shown in Tables~\ref{tab:tspmd_1LKH_small}--\ref{tab:tspmd_1LKH_medium}, LKH Config1 (executed once within 5min) does not always outperform ALNS. For instance poi-10-4 in Fig.~\ref{fig:tspmd_figcompare3}, ALNS achieves a smaller makespan than LKH.

\subsubsection{Computational performance analysis}
For small TSP-mD instances, Table~\ref{tab:tspmd_1LKH_small} provides a detailed comparison between the optimal solutions of the MIP model solved by CPLEX, and the heuristic solutions generated by ALNS and LKH Config1. The table extends the format of Table~\ref{tab:fstsp_1LKH_small} (FSTSP) by adding a new column, \#D, to specify the number of drones used. Notably, ALNS does not produce valid tours for instances \textbf{poi-7-23} and \textbf{poi-10-23} throughout all five runs. 

% TSP-mD: cplex 103/25=4.12, 229.15/25=9.166s, lkh 96/25=3.83, 27.35/25=1.094s

The transition from one-drone FSTSP to five-drone TSP-mD configurations makes it harder to optimally solve the MIP model. In Table~\ref{tab:tspmd_1LKH_small}, CPLEX consumes more time (9.12s on average) to solve poi-7- instances compared to 2.99s in the previous Table~\ref{tab:fstsp_1LKH_small} for FSTSP. Moreover, it is unable to identify optimal solutions for poi-10- instances within a runtime threshold of 4h. In contrast to CPLEX, for poi-7- instances, LKH Config1 not only matches all optimal values ($\%gap^* = 100 \cdot (z^*-\text{M}^1)/z^* = 0 $), but also requires only 1.09s on average, which is less than the 6.24s needed for FSTSP in Table~\ref{tab:fstsp_1LKH_small}. For poi-10- instances, the average time by LKH also decreases from 53.83s (FSTSP) to 20.26s (TSP-mD). This decrease in computational time from FSTSP to TSP-mD may result from the penalty term that we design in LKH. Infeasible solutions for FSTSP (e.g. 1 $<$ \# used drones $<$ 5) in early iterations could be feasible for TSP-mD, which means that LKH satisfies the drone number constraint earlier for TSP-mD. Furthermore, LKH Config1 on average deploys fewer drones (3.84) than CPLEX (4.12) for poi-7- instances, suggesting potential operational cost savings.

In addition, positive values in the column $\%gap^{\textit{ALNS}} = 100 \cdot (M^{\textit{ALNS}}-\text{M}^1)/M^{\textit{ALNS}}$ reflect makespan reductions achieved by LKH Config1 compared to ALNS. Therefore, the majority of positive gaps in Tables~\ref{tab:tspmd_1LKH_small}--\ref{tab:tspmd_1LKH_medium} demonstrate the superior performance of LKH over ALNS.

% first comparison
\begin{table}[H]
    \centering
    % \tiny
    % \scriptsize
    \footnotesize  % Slightly smaller to better fit
    % \small
    % \normalsize  % (default size) 
    % \large
    \caption{TSPmD Comparison of ALNS and LKH (Config 1) on small instances}
    \label{tab:tspmd_1LKH_small}
    \begin{tabular*}{\textwidth}{@{\extracolsep{\fill}}l@{\;}r@{\;}r@{\;}r@{\;}r@{\;}r@{\;}r@{\;}r@{\;}r@{\;}r@{\;}r@{\;}r@{\;}r@{\;}r@{\;}r@{\;}r@{\;}r@{}}
    %%
    % \begin{tabularx}{\textwidth}{l r r r r r r r r r r r r r r r}
    % \begin{tabularx}{\textwidth}{l c c c c c c c c c c c}
    % \begin{tabularx}{\textwidth}{l *{11}{>{\centering\arraybackslash}X}}
        \toprule
        \multicolumn{2}{c}{} & \multicolumn{4}{c}{\textbf{CPLEX}} & \multicolumn{4}{c}{\textbf{ALNS}} & \multicolumn{5}{c}{\textbf{LKH (Config 1)}} & \multicolumn{1}{c}{} & \multicolumn{1}{c}{} \\
        \cmidrule(lr){3-6} \cmidrule(lr){7-10} \cmidrule(lr){11-15}
        \raisebox{1.5ex}[0pt][0pt]{Instance} & \raisebox{1.5ex}[0pt][0pt] {$N_D$} & $z^*$ & Dcus & \#D & $t^{\textit{MIP}}$(s) & $M^{\textit{ALNS}}$ & Dcus & \#D & Time(s) & $\text{M}^1$ & Dcus & \#D & Iter & Time(s) & \raisebox{1.5ex}[0pt][0pt] {$\%gap^*$}  & \raisebox{1.5ex}[0pt][0pt] {$\%gap^{\textit{ALNS}}$} \\
        \midrule
poi-7-1 & 6 & 90 & 4 & 4 & 11.34 & 90.07 & 5 & 5 & 300.01 & 90 & 3 & 3 & 126 & 0.4 & 0 & 0.08 \\
poi-7-2 & 6 & 65.2 & 4 & 4 & 8.87 & 65.22 & 5 & 5 & 5.15 & 65.2 & 4 & 4 & 55 & 0.15 & 0 & 0.03 \\
poi-7-3 & 6 & 73.8 & 4 & 4 & 12 & 73.87 & 2 & 2 & 66.55 & 73.8 & 1 & 1 & 229 & 0.71 & 0 & 0.09 \\
poi-7-4 & 6 & 81.2 & 4 & 4 & 12.09 & 81.2 & 5 & 5 & 267.03 & 81.2 & 3 & 3 & 187 & 0.48 & 0 & 0 \\
poi-7-5 & 6 & 67.3 & 4 & 4 & 8.55 & 67.31 & 4 & 4 & 300 & 67.3 & 4 & 4 & 10 & 0.02 & 0 & 0.01 \\
poi-7-6 & 6 & 72.8 & 5 & 5 & 9.05 & 72.8 & 5 & 5 & 300.01 & 72.8 & 5 & 5 & 9 & 0.02 & 0 & 0 \\
poi-7-7 & 6 & 55 & 4 & 4 & 8.68 & 55.05 & 5 & 5 & 300.01 & 55 & 4 & 4 & 11 & 0.04 & 0 & 0.09 \\
poi-7-8 & 6 & 58.6 & 5 & 5 & 11.25 & 58.6 & 5 & 5 & 300 & 58.6 & 5 & 5 & 16 & 0.04 & 0 & 0 \\
poi-7-9 & 6 & 87.3 & 5 & 5 & 8.74 & 87.33 & 5 & 5 & 300.01 & 87.3 & 5 & 5 & 3 & 0.02 & 0 & 0.03 \\
poi-7-10 & 6 & 59.3 & 4 & 4 & 8.42 & 59.49 & 5 & 5 & 300.01 & 59.3 & 4 & 4 & 2309 & 5.51 & 0 & 0.32 \\
poi-7-11 & 6 & 69.3 & 3 & 3 & 11.68 & 69.3 & 5 & 5 & 300.01 & 69.3 & 3 & 3 & 9 & 0.02 & 0 & 0 \\
poi-7-12 & 6 & 66.2 & 3 & 3 & 8.86 & 66.2 & 5 & 5 & 7.67 & 66.2 & 3 & 3 & 99 & 0.28 & 0 & 0 \\
poi-7-13 & 6 & 50.2 & 4 & 4 & 9.05 & 50.2 & 5 & 5 & 300 & 50.2 & 4 & 4 & 1 & 0.01 & 0 & 0 \\
poi-7-14 & 6 & 66.8 & 4 & 4 & 11.24 & 66.8 & 5 & 5 & 300.01 & 66.8 & 4 & 4 & 34 & 0.05 & 0 & 0 \\
poi-7-15 & 6 & 75.4 & 5 & 5 & 8.44 & 75.4 & 5 & 5 & 300 & 75.4 & 4 & 4 & 10 & 0.02 & 0 & 0 \\
poi-7-16 & 6 & 59.3 & 3 & 3 & 8.56 & 59.3 & 3 & 3 & 300 & 59.3 & 3 & 3 & 161 & 0.41 & 0 & 0 \\
poi-7-17 & 6 & 64.6 & 4 & 4 & 11.73 & 64.67 & 5 & 5 & 300 & 64.6 & 4 & 4 & 1 & 0.01 & 0 & 0.11 \\
poi-7-18 & 6 & 63.2 & 4 & 4 & 8.82 & 63.22 & 5 & 5 & 4.45 & 63.2 & 4 & 4 & 59 & 0.09 & 0 & 0.03 \\
poi-7-19 & 6 & 72.2 & 4 & 4 & 10.73 & 72.26 & 5 & 5 & 5.26 & 72.2 & 4 & 4 & 37631 & 15.31 & 0 & 0.08 \\
poi-7-20 & 6 & 71.2 & 5 & 5 & 8.75 & 71.23 & 5 & 5 & 177.03 & 71.2 & 5 & 5 & 66 & 0.09 & 0 & 0.04 \\
poi-7-21 & 6 & 47.5 & 3 & 3 & 10.97 & 47.67 & 4 & 4 & 124.29 & 47.5 & 3 & 3 & 148 & 0.5 & 0 & 0.36 \\
poi-7-22 & 6 & 66.6 & 4 & 4 & 10.96 & 66.61 & 5 & 5 & 300.01 & 66.6 & 4 & 4 & 18 & 0.04 & 0 & 0.02 \\
poi-7-23 & 6 & 80.2 & 4 & 4 & 8.6 & - & - & - & - & 80.2 & 4 & 4 & 17 & 0.06 & 0 & - \\
poi-7-24 & 6 & 46.7 & 4 & 4 & 11.27 & 46.7 & 4 & 4 & 300.01 & 46.7 & 3 & 3 & 639 & 2.39 & 0 & 0 \\
poi-7-25 & 6 & 68 & 5 & 5 & 9.4 & 68.05 & 5 & 5 & 245.43 & 68 & 3 & 3 & 298 & 0.59 & 0 & 0.07 \\
Average & & & & 4.12 & 9.12 & & & & & & & 3.84 & & 1.09 & 0 & \\
poi-10-1 & 9 & & & & & 90.07 & 5 & 5 & 33.9 & 90 & 4 & 4 & 5860 & 32.87 & & 0.08 \\
poi-10-2 & 9 & & & & & 82.83 & 5 & 5 & 24.23 & 65 & 6 & 5 & 13254 & 106.51 & & 21.53 \\
poi-10-3 & 9 & & & & & 78.35 & 5 & 5 & 8.58 & 73.9 & 4 & 4 & 530 & 3.09 & & 5.68 \\
poi-10-4 & 9 & & & & & 77.4 & 5 & 5 & 30.91 & 77.5 & 5 & 5 & 274 & 1.21 & & -0.13 \\
poi-10-5 & 9 & & & & & 67.49 & 5 & 5 & 26.79 & 65.8 & 7 & 5 & 119 & 0.69 & & 2.5 \\
poi-10-6 & 9 & & & & & 90.8 & 5 & 5 & 25.38 & 77.5 & 6 & 5 & 185 & 0.96 & & 14.65 \\
poi-10-7 & 9 & & & & & 49.54 & 5 & 5 & 9.55 & 49.6 & 5 & 5 & 488 & 2.83 & & -0.12 \\
poi-10-8 & 9 & & & & & 57.19 & 5 & 5 & 34.34 & 57.2 & 4 & 4 & 219 & 1.12 & & -0.02 \\
poi-10-9 & 9 & & & & & 102.44 & 5 & 5 & 32.28 & 89.5 & 8 & 5 & 330 & 1.8 & & 12.63 \\
poi-10-10 & 9 & & & & & 56.96 & 5 & 5 & 10.3 & 56.9 & 4 & 4 & 473 & 2.77 & & 0.11 \\
poi-10-11 & 9 & & & & & 64.98 & 5 & 5 & 34.52 & 65 & 4 & 4 & 586 & 2.89 & & -0.03 \\
poi-10-12 & 9 & & & & & 66.2 & 5 & 5 & 9.52 & 66.2 & 4 & 4 & 400 & 1.99 & & 0 \\
poi-10-13 & 9 & & & & & 63.63 & 5 & 5 & 12.82 & 52.4 & 7 & 5 & 7704 & 34.82 & & 17.65 \\
poi-10-14 & 9 & & & & & 71.25 & 5 & 5 & 23.89 & 65.1 & 6 & 5 & 1373 & 9.67 & & 8.63 \\
poi-10-15 & 9 & & & & & 65.61 & 5 & 5 & 7.29 & 64.9 & 4 & 4 & 714 & 3.15 & & 1.08 \\
poi-10-16 & 9 & & & & & 69.75 & 5 & 5 & 33.08 & 64.7 & 6 & 5 & 307 & 2.25 & & 7.24 \\
poi-10-17 & 9 & & & & & 84.18 & 5 & 5 & 27.49 & 64.6 & 6 & 5 & 5158 & 33.46 & & 23.26 \\
poi-10-18 & 9 & & & & & 74.73 & 5 & 5 & 9.83 & 63.2 & 7 & 5 & 27804 & 163.16 & & 15.43 \\
poi-10-19 & 9 & & & & & 76.62 & 5 & 5 & 12.73 & 73 & 6 & 5 & 5797 & 40.14 & & 4.72 \\
poi-10-20 & 9 & & & & & 85.57 & 5 & 5 & 21.31 & 74.5 & 7 & 5 & 1955 & 11.9 & & 12.94 \\
poi-10-21 & 9 & & & & & 54.89 & 5 & 5 & 29.52 & 47.5 & 6 & 5 & 120 & 0.43 & & 13.46 \\
poi-10-22 & 9 & & & & & 63.04 & 6 & 5 & 9.28 & 61.9 & 4 & 4 & 739 & 4.24 & & 1.81 \\
poi-10-23 & 9 & & & & & - & - & - & - & 80.2 & 5 & 5 & 305 & 1.67 & & - \\
poi-10-24 & 9 & & & & & 76.62 & 5 & 5 & 9.09 & 73.9 & 4 & 4 & 991 & 5.27 & & 3.55 \\
poi-10-25 & 9 & & & & & 88.38 & 5 & 5 & 26.85 & 74.9 & 7 & 5 & 166 & 0.61 & & 15.25 \\
Average & & & & & & & & & & & & & & 20.26 & & \\
        \bottomrule
    \end{tabular*}
\end{table}

\begin{table}[H]
    \centering
    % \tiny
    % \scriptsize
    \footnotesize  % Slightly smaller to better fit
    % \small
    % \normalsize  % (default size) 
    % \large
    \caption{TSPmD Comparison of ALNS and LKH (Config 1) on medium instances}
    \label{tab:tspmd_1LKH_medium}
    \begin{tabular*}{\textwidth}{@{\extracolsep{\fill}}l@{\;}r@{\;}r@{\;}r@{\;}r@{\;}r@{\;}r@{\;}r@{\;}r@{\;}r@{\;}r@{\;}r@{}}
    %%
    % \begin{tabularx}{\textwidth}{l r r r r r r r r r r r}
    % \begin{tabularx}{\textwidth}{l c c c c c c c c c c c}
    % \begin{tabularx}{\textwidth}{l *{11}{>{\centering\arraybackslash}X}}
        \toprule
        \multicolumn{2}{c}{} & \multicolumn{4}{c}{\textbf{ALNS}} & \multicolumn{5}{c}{\textbf{LKH (Config 1)}} & \multicolumn{1}{c}{}\\
        \cmidrule(lr){3-6} \cmidrule(lr){7-11}
        \raisebox{1.5ex}[0pt][0pt]{Instance} & \raisebox{1.5ex}[0pt][0pt] {$N_D$} & $M^{\textit{ALNS}}$ & Dcus & \#D & Time(s) & $\text{M}^1$ & Dcus & \#D & Iter & Time(s) & \raisebox{1.5ex}[0pt][0pt] {$\%gap^{\textit{ALNS}}$} \\
        \midrule
poi-20-1 & 19 & 129.74 & 10 & 5 & 300 & 106.4 & 12 & 5 & 1642 & 35.41 & 17.99 \\
poi-20-2 & 19 & 130.96 & 10 & 5 & 131.38 & 103 & 13 & 5 & 6365 & 96.96 & 21.35 \\
poi-20-3 & 19 & 152.36 & 11 & 5 & 300.01 & 108.4 & 14 & 5 & 16139 & 253.29 & 28.85 \\
poi-20-4 & 19 & 131.32 & 10 & 5 & 141.7 & 89.6 & 12 & 5 & 12208 & 238.58 & 31.77 \\
poi-20-5 & 19 & 118.16 & 10 & 5 & 266.84 & 94.1 & 15 & 5 & 10451 & 226.77 & 20.36 \\
poi-20-6 & 19 & 107.76 & 13 & 5 & 144.38 & 85.5 & 14 & 5 & 15968 & 276.3 & 20.66 \\
poi-20-7 & 19 & 129.33 & 10 & 5 & 46.86 & 111.1 & 11 & 5 & 1804 & 27 & 14.1 \\
poi-20-8 & 19 & 128.28 & 11 & 5 & 142.1 & 118.7 & 10 & 5 & 7080 & 143.1 & 7.47 \\
poi-20-9 & 19 & 128.04 & 9 & 5 & 300.01 & 113.9 & 12 & 5 & 15617 & 255.69 & 11.04 \\
poi-20-10 & 19 & 129.34 & 10 & 5 & 300.01 & 97.5 & 14 & 5 & 11607 & 217.25 & 24.62 \\
poi-20-11 & 19 & 129.39 & 11 & 5 & 65.07 & 109.3 & 13 & 5 & 6652 & 116.48 & 15.53 \\
poi-20-12 & 19 & 110.16 & 12 & 5 & 128.03 & 99.9 & 13 & 5 & 6352 & 93.69 & 9.31 \\
poi-20-13 & 19 & 116.94 & 8 & 5 & 67.08 & 88 & 13 & 5 & 13142 & 297.71 & 24.75 \\
poi-20-14 & 19 & 122.99 & 10 & 5 & 152.33 & 103.7 & 12 & 5 & 9621 & 183.33 & 15.68 \\
poi-20-15 & 19 & 138.82 & 11 & 5 & 81.42 & 100 & 13 & 5 & 676 & 12.71 & 27.96 \\
poi-20-16 & 19 & 116.24 & 11 & 5 & 172.05 & 99.2 & 14 & 5 & 14084 & 231.6 & 14.66 \\
poi-20-17 & 19 & 116.42 & 9 & 5 & 193.94 & 100.1 & 11 & 5 & 7958 & 168.22 & 14.02 \\
poi-20-18 & 19 & 141.74 & 6 & 5 & 300 & 109.5 & 13 & 5 & 15156 & 277.12 & 22.75 \\
poi-20-19 & 19 & 115.71 & 8 & 5 & 256.68 & 93.9 & 15 & 5 & 4739 & 84.35 & 18.85 \\
poi-20-20 & 19 & 149.16 & 10 & 5 & 162.72 & 120.9 & 13 & 5 & 8623 & 132.6 & 18.95 \\
poi-20-21 & 19 & 134.55 & 11 & 5 & 224.66 & 99.6 & 12 & 5 & 2407 & 41.12 & 25.98 \\
poi-20-22 & 19 & 123.67 & 11 & 5 & 129.14 & 98.1 & 11 & 5 & 9452 & 204.49 & 20.68 \\
poi-20-23 & 19 & 117.12 & 11 & 5 & 218.45 & 107.5 & 11 & 5 & 5480 & 98.55 & 8.21 \\
poi-20-24 & 19 & 123.92 & 8 & 5 & 67.63 & 97.8 & 13 & 5 & 4250 & 70.94 & 21.08 \\
poi-20-25 & 19 & 123.89 & 10 & 5 & 72.43 & 98.1 & 14 & 5 & 963 & 12.48 & 20.82 \\
poi-30-1 & 29 & 180.29 & 18 & 5 & 300 & 127 & 16 & 5 & 8925 & 285.13 & 29.56 \\
poi-30-2 & 29 & 147.83 & 20 & 5 & 200.17 & 118.5 & 20 & 5 & 7652 & 284.62 & 19.84 \\
poi-30-3 & 29 & 136.16 & 18 & 5 & 300.03 & 124.9 & 20 & 5 & 6341 & 232.16 & 8.27 \\
poi-30-4 & 29 & 147.49 & 18 & 5 & 300.04 & 128.5 & 20 & 5 & 3962 & 112.69 & 12.88 \\
poi-30-5 & 29 & 142.86 & 17 & 5 & 300.04 & 127.3 & 19 & 5 & 4593 & 198.29 & 10.89 \\
poi-30-6 & 29 & 154.07 & 15 & 5 & 300.06 & 128.7 & 19 & 5 & 7084 & 268.87 & 16.47 \\
poi-30-7 & 29 & 151.06 & 20 & 5 & 267.86 & 108.8 & 21 & 5 & 7335 & 198.3 & 27.98 \\
poi-30-8 & 29 & 152.47 & 15 & 5 & 300.03 & 141.3 & 22 & 5 & 10667 & 297.09 & 7.33 \\
poi-30-9 & 29 & 152.8 & 13 & 5 & 300.02 & 124.3 & 20 & 5 & 5880 & 266.11 & 18.65 \\
poi-30-10 & 29 & 191.26 & 8 & 5 & 300.01 & 142.5 & 21 & 5 & 6740 & 295.46 & 25.49 \\
poi-30-11 & 29 & 175.25 & 19 & 5 & 162.63 & 133.3 & 19 & 5 & 8533 & 293.55 & 23.94 \\
poi-30-12 & 29 & 178.27 & 16 & 5 & 300.03 & 132.2 & 20 & 5 & 9864 & 281.44 & 25.84 \\
poi-30-13 & 29 & 155.27 & 19 & 5 & 300.01 & 145.9 & 19 & 5 & 4932 & 167.33 & 6.03 \\
poi-30-14 & 29 & 165.11 & 19 & 5 & 300.05 & 130.8 & 23 & 5 & 4965 & 277.04 & 20.78 \\
poi-30-15 & 29 & 155.32 & 20 & 5 & 300.06 & 112.2 & 20 & 5 & 7168 & 287.62 & 27.76 \\
poi-30-16 & 29 & 153.76 & 15 & 5 & 300.02 & 116.6 & 21 & 5 & 5199 & 225.32 & 24.17 \\
poi-30-17 & 29 & 163.9 & 19 & 5 & 300.1 & 128.4 & 18 & 5 & 7010 & 263.85 & 21.66 \\
poi-30-18 & 29 & 163.16 & 10 & 5 & 300.02 & 121.5 & 23 & 5 & 1929 & 102.92 & 25.53 \\
poi-30-19 & 29 & 163.29 & 16 & 5 & 300.02 & 142.9 & 19 & 5 & 4262 & 186.43 & 12.49 \\
poi-30-20 & 29 & 169.53 & 4 & 4 & 300.02 & 122.2 & 18 & 5 & 3986 & 173.89 & 27.92 \\
poi-30-21 & 29 & 163.14 & 19 & 5 & 300.02 & 150.8 & 23 & 5 & 5681 & 222.18 & 7.56 \\
poi-30-22 & 29 & 171.77 & 11 & 5 & 300.02 & 142.2 & 20 & 5 & 4270 & 167.61 & 17.21 \\
poi-30-23 & 29 & 148.16 & 17 & 5 & 300 & 117 & 19 & 5 & 5308 & 243.85 & 21.03 \\
poi-30-24 & 29 & 145.95 & 18 & 5 & 300.01 & 129.6 & 20 & 5 & 8142 & 299.73 & 11.2 \\
poi-30-25 & 29 & 195.95 & 10 & 5 & 300.01 & 155.8 & 19 & 5 & 5945 & 229.09 & 20.49 \\
        \bottomrule
        \multicolumn{12}{r}{Continued on next page}
    \end{tabular*}
\end{table}

\clearpage % Force a page break here

\begin{table}[!t] % Force placement at the absolute top
    \vspace*{0pt} % Remove any vertical space at the top
    \centering
    % \tiny
    % \scriptsize
    \footnotesize  % Slightly smaller to better fit
    % \small
    % \normalsize  % (default size) 
    % \large
    % Use the caption command but suppress the numbering
    \refstepcounter{table}
    \addtocounter{table}{-1} % Keep the table number the same
    \makeatletter
    \def\@captiontype{table}
    {\@makecaption{\tablename~\thetable}{(continued)}}
    \vspace{0.5em}
    \makeatother
    \begin{tabular*}{\textwidth}{@{\extracolsep{\fill}}l@{\;}r@{\;}r@{\;}r@{\;}r@{\;}r@{\;}r@{\;}r@{\;}r@{\;}r@{\;}r@{\;}r@{}}
    %%
    % \begin{tabularx}{\textwidth}{l r r r r r r r r r r r}
    % \begin{tabularx}{\textwidth}{l c c c c c c c c c c c}
    % \begin{tabularx}{\textwidth}{l *{11}{>{\centering\arraybackslash}X}}
        \toprule
        \multicolumn{2}{c}{} & \multicolumn{4}{c}{\textbf{ALNS}} & \multicolumn{5}{c}{\textbf{LKH (Config 1)}} & \multicolumn{1}{c}{}\\
        \cmidrule(lr){3-6} \cmidrule(lr){7-11}
        \raisebox{1.5ex}[0pt][0pt]{Instance} & \raisebox{1.5ex}[0pt][0pt] {$N_D$} & $M^{\textit{ALNS}}$ & Dcus & \#D & Time(s) & $\text{M}^1$ & Dcus & \#D & Iter & Time(s) & \raisebox{1.5ex}[0pt][0pt] {$\%gap^{\textit{ALNS}}$} \\
        \midrule
poi-40-1 & 39 & 223.48 & 10 & 5 & 300.08 & 158.4 & 27 & 5 & 3060 & 203.2 & 29.12 \\
poi-40-2 & 39 & 232.72 & 5 & 5 & 300.07 & 165.1 & 28 & 5 & 3712 & 242.52 & 29.06 \\
poi-40-3 & 39 & 256.32 & 0 & 0 & 300 & 178.2 & 29 & 5 & 9707 & 283.86 & 30.48 \\
poi-40-4 & 39 & 179.55 & 17 & 5 & 300.02 & 186.1 & 27 & 5 & 2757 & 216.73 & -3.65 \\
poi-40-5 & 39 & 243.8 & 8 & 5 & 300.07 & 178.6 & 28 & 5 & 3237 & 200.05 & 26.74 \\
poi-40-6 & 39 & 218.82 & 9 & 5 & 300.02 & 154.8 & 28 & 5 & 2923 & 237.9 & 29.26 \\
poi-40-7 & 39 & 217.17 & 9 & 5 & 300.04 & 153.8 & 27 & 5 & 6652 & 278.1 & 29.18 \\
poi-40-8 & 39 & 197.41 & 23 & 5 & 300.02 & 182.1 & 30 & 5 & 4291 & 269.63 & 7.76 \\
poi-40-9 & 39 & 216.99 & 18 & 5 & 300.16 & 166.1 & 29 & 5 & 6585 & 281.77 & 23.45 \\
poi-40-10 & 39 & 218.64 & 5 & 5 & 300.1 & 160.9 & 29 & 5 & 4806 & 298.25 & 26.41 \\
poi-40-11 & 39 & 258.69 & 0 & 0 & 300.01 & 168.2 & 26 & 5 & 6275 & 286.17 & 34.98 \\
poi-40-12 & 39 & 236.97 & 5 & 5 & 300.18 & 180.6 & 30 & 5 & 4434 & 282.35 & 23.79 \\
poi-40-13 & 39 & 235.17 & 5 & 5 & 300.01 & 173.5 & 24 & 5 & 4926 & 298.79 & 26.22 \\
poi-40-14 & 39 & 241.43 & 5 & 5 & 300.01 & 154 & 27 & 5 & 5418 & 287.13 & 36.21 \\
poi-40-15 & 39 & 233.91 & 5 & 5 & 300.01 & 183.6 & 29 & 5 & 3721 & 267.12 & 21.51 \\
poi-40-16 & 39 & 222.42 & 26 & 5 & 300.04 & 161.1 & 27 & 5 & 3890 & 228.55 & 27.57 \\
poi-40-17 & 39 & 239.06 & 5 & 5 & 300.04 & 188.9 & 30 & 5 & 7845 & 248.01 & 20.98 \\
poi-40-18 & 39 & 239.44 & 10 & 5 & 300.01 & 195.9 & 28 & 5 & 3651 & 95.91 & 18.18 \\
poi-40-19 & 39 & 218.5 & 5 & 5 & 300.4 & 175.3 & 31 & 5 & 3144 & 244.05 & 19.77 \\
poi-40-20 & 39 & 211.9 & 9 & 5 & 300.01 & 176.8 & 26 & 5 & 4044 & 249.82 & 16.56 \\
poi-40-21 & 39 & 188.04 & 21 & 5 & 300.14 & 165.6 & 29 & 5 & 5632 & 275.79 & 11.93 \\
poi-40-22 & 39 & 250.43 & 5 & 5 & 300.01 & 177.3 & 28 & 5 & 4124 & 243.65 & 29.2 \\
poi-40-23 & 39 & 259.04 & 10 & 5 & 300.03 & 190.4 & 28 & 5 & 6172 & 287.4 & 26.5 \\
poi-40-24 & 39 & 226.48 & 5 & 5 & 300.02 & 167.4 & 29 & 5 & 4643 & 193.15 & 26.09 \\
poi-40-25 & 39 & 210.93 & 12 & 5 & 300.07 & 171.9 & 27 & 5 & 4330 & 297.94 & 18.5 \\
        \bottomrule
    \end{tabular*}
\end{table}

Regarding the negative gap for poi-40-4 in Table~\ref{tab:tspmd_1LKH_medium}, we extend the runtime of LKH Config1 from 5min to 4h. The resulting reduced makespan of 150.6 is shown in Table~\ref{tab:tspmd_1LKH_longer}, indicating that LKH can generate better solutions than ALNS under longer runtime. Additionally, the improved solution reduces Dcus from 27 to 26, suggesting that fewer drone assignments may derive more efficient truck-and-drone routing coordination.

\begin{table}[H]
    \centering
    % \tiny
    % \scriptsize
    \footnotesize  % Slightly smaller to better fit
    % \small
    % \normalsize  % (default size) 
    % \large
    \caption{TSPmD Comparison of ALNS and LKH (Config 1 4h)}
    \label{tab:tspmd_1LKH_longer}
    \begin{tabular*}{\textwidth}{@{\extracolsep{\fill}}l@{\;}r@{\;}r@{\;}r@{\;}r@{\;}r@{\;}r@{\;}r@{\;}r@{\;}r@{\;}r@{\;}r@{}}
    %%
    % \begin{tabularx}{\textwidth}{l r r r r r r r r r r r}
    % \begin{tabularx}{\textwidth}{l c c c c c c c c c c c}
    % \begin{tabularx}{\textwidth}{l *{11}{>{\centering\arraybackslash}X}}
        \toprule
        \multicolumn{2}{c}{} & \multicolumn{4}{c}{\textbf{ALNS}} & \multicolumn{5}{c}{\textbf{LKH (Config 1 4h)}} & \multicolumn{1}{c}{}\\
        \cmidrule(lr){3-6} \cmidrule(lr){7-11}
        \raisebox{1.5ex}[0pt][0pt]{Instance} & \raisebox{1.5ex}[0pt][0pt] {$N_D$} & $M^{\textit{ALNS}}$ & Dcus & \#D & Time(s) & $\text{M}^1$ & Dcus & \#D & Iter & Time(s) & \raisebox{1.5ex}[0pt][0pt] {$\%gap^{\textit{ALNS}}$} \\
        \midrule
\textbf{poi-40-4} & 39 & \textbf{179.55} & 17 & 5 & 300.02 & \textbf{150.6} & 26 & 5 & 139760 & 13287.69 & \textbf{16.12} \\
        \bottomrule
    \end{tabular*}
\end{table}

Then, we compapre LKH Config1 and Config2 in Tables~\ref{tab:tspmd_2LKH_small}--\ref{tab:tspmd_2LKH_medium}. The latter configuration enables the truck to visit nodes more than once, but this privilege is only utilized in some medium instances, highlighted in bold. The resulting gaps, $\%gap^{1,2} = 100 \cdot (\text{M}^1-\text{M}^2)/\text{M}^1$, vary across instances, implying that the added flexibility in Config2 does not always enhance LKH solutions within the 5min runtime. This variability may be due to the convergence behavior that is affected by the expanded solution space in LKH Config2.

% second comparison
\begin{table}[H]
    \centering
    % \tiny
    % \scriptsize
    \footnotesize  % Slightly smaller to better fit
    % \small
    % \normalsize  % (default size) 
    % \large
    \caption{TSPmD Comparison of LKH (Config 1) and LKH (Config 2) on small instances}
    \label{tab:tspmd_2LKH_small}
    \begin{tabular*}{\textwidth}{@{\extracolsep{\fill}}l@{\;}r@{\;}r@{\;}r@{\;}r@{\;}r@{\;}r@{\;}r@{\;}r@{\;}r@{\;}r@{\;}r@{\;}r@{}}
    %%
    % \begin{tabularx}{\textwidth}{l r r r r r r r r r r r r}
    % \begin{tabularx}{\textwidth}{l c c c c c c c c c c c c}
    % \begin{tabularx}{\textwidth}{l *{12}{>{\centering\arraybackslash}X}}
        \toprule
        \multicolumn{2}{c}{} & \multicolumn{5}{c}{\textbf{LKH (Config 1)}} & \multicolumn{5}{c}{\textbf{LKH (Config 2)}} & \multicolumn{1}{c}{}\\
        \cmidrule(lr){3-7} \cmidrule(lr){8-12}
        \raisebox{1.5ex}[0pt][0pt]{Instance} & \raisebox{1.5ex}[0pt][0pt] {$N_D$} & $\text{M}^1$ & Dcus & \#D & Iter & Time & $\text{M}^2$ & Dcus & \#D & Iter & Time & \raisebox{1.5ex}[0pt][0pt] {$\%gap^{1,2}$} \\
        \midrule
poi-7-1 & 6 & 90 & 3 & 3 & 126 & 0.4 & 90 & 3 & 3 & 189 & 0.4 & 0 \\
poi-7-2 & 6 & 65.2 & 4 & 4 & 55 & 0.15 & 65.2 & 4 & 4 & 50 & 0.1 & 0 \\
poi-7-3 & 6 & 73.8 & 1 & 1 & 229 & 0.71 & 73.8 & 1 & 1 & 109 & 0.46 & 0 \\
poi-7-4 & 6 & 81.2 & 3 & 3 & 187 & 0.48 & 81.2 & 3 & 3 & 1051 & 2.51 & 0 \\
poi-7-5 & 6 & 67.3 & 4 & 4 & 10 & 0.02 & 67.3 & 4 & 4 & 67 & 0.09 & 0 \\
poi-7-6 & 6 & 72.8 & 5 & 5 & 9 & 0.02 & 72.8 & 5 & 5 & 1 & 0.01 & 0 \\
poi-7-7 & 6 & 55 & 4 & 4 & 11 & 0.04 & 55 & 4 & 4 & 312 & 0.67 & 0 \\
poi-7-8 & 6 & 58.6 & 5 & 5 & 16 & 0.04 & 58.6 & 5 & 5 & 76 & 0.29 & 0 \\
poi-7-9 & 6 & 87.3 & 5 & 5 & 3 & 0.02 & 87.3 & 5 & 5 & 2 & 0.01 & 0 \\
poi-7-10 & 6 & 59.3 & 4 & 4 & 2309 & 5.51 & 59.3 & 4 & 4 & 2776 & 9.86 & 0 \\
poi-7-11 & 6 & 69.3 & 3 & 3 & 9 & 0.02 & 69.3 & 3 & 3 & 97 & 0.27 & 0 \\
poi-7-12 & 6 & 66.2 & 3 & 3 & 99 & 0.28 & 66.2 & 3 & 3 & 135 & 0.36 & 0 \\
poi-7-13 & 6 & 50.2 & 4 & 4 & 1 & 0.01 & 50.2 & 4 & 4 & 1 & 0.01 & 0 \\
poi-7-14 & 6 & 66.8 & 4 & 4 & 34 & 0.05 & 66.8 & 4 & 4 & 50 & 0.1 & 0 \\
poi-7-15 & 6 & 75.4 & 4 & 4 & 10 & 0.02 & 75.4 & 4 & 4 & 4 & 0.02 & 0 \\
poi-7-16 & 6 & 59.3 & 3 & 3 & 161 & 0.41 & 59.3 & 3 & 3 & 99 & 0.37 & 0 \\
poi-7-17 & 6 & 64.6 & 4 & 4 & 1 & 0.01 & 64.6 & 4 & 4 & 17 & 0.04 & 0 \\
poi-7-18 & 6 & 63.2 & 4 & 4 & 59 & 0.09 & 63.2 & 4 & 4 & 18 & 0.03 & 0 \\
poi-7-19 & 6 & 72.2 & 4 & 4 & 37631 & 15.31 & 72.2 & 4 & 4 & 100 & 0.14 & 0 \\
poi-7-20 & 6 & 71.2 & 5 & 5 & 66 & 0.09 & 71.2 & 5 & 5 & 20 & 0.04 & 0 \\
poi-7-21 & 6 & 47.5 & 3 & 3 & 148 & 0.5 & 47.5 & 3 & 3 & 963 & 3.12 & 0 \\
poi-7-22 & 6 & 66.6 & 4 & 4 & 18 & 0.04 & 66.6 & 4 & 4 & 109 & 0.32 & 0 \\
poi-7-23 & 6 & 80.2 & 4 & 4 & 17 & 0.06 & 80.2 & 4 & 4 & 120 & 0.2 & 0 \\
poi-7-24 & 6 & 46.7 & 3 & 3 & 639 & 2.39 & 46.7 & 3 & 3 & 80 & 0.21 & 0 \\
poi-7-25 & 6 & 68 & 3 & 3 & 298 & 0.59 & 68 & 3 & 3 & 175 & 0.44 & 0 \\
poi-10-1 & 9 & 90 & 4 & 4 & 5860 & 32.87 & 90 & 4 & 4 & 12130 & 67.47 & 0 \\
poi-10-2 & 9 & 65 & 6 & 5 & 13254 & 106.51 & 65 & 6 & 5 & 2667 & 19.4 & 0 \\
poi-10-3 & 9 & 73.9 & 4 & 4 & 530 & 3.09 & 73.9 & 4 & 4 & 430 & 2.49 & 0 \\
poi-10-4 & 9 & 77.5 & 5 & 5 & 274 & 1.21 & 77.5 & 5 & 5 & 1195 & 8.62 & 0 \\
poi-10-5 & 9 & 65.8 & 7 & 5 & 119 & 0.69 & 65.8 & 7 & 5 & 237 & 1.4 & 0 \\
poi-10-6 & 9 & 77.5 & 6 & 5 & 185 & 0.96 & 77.5 & 6 & 5 & 171 & 1 & 0 \\
poi-10-7 & 9 & 49.6 & 5 & 5 & 488 & 2.83 & 49.6 & 5 & 5 & 597 & 3.3 & 0 \\
poi-10-8 & 9 & 57.2 & 4 & 4 & 219 & 1.12 & 57.2 & 4 & 4 & 991 & 5.31 & 0 \\
poi-10-9 & 9 & 89.5 & 8 & 5 & 330 & 1.8 & 89.5 & 8 & 5 & 434 & 2.37 & 0 \\
poi-10-10 & 9 & 56.9 & 4 & 4 & 473 & 2.77 & 56.9 & 4 & 4 & 865 & 4.56 & 0 \\
poi-10-11 & 9 & 65 & 4 & 4 & 586 & 2.89 & 65 & 4 & 4 & 3024 & 16.25 & 0 \\
poi-10-12 & 9 & 66.2 & 4 & 4 & 400 & 1.99 & 66.2 & 4 & 4 & 7427 & 33.03 & 0 \\
poi-10-13 & 9 & 52.4 & 7 & 5 & 7704 & 34.82 & 52.4 & 7 & 5 & 668 & 4.14 & 0 \\
poi-10-14 & 9 & 65.1 & 6 & 5 & 1373 & 9.67 & 65.1 & 6 & 5 & 2323 & 15.8 & 0 \\
poi-10-15 & 9 & 64.9 & 4 & 4 & 714 & 3.15 & 64.9 & 4 & 4 & 62 & 0.39 & 0 \\
poi-10-16 & 9 & 64.7 & 6 & 5 & 307 & 2.25 & 64.7 & 6 & 5 & 25855 & 151.63 & 0 \\
poi-10-17 & 9 & 64.6 & 6 & 5 & 5158 & 33.46 & 64.6 & 6 & 5 & 4599 & 36.59 & 0 \\
poi-10-18 & 9 & 63.2 & 7 & 5 & 27804 & 163.16 & 63.2 & 7 & 5 & 673 & 4.12 & 0 \\
poi-10-19 & 9 & 73 & 6 & 5 & 5797 & 40.14 & 73 & 6 & 5 & 8829 & 62.72 & 0 \\
poi-10-20 & 9 & 74.5 & 7 & 5 & 1955 & 11.9 & 74.5 & 7 & 5 & 4643 & 28.18 & 0 \\
poi-10-21 & 9 & 47.5 & 6 & 5 & 120 & 0.43 & 47.6 & 7 & 5 & 99 & 0.43 & -0.21 \\
poi-10-22 & 9 & 61.9 & 4 & 4 & 739 & 4.24 & 61.9 & 4 & 4 & 37873 & 223.01 & 0 \\
poi-10-23 & 9 & 80.2 & 5 & 5 & 305 & 1.67 & 80.2 & 5 & 5 & 24822 & 142.74 & 0 \\
poi-10-24 & 9 & 73.9 & 4 & 4 & 991 & 5.27 & 73.9 & 4 & 4 & 802 & 4.91 & 0 \\
poi-10-25 & 9 & 74.9 & 7 & 5 & 166 & 0.61 & 74.9 & 7 & 5 & 153 & 0.79 & 0 \\
        \bottomrule
    \end{tabular*}
\end{table}

\begin{table}[H]
    \centering
    % \tiny
    % \scriptsize
    \footnotesize  % Slightly smaller to better fit
    % \small
    % \normalsize  % (default size) 
    % \large
    \caption{TSPmD Comparison of LKH (Config 1) and LKH (Config 2) on medium instances}
    \label{tab:tspmd_2LKH_medium}
    \begin{tabular*}{\textwidth}{@{\extracolsep{\fill}}l@{\;}r@{\;}r@{\;}r@{\;}r@{\;}r@{\;}r@{\;}r@{\;}r@{\;}r@{\;}r@{\;}r@{\;}r@{}}
    %%
    % \begin{tabularx}{\textwidth}{l r r r r r r r r r r r r}
    % \begin{tabularx}{\textwidth}{l c c c c c c c c c c c c}
    % \begin{tabularx}{\textwidth}{l *{12}{>{\centering\arraybackslash}X}}
        \toprule
        \multicolumn{2}{c}{} & \multicolumn{5}{c}{\textbf{LKH (Config 1)}} & \multicolumn{5}{c}{\textbf{LKH (Config 2)}} & \multicolumn{1}{c}{}\\
        \cmidrule(lr){3-7} \cmidrule(lr){8-12}
        \raisebox{1.5ex}[0pt][0pt]{Instance} & \raisebox{1.5ex}[0pt][0pt] {$N_D$} & $\text{M}^1$ & Dcus & \#D & Iter & Time(s) & $\text{M}^2$ & Dcus & \#D & Iter & Time(s) & \raisebox{1.5ex}[0pt][0pt] {$\%gap^{1,2}$} \\
        \midrule
\textbf{poi-20-1} & 19 & \textbf{106.4} & 12 & 5 & 1642 & 35.41 & \textbf{97.8 revisit} & 15 & 5 & 7506 & 192.86 & \textbf{8.08} \\
poi-20-2 & 19 & 103 & 13 & 5 & 6365 & 96.96 & 98.7 & 13 & 5 & 13442 & 271.67 & 4.17 \\
poi-20-3 & 19 & 108.4 & 14 & 5 & 16139 & 253.29 & 108.4 & 14 & 5 & 10698 & 197.66 & 0 \\
poi-20-4 & 19 & 89.6 & 12 & 5 & 12208 & 238.58 & 97.1 & 13 & 5 & 2177 & 52.52 & -8.37 \\
poi-20-5 & 19 & 94.1 & 15 & 5 & 10451 & 226.77 & 93.2 & 13 & 5 & 10354 & 245.5 & 0.96 \\
poi-20-6 & 19 & 85.5 & 14 & 5 & 15968 & 276.3 & 85.7 & 13 & 5 & 4266 & 110.6 & -0.23 \\
poi-20-7 & 19 & 111.1 & 11 & 5 & 1804 & 27 & 104.5 & 11 & 5 & 1927 & 34.73 & 5.94 \\
poi-20-8 & 19 & 118.7 & 10 & 5 & 7080 & 143.1 & 118.8 & 9 & 5 & 9110 & 269.49 & -0.08 \\
poi-20-9 & 19 & 113.9 & 12 & 5 & 15617 & 255.69 & 117.1 & 10 & 5 & 4644 & 83.66 & -2.81 \\
poi-20-10 & 19 & 97.5 & 14 & 5 & 11607 & 217.25 & 94.2 & 13 & 5 & 6341 & 157.83 & 3.38 \\
poi-20-11 & 19 & 109.3 & 13 & 5 & 6652 & 116.48 & 116.9 & 14 & 5 & 10785 & 282.04 & -6.95 \\
\textbf{poi-20-12} & 19 & \textbf{99.9} & 13 & 5 & 6352 & 93.69 & \textbf{86.8 revisit} & 12 & 5 & 7609 & 181.47 & \textbf{13.11} \\
poi-20-13 & 19 & 88 & 13 & 5 & 13142 & 297.71 & 99 & 14 & 5 & 11242 & 280.32 & -12.5 \\
poi-20-14 & 19 & 103.7 & 12 & 5 & 9621 & 183.33 & 103.4 & 14 & 5 & 2040 & 34.92 & 0.29 \\
poi-20-15 & 19 & 100 & 13 & 5 & 676 & 12.71 & 98.7 & 13 & 5 & 4462 & 108.95 & 1.3 \\
poi-20-16 & 19 & 99.2 & 14 & 5 & 14084 & 231.6 & 98 & 15 & 5 & 3747 & 106.47 & 1.21 \\
poi-20-17 & 19 & 100.1 & 11 & 5 & 7958 & 168.22 & 97 & 10 & 5 & 3909 & 77.16 & 3.1 \\
\textbf{poi-20-18} & 19 & \textbf{109.5} & 13 & 5 & 15156 & 277.12 & \textbf{110.7 revisit} & 13 & 5 & 3836 & 79.92 & \textbf{-1.1} \\
\textbf{poi-20-19} & 19 & \textbf{93.9} & 15 & 5 & 4739 & 84.35 & \textbf{93.9 revisit} & 14 & 5 & 3064 & 45.55 & \textbf{0} \\
poi-20-20 & 19 & 120.9 & 13 & 5 & 8623 & 132.6 & 124 & 12 & 5 & 3111 & 60.8 & -2.56 \\
poi-20-21 & 19 & 99.6 & 12 & 5 & 2407 & 41.12 & 94.8 & 13 & 5 & 13372 & 286.42 & 4.82 \\
poi-20-22 & 19 & 98.1 & 11 & 5 & 9452 & 204.49 & 96.6 & 12 & 5 & 3987 & 94.24 & 1.53 \\
poi-20-23 & 19 & 107.5 & 11 & 5 & 5480 & 98.55 & 108.2 & 10 & 5 & 5990 & 119.75 & -0.65 \\
poi-20-24 & 19 & 97.8 & 13 & 5 & 4250 & 70.94 & 92.1 & 13 & 5 & 878 & 16.86 & 5.83 \\
\textbf{poi-20-25} & 19 & \textbf{98.1} & 14 & 5 & 963 & 12.48 & \textbf{94 revisit} & 14 & 5 & 9092 & 268.02 & \textbf{4.18} \\
poi-30-1 & 29 & 127 & 16 & 5 & 8925 & 285.13 & 152.5 & 18 & 5 & 7905 & 260.59 & -20.08 \\
poi-30-2 & 29 & 118.5 & 20 & 5 & 7652 & 284.62 & 125 & 22 & 5 & 5658 & 273.25 & -5.49 \\
\textbf{poi-30-3} & 29 & \textbf{4.18} & 20 & 5 & 6341 & 232.16 & \textbf{121 revisit} & 23 & 5 & 1722 & 78.37 & \textbf{3.12} \\
poi-30-4 & 29 & 128.5 & 20 & 5 & 3962 & 112.69 & 131.1 & 21 & 5 & 7197 & 264.99 & -2.02 \\
\textbf{poi-30-5} & 29 & \textbf{127.3} & 19 & 5 & 4593 & 198.29 & \textbf{157.7 revisit} & 21 & 5 & 3113 & 118.19 & \textbf{-23.88} \\
poi-30-6 & 29 & 128.7 & 19 & 5 & 7084 & 268.87 & 127.8 & 19 & 5 & 5333 & 167.39 & 0.7 \\
poi-30-7 & 29 & 108.8 & 21 & 5 & 7335 & 198.3 & 113.4 & 20 & 5 & 3511 & 109.8 & -4.23 \\
poi-30-8 & 29 & 141.3 & 22 & 5 & 10667 & 297.09 & 122.2 & 20 & 5 & 5638 & 246.58 & 13.52 \\
\textbf{poi-30-9} & 29 & \textbf{124.3} & 20 & 5 & 5880 & 266.11 & \textbf{124.2 revisit} & 19 & 5 & 4705 & 260.47 & \textbf{0.08} \\
poi-30-10 & 29 & 142.5 & 21 & 5 & 6740 & 295.46 & 139.6 & 22 & 5 & 6890 & 258.34 & 2.04 \\
poi-30-11 & 29 & 133.3 & 19 & 5 & 8533 & 293.55 & 127.5 & 21 & 5 & 7680 & 272.34 & 4.35 \\
poi-30-12 & 29 & 132.2 & 20 & 5 & 9864 & 281.44 & 140.8 & 19 & 5 & 4031 & 160.15 & -6.51 \\
poi-30-13 & 29 & 145.9 & 19 & 5 & 4932 & 167.33 & 147.2 & 19 & 5 & 6778 & 208.25 & -0.89 \\
\textbf{poi-30-14} & 29 & \textbf{130.8} & 23 & 5 & 4965 & 277.04 & \textbf{143 revisit} & 24 & 5 & 3470 & 207.84 & \textbf{-9.33} \\
poi-30-15 & 29 & 112.2 & 20 & 5 & 7168 & 287.62 & 128.4 & 21 & 5 & 6039 & 212.82 & -14.44 \\
\textbf{poi-30-16} & 29 & \textbf{116.6} & 21 & 5 & 5199 & 225.32 & \textbf{138.4 revisit} & 24 & 5 & 11518 & 294 & \textbf{-18.7} \\
poi-30-17 & 29 & 128.4 & 18 & 5 & 7010 & 263.85 & 130.3 & 19 & 5 & 4991 & 185.91 & -1.48 \\
\textbf{poi-30-18} & 29 & \textbf{121.5} & 23 & 5 & 1929 & 102.92 & \textbf{116 revisit} & 23 & 5 & 4081 & 285.69 & \textbf{4.53} \\
\textbf{poi-30-19} & 29 & \textbf{142.9} & 19 & 5 & 4262 & 186.43 & \textbf{155 revisit} & 22 & 5 & 5851 & 273.28 & \textbf{-8.47} \\
poi-30-20 & 29 & 122.2 & 18 & 5 & 3986 & 173.89 & 120.9 & 22 & 5 & 6093 & 218.73 & 1.06 \\
\textbf{poi-30-21} & 29 & \textbf{150.8} & 23 & 5 & 5681 & 222.18 & \textbf{118.4 revisit} & 20 & 5 & 6841 & 290.29 & \textbf{21.49} \\
poi-30-22 & 29 & 142.2 & 20 & 5 & 4270 & 167.61 & 158.1 & 21 & 5 & 1941 & 114.95 & -11.18 \\
poi-30-23 & 29 & 117 & 19 & 5 & 5308 & 243.85 & 135.5 & 19 & 5 & 5711 & 260.17 & -15.81 \\
poi-30-24 & 29 & 129.6 & 20 & 5 & 8142 & 299.73 & 121.4 & 23 & 5 & 6841 & 264.65 & 6.33 \\
poi-30-25 & 29 & 155.8 & 19 & 5 & 5945 & 229.09 & 156.3 & 18 & 5 & 7509 & 287.64 & -0.32 \\
        \bottomrule
        \multicolumn{13}{r}{Continued on next page}
    \end{tabular*}
\end{table}

\clearpage % Force a page break here

\begin{table}[!t] % Force placement at the absolute top
    \vspace*{0pt} % Remove any vertical space at the top
    \centering
    % \tiny
    % \scriptsize
    \footnotesize  % Slightly smaller to better fit
    % \small
    % \normalsize  % (default size) 
    % \large
    \refstepcounter{table}
    \addtocounter{table}{-1} % Keep the table number the same
    \makeatletter
    \def\@captiontype{table}
    {\@makecaption{\tablename~\thetable}{(continued)}}
    \vspace{0.5em}
    \makeatother
    \begin{tabular*}{\textwidth}{@{\extracolsep{\fill}}l@{\;}r@{\;}r@{\;}r@{\;}r@{\;}r@{\;}r@{\;}r@{\;}r@{\;}r@{\;}r@{\;}r@{\;}r@{}}
    %%
    % \begin{tabularx}{\textwidth}{l r r r r r r r r r r r r}
    % \begin{tabularx}{\textwidth}{l c c c c c c c c c c c c}
    % \begin{tabularx}{\textwidth}{l *{12}{>{\centering\arraybackslash}X}}
        \toprule
        \multicolumn{2}{c}{} & \multicolumn{5}{c}{\textbf{LKH (Config 1)}} & \multicolumn{5}{c}{\textbf{LKH (Config 2)}} & \multicolumn{1}{c}{}\\
        \cmidrule(lr){3-7} \cmidrule(lr){8-12}
        \raisebox{1.5ex}[0pt][0pt]{Instance} & \raisebox{1.5ex}[0pt][0pt] {$N_D$} & $\text{M}^1$ & Dcus & \#D & Iter & Time(s) & $\text{M}^2$ & Dcus & \#D & Iter & Time(s) & \raisebox{1.5ex}[0pt][0pt] {$\%gap^{1,2}$} \\
        \midrule
\textbf{poi-40-1} & 39 & \textbf{158.4} & 27 & 5 & 3060 & 203.2 & \textbf{180.1 revisit} & 30 & 5 & 4809 & 263.39 & \textbf{-13.7} \\
\textbf{poi-40-2} & 39 & \textbf{165.1} & 28 & 5 & 3712 & 242.52 & \textbf{183.1 revisit} & 32 & 5 & 3894 & 282.72 & \textbf{-10.9} \\
\textbf{poi-40-3} & 39 & \textbf{178.2} & 29 & 5 & 9707 & 283.86 & \textbf{184.5 revisit} & 26 & 5 & 5118 & 289.27 & \textbf{-3.54} \\
\textbf{poi-40-4} & 39 & \textbf{186.1} & 27 & 5 & 2757 & 216.73 & \textbf{172.9 revisit} & 28 & 5 & 4175 & 281.57 & \textbf{7.09} \\
\textbf{poi-40-5} & 39 & \textbf{178.6} & 28 & 5 & 3237 & 200.05 & \textbf{198.2 revisit} & 30 & 5 & 8677 & 292.41 & \textbf{-10.97} \\
\textbf{poi-40-6} & 39 & \textbf{154.8} & 28 & 5 & 2923 & 237.9 & \textbf{174.3 revisit} & 31 & 5 & 4207 & 247.49 & \textbf{-12.6} \\
\textbf{poi-40-7} & 39 & \textbf{153.8} & 27 & 5 & 6652 & 278.1 & \textbf{194 revisit} & 29 & 5 & 4854 & 215.28 & \textbf{-26.14} \\
\textbf{poi-40-8} & 39 & \textbf{182.1} & 30 & 5 & 4291 & 269.63 & \textbf{189.9 revisit} & 31 & 5 & 4159 & 215.75 & \textbf{-4.28} \\
poi-40-9 & 39 & 166.1 & 29 & 5 & 6585 & 281.77 & 198.3 & 31 & 5 & 2784 & 287.88 & -19.39 \\
\textbf{poi-40-10} & 39 & \textbf{160.9} & 29 & 5 & 4806 & 298.25 & \textbf{198.2 revisit} & 34 & 5 & 5009 & 286.74 & \textbf{-23.18} \\
poi-40-11 & 39 & 168.2 & 26 & 5 & 6275 & 286.17 & 183.7 & 29 & 5 & 3506 & 183.94 & -9.22 \\
\textbf{poi-40-12} & 39 & \textbf{180.6} & 30 & 5 & 4434 & 282.35 & \textbf{186.6 revisit} & 32 & 5 & 5458 & 264.88 & \textbf{-3.32} \\
\textbf{poi-40-13} & 39 & \textbf{173.5} & 24 & 5 & 4926 & 298.79 & \textbf{175.1 revisit} & 31 & 5 & 2238 & 136.79 & \textbf{-0.92} \\
\textbf{poi-40-14} & 39 & \textbf{154} & 27 & 5 & 5418 & 287.13 & \textbf{187 revisit} & 28 & 5 & 3267 & 279.12 & \textbf{-21.43} \\
\textbf{poi-40-15} & 39 & \textbf{183.6} & 29 & 5 & 3721 & 267.12 & \textbf{176.6 revisit} & 31 & 5 & 2850 & 216.33 & \textbf{3.81} \\
\textbf{poi-40-16} & 39 & \textbf{161.1} & 27 & 5 & 3890 & 228.55 & \textbf{177.3 revisit} & 26 & 5 & 4491 & 296.46 & \textbf{-10.06} \\
poi-40-17 & 39 & 188.9 & 30 & 5 & 7845 & 248.01 & 166.9 & 28 & 5 & 3228 & 177.65 & 11.65 \\
\textbf{poi-40-18} & 39 & \textbf{195.9} & 28 & 5 & 3651 & 95.91 & \textbf{203.8 revisit} & 30 & 5 & 4419 & 260.33 & \textbf{-4.03} \\
poi-40-19 & 39 & 175.3 & 31 & 5 & 3144 & 244.05 & 180.3 & 34 & 5 & 4258 & 299.92 & -2.85 \\
poi-40-20 & 39 & 176.8 & 26 & 5 & 4044 & 249.82 & 190.1 & 32 & 5 & 4850 & 294 & -7.52 \\
\textbf{poi-40-21} & 39 & \textbf{165.6} & 29 & 5 & 5632 & 275.79 & \textbf{169 revisit} & 27 & 5 & 3097 & 214.32 & \textbf{-2.05} \\
\textbf{poi-40-22} & 39 & \textbf{177.3} & 28 & 5 & 4124 & 243.65 & \textbf{185.7 revisit} & 29 & 5 & 2979 & 297.14 & \textbf{-4.74} \\
poi-40-23 & 39 & 190.4 & 28 & 5 & 6172 & 287.4 & 191.8 & 29 & 5 & 4911 & 288.29 & -0.74 \\
poi-40-24 & 39 & 167.4 & 29 & 5 & 4643 & 193.15 & 191 & 28 & 5 & 4384 & 257.1 & -14.1 \\
\textbf{poi-40-25} & 39 & \textbf{171.9} & 27 & 5 & 4330 & 297.94 & \textbf{206.6 revisit} & 34 & 5 & 2959 & 279.66 & \textbf{-20.19} \\
        \bottomrule
    \end{tabular*}
\end{table}

In brief, within one execution in 5min shown in Tables~\ref{tab:tspmd_1LKH_small}--\ref{tab:tspmd_1LKH_medium}, LKH Config1 identifies optimal solutions for small Poikonen instances, and outperforms the ALNS algorithm (five 5min runs) from \cite{tu2018traveling} on all Poikonen instances except one poi-40-4. However, this exception is resolved when LKH is given an extended runtime of 4h, as shown in Table~\ref{tab:tspmd_1LKH_longer}.

\subsection{VRP-D}
% For each problem instance, we can obtain the customer demands ($q_i$) and locations ($(x_i, y_i)$) by reading its ``.txt'' instance document. They are generated by randomly uniform distribution and provided by \cite{sacramento2019adaptive}. Input: $U(0,2.27)$ kg for drone-eligible (86\%), $U(2.27,68)$ kg for truck-only (14\%)

Aligned to the setup in Section 4.3.2, we benchmark LKH Config1 against the ALNS results reported in \cite{sacramento2019adaptive}, using the problem instances they generated, i.e., ``Set 6": 48 min-cost instances of 6/10/12/20 customers. We also apply LKH Config2 to the same set to investigate the impact of permitting node revisit. Notably, this VRP-D variant involves constraints 1-4 described in \hyperref[sec:introduction]{Introduction}: drone endurance ($e$), drone operation time ($D_{L/R}$), service time ($\textit{Serv}$), and drone capacity ($Q_D$).

Specifically, \cite{sacramento2019adaptive} reports the optimal solution in CPLEX, and the best solution found by ALNS over 10 independent 5min runs. Following their parameters provided in Table~\ref{tab:vrpd_parameter}, only 86\% of the input customer demands fall below the drone payload limit of 5 kg ($Q_D$), making them eligible for drone delivery. Each truck has a payload capacity ($Q_T$) of 1300 kg when carrying drones and is constrained by a maximum operational time ($T_\textit{max}$) of 8h. The truck and drone speeds are set to 35 mph ($v_T$) and 50 mph ($v_D$), respectively. For customers, making deliveries requires 2min by truck ($\textit{Serv}_T$) and 1min by drone ($\textit{Serv}_D$). $D_{L/R}$ implies a drone needs 1 minute for preparation and takeoff at a launch point, and 1 minute for recovery and recharging at a retrieval point. The drone battery life ($e$) is capped at 30 minutes \cite{kharpal2016firm}.

In addition, \cite{sacramento2019adaptive} use a truck fuel price of 1.13 €/L with a consumption rate of 0.07 L/km, converting miles to kilometers using a factor of 1.61 km/mile. The drone-to-truck cost ratio ($\alpha$) is set at 0.1, reflecting the lower operational cost of drones compared to trucks on equivalent distances. Based on these parameters, we estimate the operating costs of 0.127351 €/mile and 0.7428808 €/min for trucks, versus just 0.0127351 €/mile and 0.01061258 €/min for drones, as listed in the last four rows.

% Distances ($d_{ij}$) are calculated as Euclidean distances in miles between any two nodes. With a maximum drone payload of 5 kg ($Q_D$), customers with demands exceeding this threshold ($q_i > Q_D$) must be serviced by trucks. Besides, the drone's battery life $e$ is limited to 30 minutes \cite{kharpal2016firm}, and the truck's capacity $Q_T$ when carrying drones is 1300 kg \cite{sacramento2019adaptive}.

% Moreover, stricter operational constraints than typical real-world conditions are imposed in this subsection to rigorously test our approach's viability: drone battery life is limited to 10 minutes (rather than 30 minutes \cite{kharpal2016firm}) and truck capacity is capped at 100 kg (instead of 1300 kg \cite{sacramento2019adaptive}). 

\begin{table}[H]
    \centering
    % \tiny
    \scriptsize
    % \footnotesize  % Slightly smaller to better fit
    % \small
    % \normalsize  % (default size) 
    % \large
    \renewcommand{\arraystretch}{1.1} 
    \caption{Parameter definitions and values}
    \label{tab:vrpd_parameter}
    \begin{tabularx}{\textwidth}{>{\centering\arraybackslash}p{3.0cm} >{\centering\arraybackslash}p{0.6cm} >{\centering\arraybackslash}p{6.9cm} >{\centering\arraybackslash}X} 
        \toprule
        \textbf{Parameters} & \textbf{Notation} & \textbf{Value} & \textbf{Reference} \\
        \midrule
Customer demand & $q_i$ & ``.txt". (86\% drone-eligible) & \cite{sacramento2019adaptive} \\
Location & $(x_i, y_i)$ & ``.txt". (distance $d_{ij}$ measured in mile) & - \\
Route duration & $T_\textit{max}$ & Standard working hours: 8 hours & \cite{sacramento2019adaptive} \\
Truck capacity & $Q_T$ & With drones 1300 kg, without 1400 kg & \cite{sacramento2019adaptive} \\
Drone capacity & $Q_D$ & 5 kg & \cite{sacramento2019adaptive} \\
Truck speed & $v_T$ & 35 mph & \cite{ponza2016optimization} \\
Drone speed & $v_D$ & 50 mph & \cite{sacramento2019adaptive} \\
Truck service time & $\textit{Serv}_T$ & 2 minutes & \cite{sacramento2019adaptive} \\
Drone service time & $\textit{Serv}_D$ & 1 minute & \cite{sacramento2019adaptive} \\
Drone L/R time & $D_{L/R}$ & 1 minute & \cite{murray2015flying} \\
Flight endurance & $e$ & 30 minutes & \cite{kharpal2016firm} \\
\midrule
\multicolumn{4}{c}{Related to cost calculation.} \\
Fuel price & $c_1$ & 1.13 €/L & - \\
Fuel consumption & $c_2$ & 0.07 L/km & - \\
Converter & $c_3$ & 1.61 km/mile & - \\
Cost factor & $\alpha$ & 0.1 & \cite{kharpal2016firm} \\
Truck cost1 &  & 0.127351 €/mile & - \\
Truck cost2 &  & 0.7428808 €/min & - \\
Drone cost1 &  & 0.0127351 €/mile & - \\
Drone cost2 &  & 0.01061258 €/min & - \\
        \bottomrule
    \end{tabularx}
\end{table}

Before analyzing the comparison tables: Tables~\ref{tab:vrpd_1LKH}--\ref{tab:vrpd_1LKH_long} between ALNS and LKH Config1, and Table~\ref{tab:vrpd_2LKH} between LKH Config1 and Config2, we exhibit several typical instances as follows.

\subsubsection{Instances}

\begin{figure}[H]
    \centering
    \setlength{\tabcolsep}{1pt} % Space between columns
    \begin{tabular}{ccc}        \includegraphics[width=0.33\textwidth]{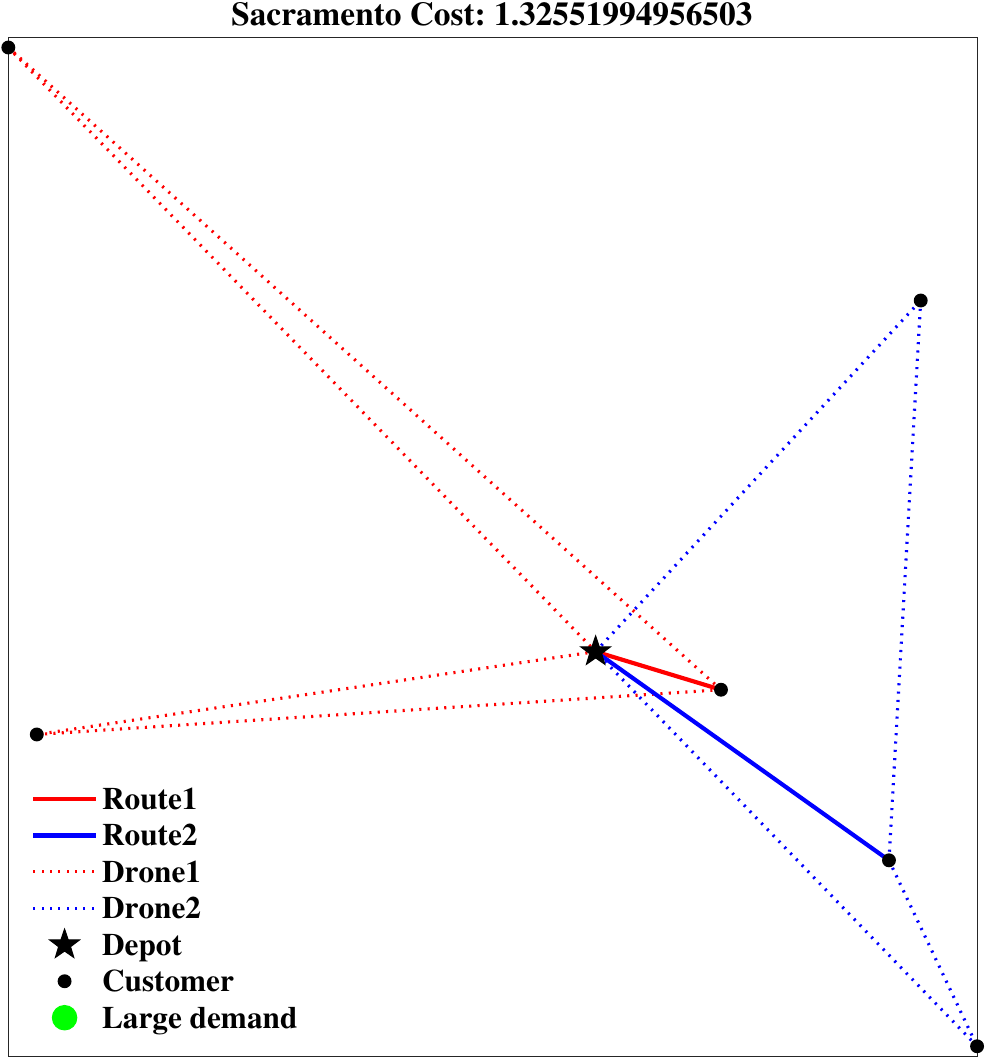} &        
    \includegraphics[width=0.33\textwidth]{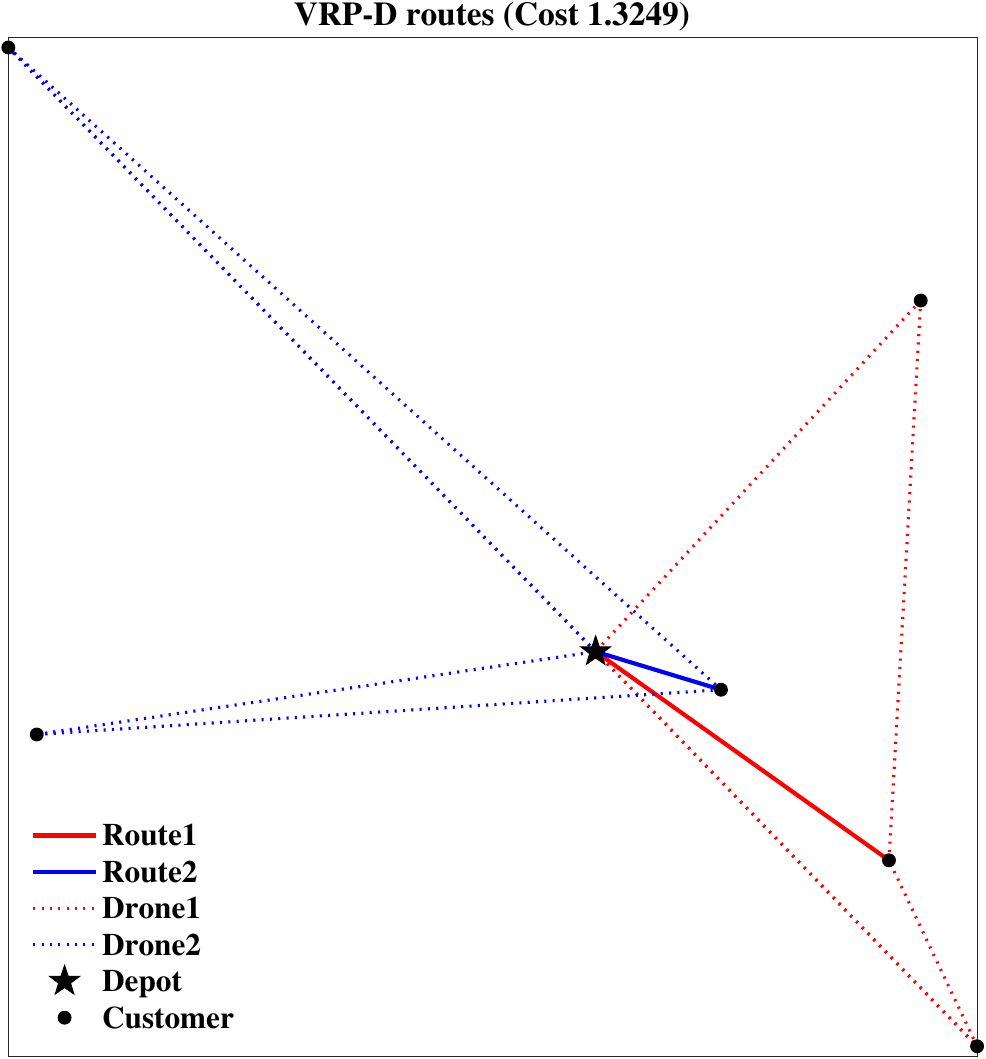} &        
    \includegraphics[width=0.33\textwidth]{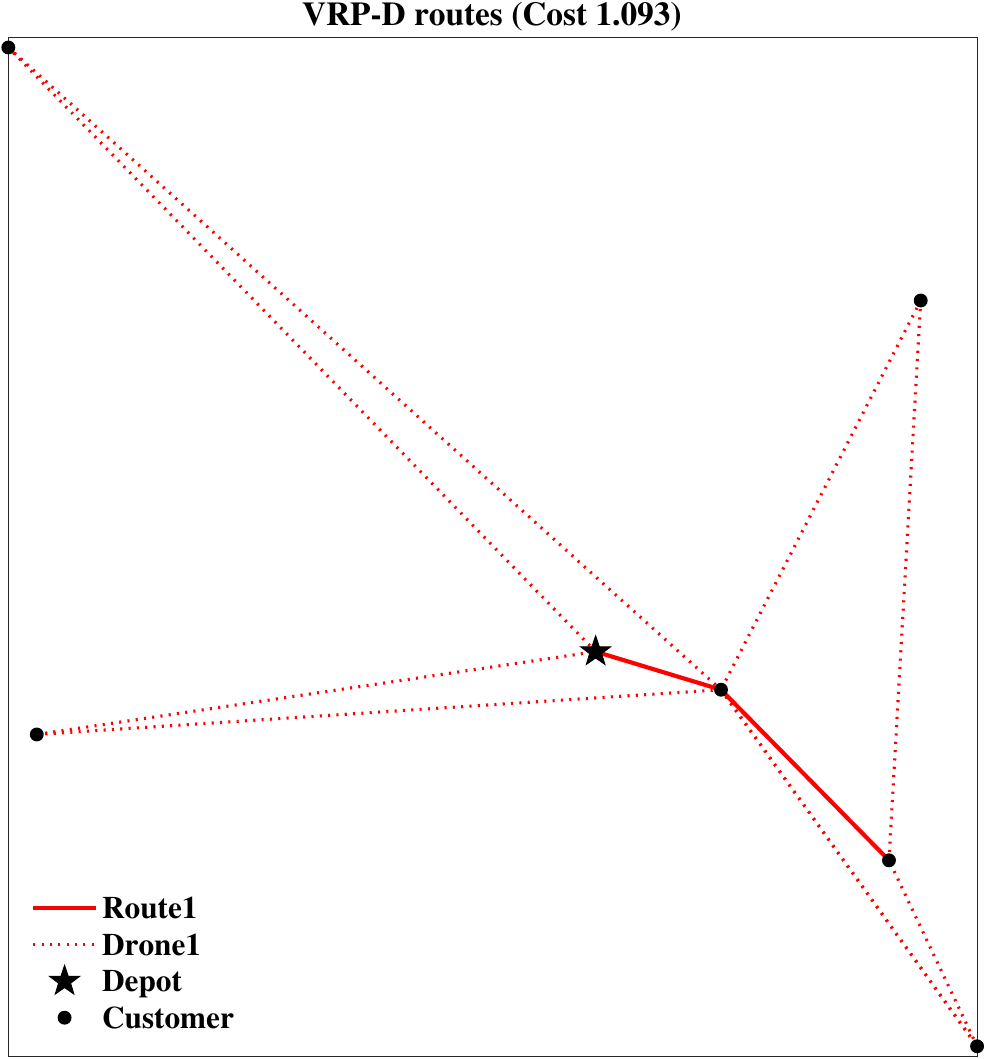} \\
        \multicolumn{1}{c}{(a) ALNS.} & \multicolumn{1}{c}{(b) LKH Config1.} & \multicolumn{1}{c}{(c) LKH Config2.}
    \end{tabular}
    \caption{VRPD tours on instance 6.10.3.}
    \label{fig:vrpd_figcompare1}
\end{figure}

Firstly, according to the comparison between two configurations in Table~\ref{tab:vrpd_2LKH}, only the bold instances utilize the node-revisit capability in LKH Config2. This flexibility enables the lowest cost of 1.093 for instance 6.10.3 as shown in Fig.~\ref{fig:vrpd_figcompare1}. Note that LKH Config1 generates an identical tour to ALNS, but their objective values differ slightly. This difference might be because we prepare integers in the input file \texttt{.drone} for LKH, as it requires integer inputs such as distances and coordinates.

% However, for some instances where revisiting nodes is utilized by Config2 to produce the best solution among the three configurations, Config1 may perform worse than ALNS due to the 5min runtime limit. 

% width=0.33\textwidth
% height=0.24\textheight

\begin{figure}[H]
    \centering
    \setlength{\tabcolsep}{1pt} % Space between columns
    \begin{tabular}{ccc}        \includegraphics[width=0.33\textwidth]{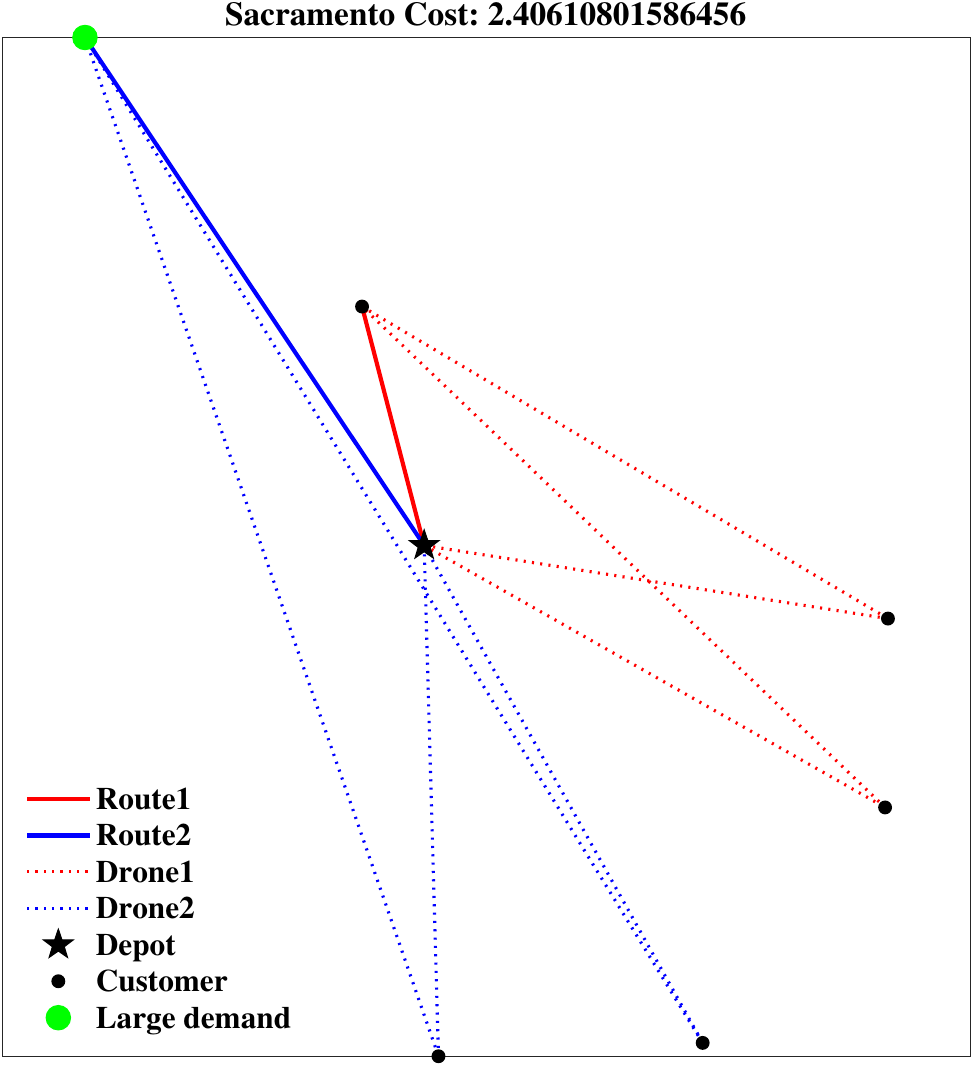} &        \includegraphics[width=0.33\textwidth]{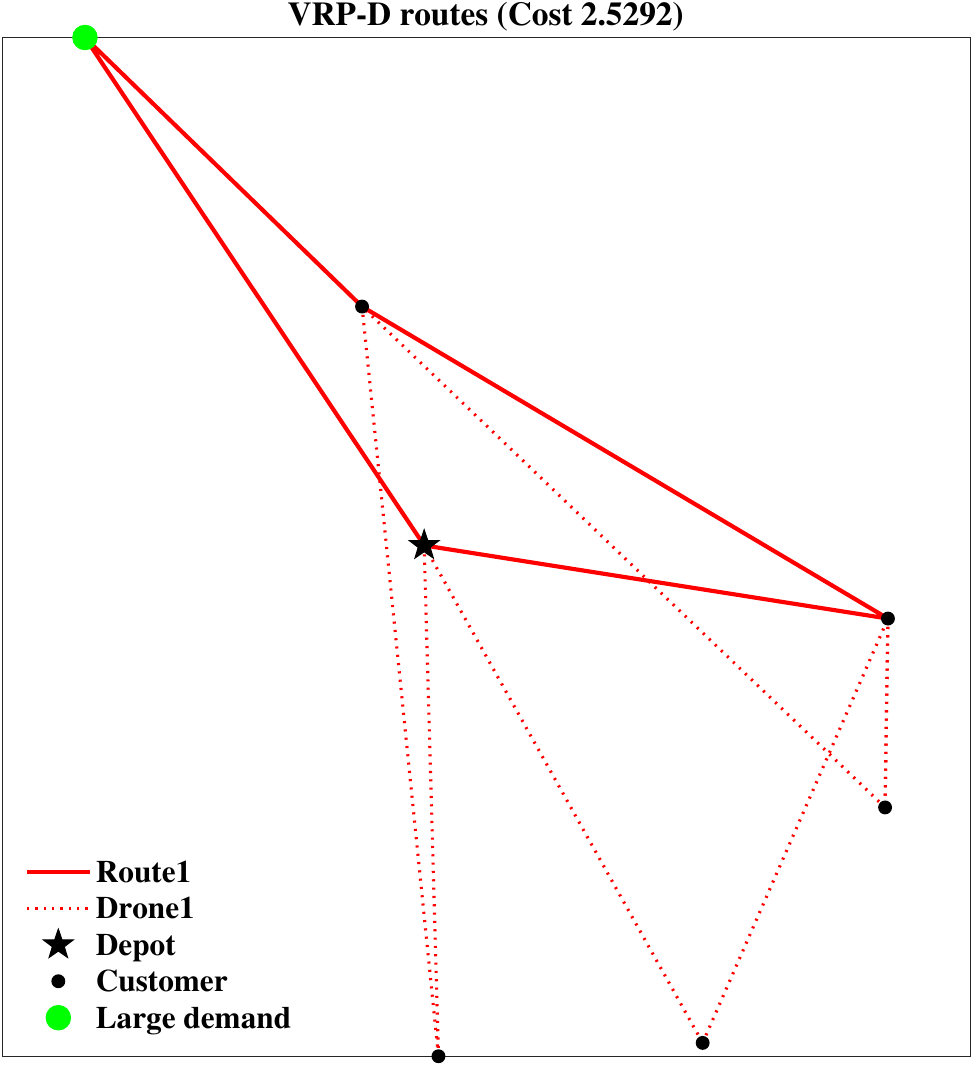} &        \includegraphics[width=0.33\textwidth]{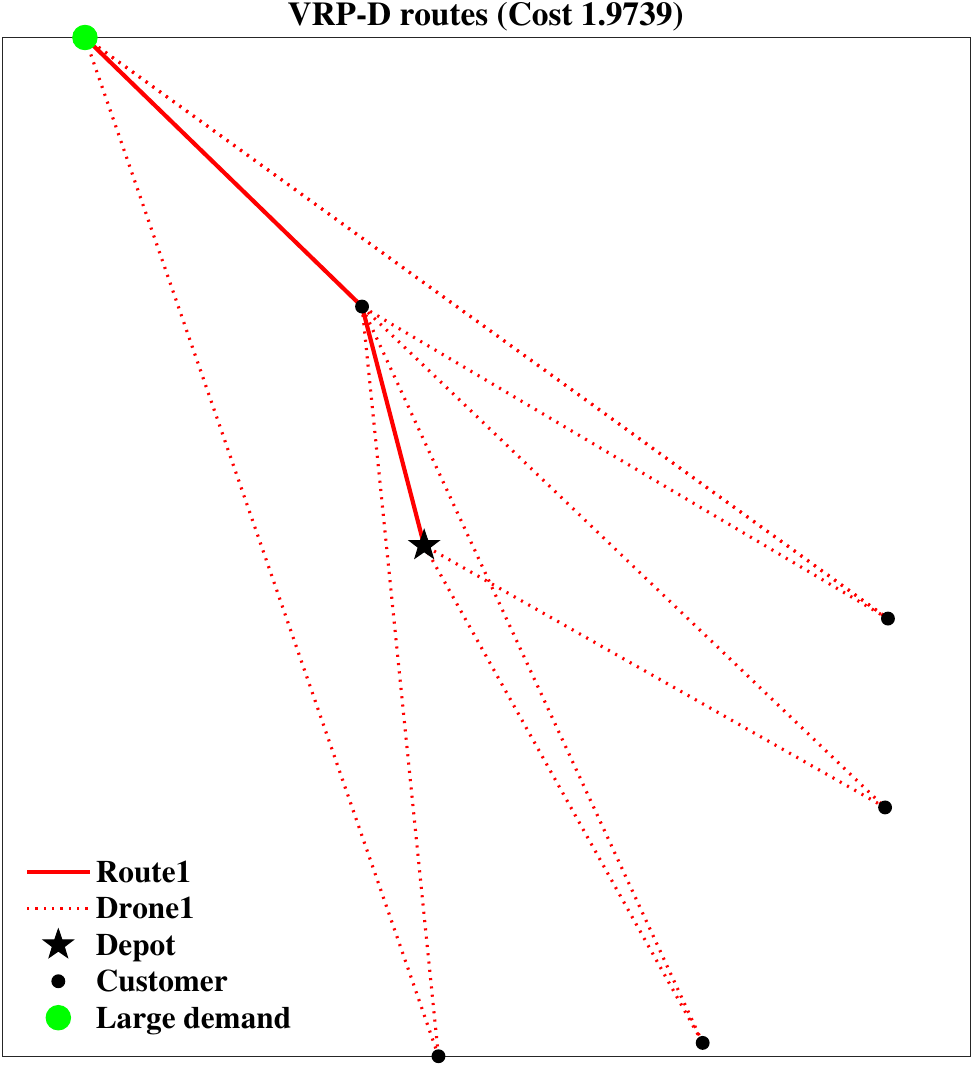} \\
        \multicolumn{1}{c}{(a) ALNS.} & \multicolumn{1}{c}{(b) LKH Config1.} & \multicolumn{1}{c}{(c) LKH Config2.}
    \end{tabular}
    \caption{VRPD tours on instance 6.10.1.}
    \label{fig:vrpd_figcompare4}
\end{figure}

\begin{figure}[H]
    \centering
    \setlength{\tabcolsep}{1pt} % Space between columns
    \begin{tabular}{ccc}        \includegraphics[width=0.33\textwidth]{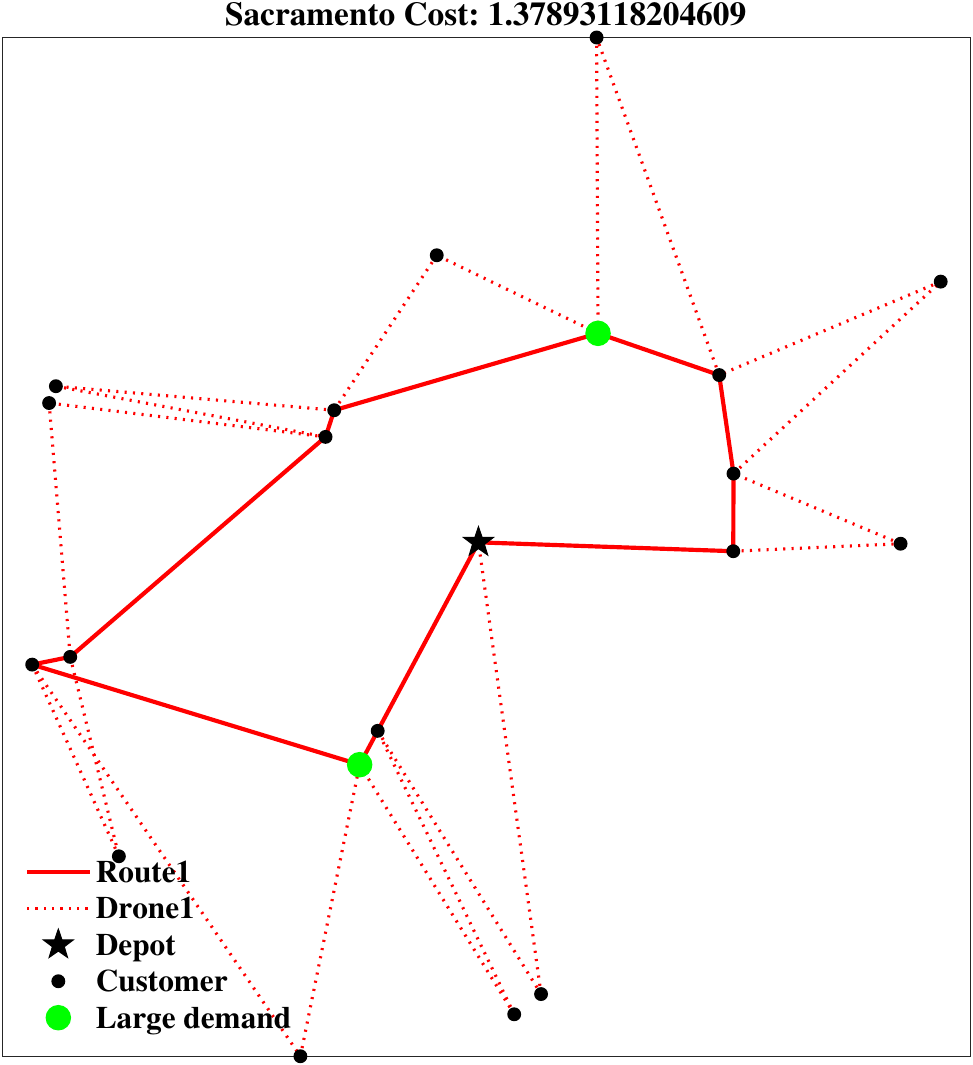} &        
    \includegraphics[width=0.33\textwidth]{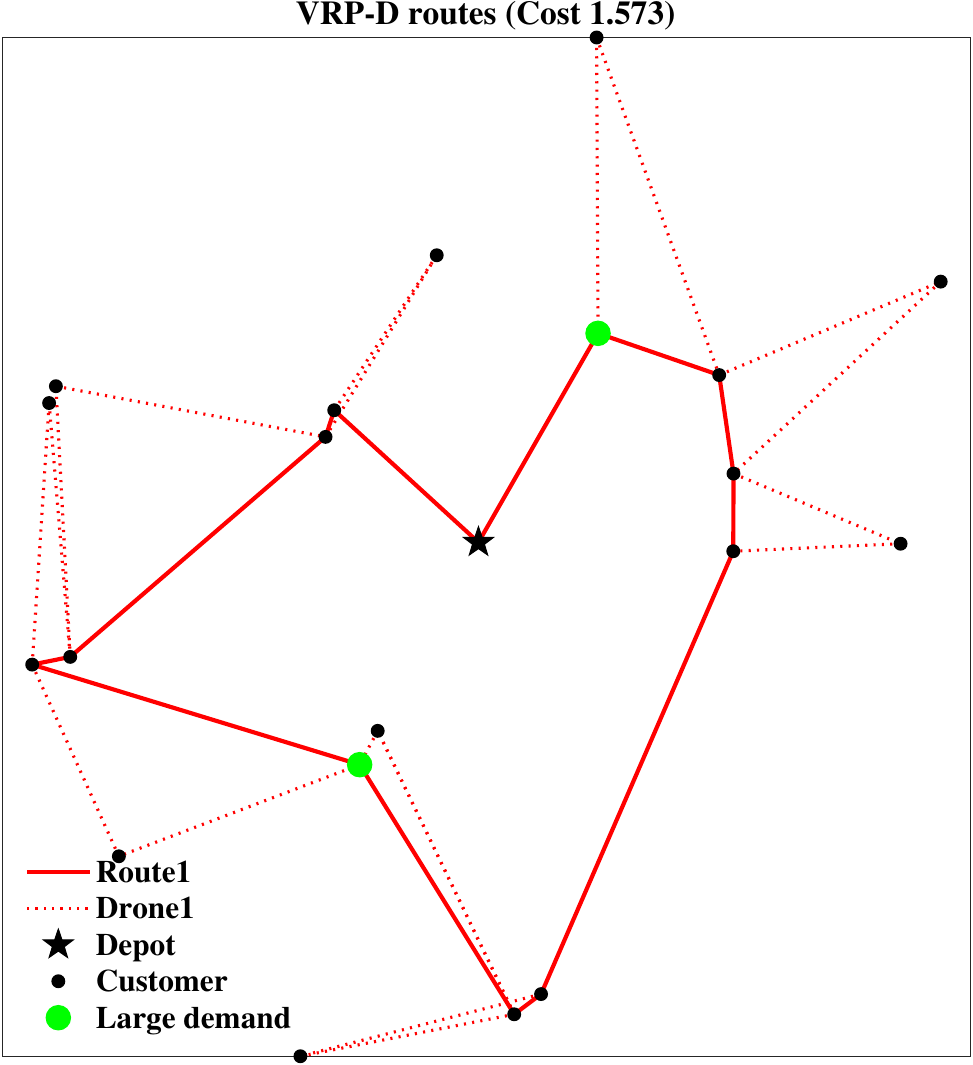} &        
    \includegraphics[width=0.33\textwidth]{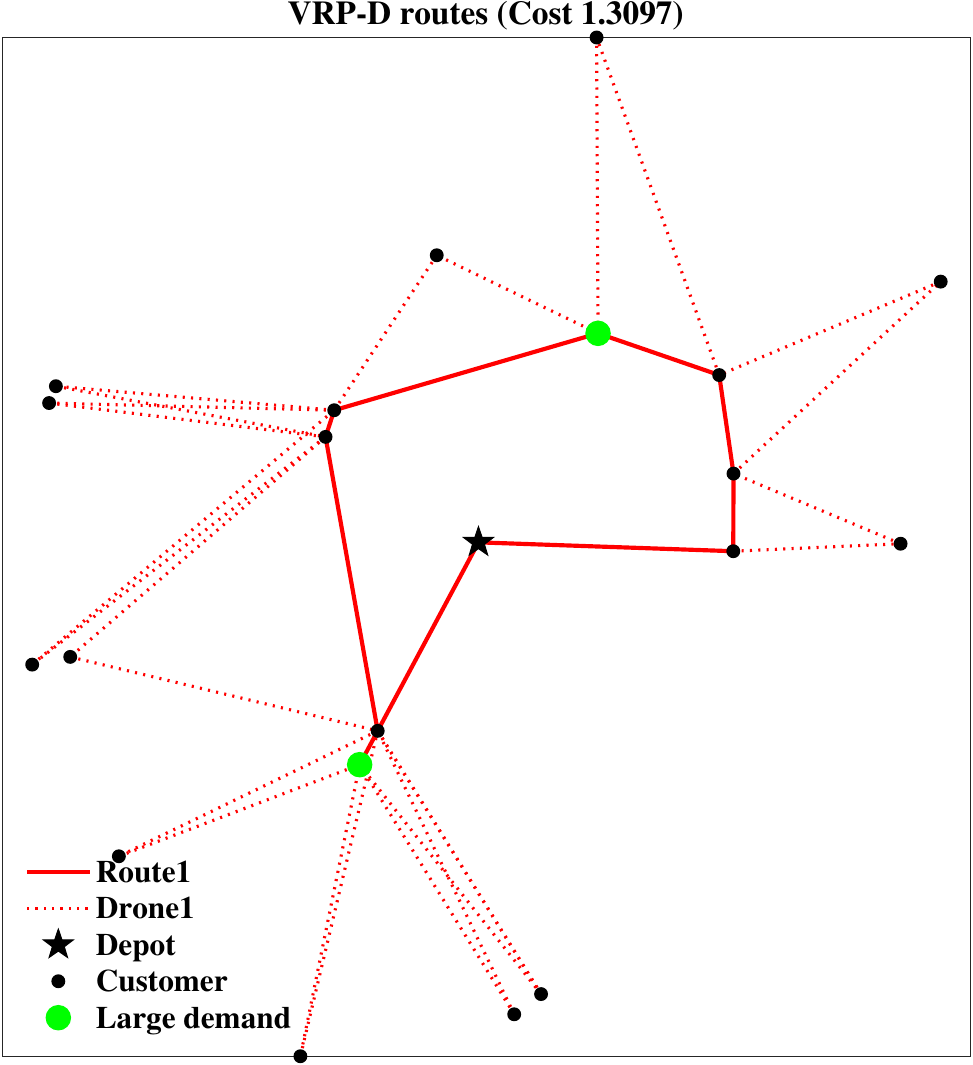} \\
        \multicolumn{1}{c}{(a) ALNS.} & \multicolumn{1}{c}{(b) LKH Config1.} & \multicolumn{1}{c}{(c) LKH Config2.}
    \end{tabular}
    \caption{VRPD tours on instance 20.5.4.}
    \label{fig:vrpd_figcompare2}
\end{figure}

Furthermore, LKH Config1 may generate inferior solutions than ALNS under the 5min runtime limit. For example, in instance 6.10.1 (Fig.~\ref{fig:vrpd_figcompare4}), LKH Config1 uses only one truck-drone group, resulting in a higher cost than ALNS. Similarly, for instance 20.5.4 (Fig.~\ref{fig:vrpd_figcompare2}), LKH Config1 produces a tour with cost of 1.573 compared to 1.37893 from ALNS. Nonetheless, LKH Config2 achieves the lowest cost for both instances by revisiting nodes.

\begin{figure}[H]
    \centering
    \setlength{\tabcolsep}{1pt} % Space between columns
    \begin{tabular}{ccc}        \includegraphics[width=0.33\textwidth]{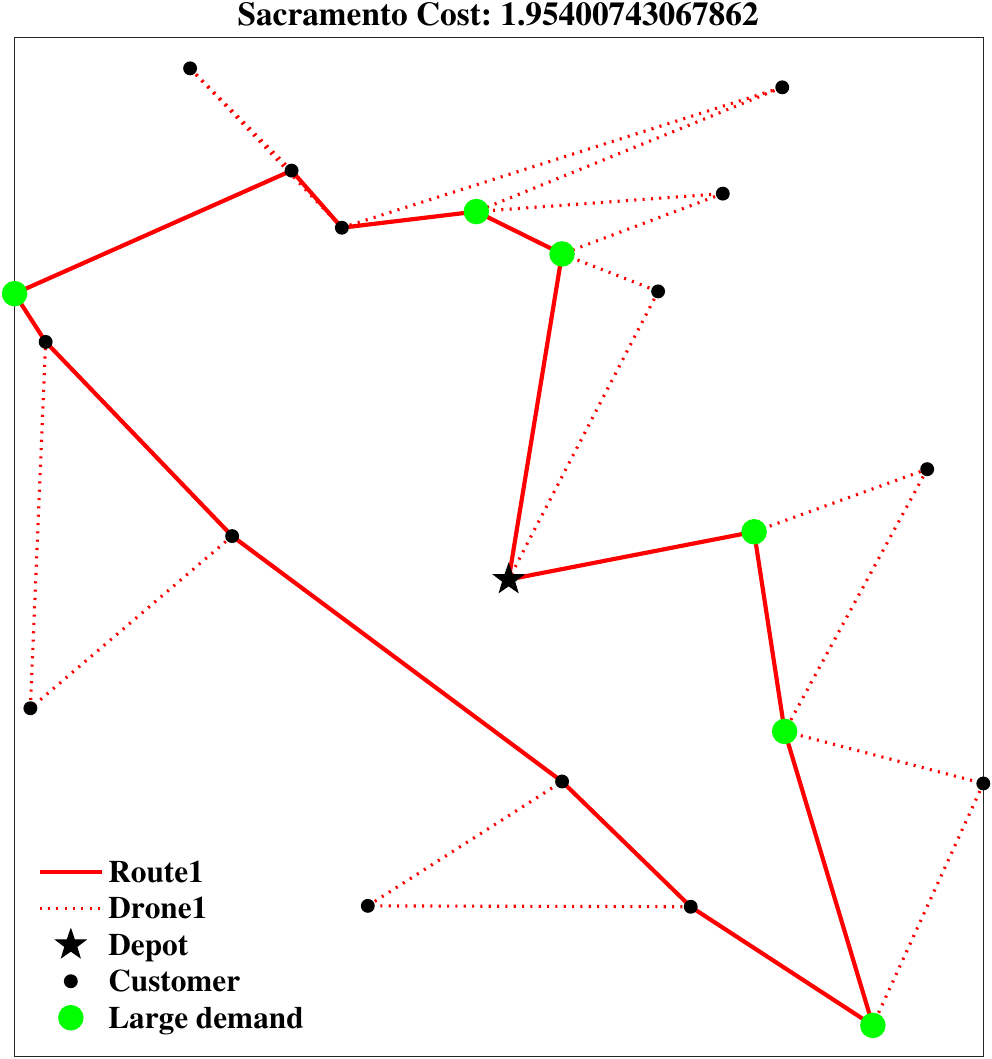} &        \includegraphics[width=0.33\textwidth]{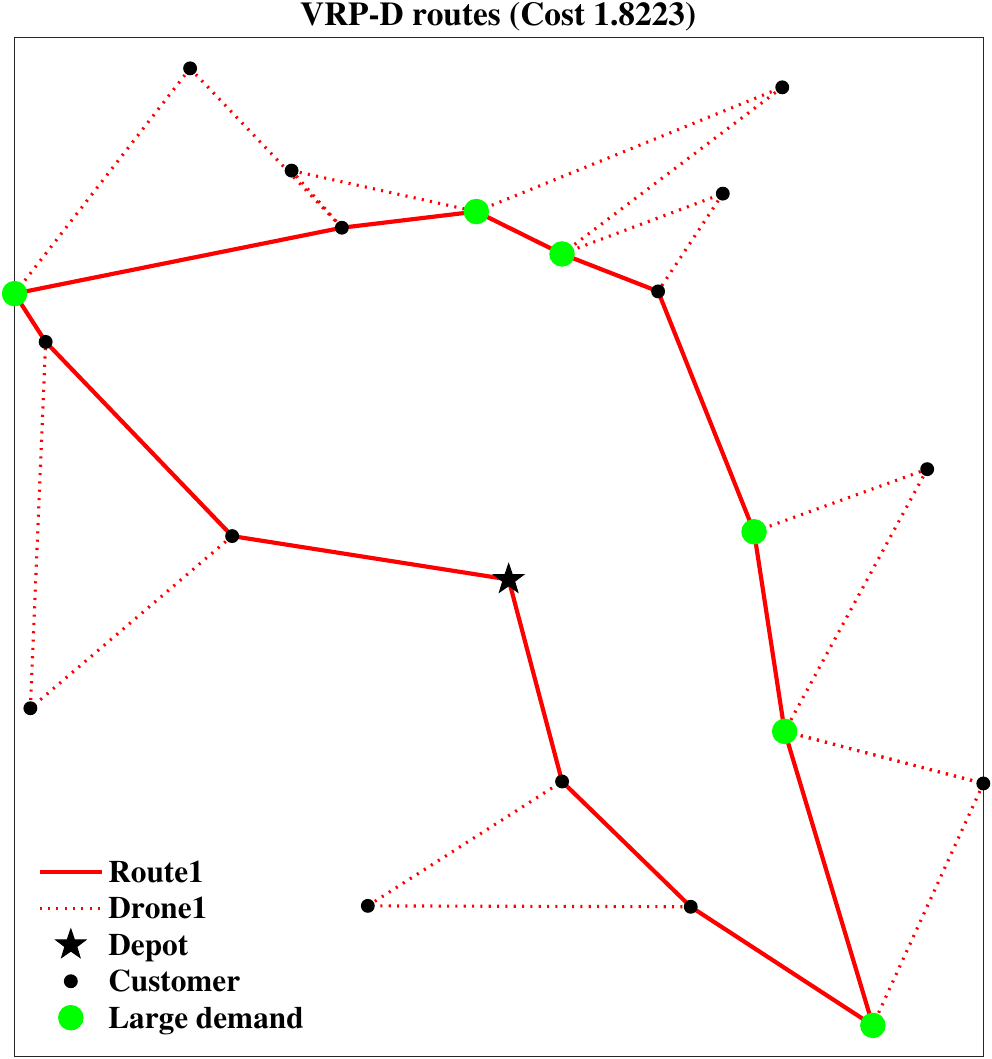} &        \includegraphics[width=0.33\textwidth]{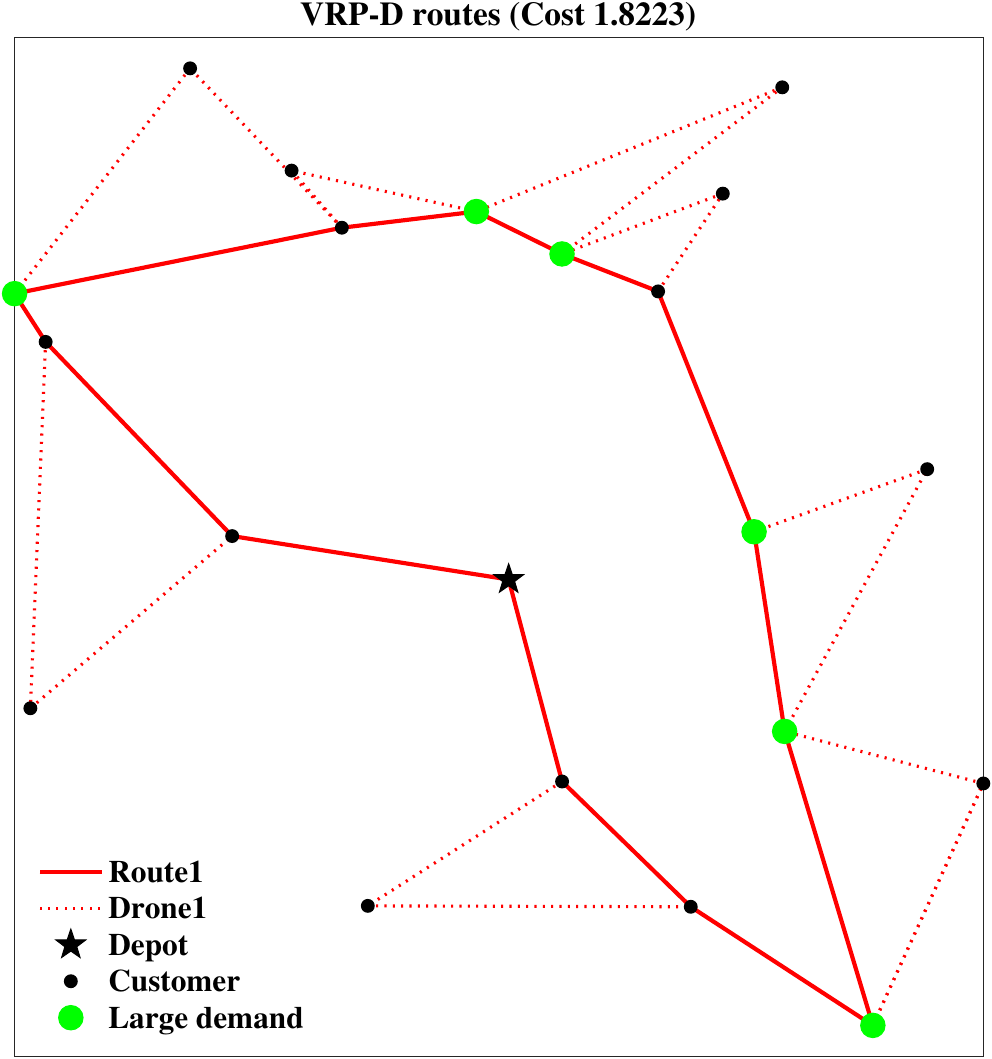} \\
        \multicolumn{1}{c}{(a) ALNS.} & \multicolumn{1}{c}{(b) LKH Config1.} & \multicolumn{1}{c}{(c) LKH Config2.}
    \end{tabular}
    \caption{VRPD tours on instance 20.5.2.}
    \label{fig:vrpd_figcompare3}
\end{figure}

Secondly, Fig.~\ref{fig:vrpd_figcompare3} visualizes instance 20.5.2, where both LKH Config1 and Config2 outperform the solution obtained from ALNS, though Config2 does not revisit nodes.

\subsubsection{Computational performance analysis}
% \subsection{Benchmark Results and Statistical Comparison}
% Thus, with increasing scale of instances, it becomes harder for LKH to catch up with ALNS.

Table~\ref{tab:vrpd_1LKH} compares LKH Config1 to the ALNS proposed by \cite{sacramento2019adaptive}, which provides both the exact (\textbf{CPLEX}) and heuristic (\textbf{ALNS}) results. For example, instance 6.5.1 represents the first location configuration among the four 6.5. instances, each with $6$ customers and one depot distributed in a $5\times5$ mile area. $N_D$ deontes the number of customers whose demands can be carried by a drone. 

According to \cite{sacramento2019adaptive}, the ALNS algorithm is executed 10 times independently with each run limited to 5min, and the best objective found is reported as $C^{\textit{ALNS}}$. Table~\ref{tab:vrpd_1LKH} also records the average number of drone-serviced customers (Dc), the average number of routes (V), and the computational time in seconds. For CPLEX, the optimal value ($z^*$) and the solution time ($t^{\textit{MIP}}$) by solving the MIP model formulated in \cite{sacramento2019adaptive} are included. Notably, $C^{\textit{ALNS}}$ and $z^*$ are identical for the first 36 instances.

For LKH Config1, run once in 5min, Table~\ref{tab:vrpd_1LKH} lists the objective value ($\text{C}^1$), the iteration count in a run (Iter) at which the reported solution is found, as well as Dc and V with no need for average. In column $\%gap^{\textit{ALNS}} = 100 \cdot (C^{\textit{ALNS}}-\text{C}^1)/C^{\textit{ALNS}}$, a positive value indicates that LKH Config1 produces a tour with a lower cost compared to ALNS. To further validate results, we re-execute the ALNS algorithm and compare its tours to those from LKH Config1 by illustrative figures. The 24 instances with visual difference are marked with $\ast$ in Table~\ref{tab:vrpd_1LKH}. Out of the 24 small instances of 6 or 10 customers, 18 exhibit zero gaps ($\%gap^{\textit{ALNS}} = 0$) without $\ast$, indicating that LKH Config1 achieves optimality in $75\%$ of these cases. Among the 24 medium instances with 12 or 20 customers, 6 are solved to optimality and 1 instance (20.5.2) is improved by LKH. However, LKH perform worse than ALNS in the remaining 17 medium instances, as evidenced by negative gaps. One possible explanation is that the ALNS is specifically designed for this VRP-D cost-minimization problem, while our generalized penalty-based solution framework is currently less specialized. 

% the average number of iterations (iter),  the Makespan of the delivery network, the number of drone-serviced customers (Dc), the number of “one truck-one drone” combinations (V), and the computational time (in seconds). Those instances with positive gaps are emphasized in bold. 

% For example, in Fig.~\ref{fig:vrpd_figcompare1}, the same tour is found for instance 6.10.3, yielding $\%gap^{\textit{ALNS}}=0.05$ in Table~\ref{tab:vrpd_1LKH}. This positive gap without improvement might be because LKH works with integer edge weights, which requires scaling and rounding data when generating the input problem file \texttt{.tsppd} in Matlab. Therefore, 

To address this limitation, we extend the LKH Config1's runtime from 5min to 24h for all 24 $\ast$ instances. Except for 10.20.4 and 20.5.2, Table~\ref{tab:vrpd_1LKH_long} shows that the other 22 $\ast$ instances are improved, as indicated by positive $\%gap^{\text{5min}} = 100 \cdot (\text{C}^1_{\text{5min}}-\text{C}^1_{\text{24h}})/\text{C}^1_{\text{5min}}$. Meanwhile, the increased number of drone-serviced customers (Dc) suggests a potential negative correlation with the min-cost objective value, whereas no such relationship is observed for the min-makespan objective in FSTSP and TSPmD.

We compare the figures obtained under 24h runtime with ALNS, finding that although 4 instances (20.10.1 - 20.10.4) do not match ALNS, other instances in Table~\ref{tab:vrpd_1LKH_long} at least achieve equivalent solutions, which is also evidenced by nonnegative $\%gap^{\textit{ALNS}} = 100 \cdot (C^{\textit{ALNS}}-\text{C}^1_{\text{24h}})/C^{\textit{ALNS}}$. In particular, the BKS results from \cite{sacramento2019adaptive} are improved for instances 20.5.2, 20.20.1 and 20.20.2, highlighted in bold. And significantly, LKH only uses 54s to achieve this improvement in instance 20.5.2. This breakthrough is visualized in Fig.~\ref{fig:vrpd_figcompare2}.

% first comparison
\begin{table}[H]
    \centering
    % \tiny
    % \scriptsize
    \footnotesize  % Slightly smaller to better fit
    % \small
    % \normalsize
    \caption{Comparison of ALNS and LKH (Config 1) on small instances}
    \label{tab:vrpd_1LKH}
    \begin{tabular*}{\textwidth}{@{\extracolsep{\fill}}l@{\;}r@{\;}r@{\;}r@{\;}r@{\;}r@{\;}r@{\;}r@{\;}r@{\;}r@{\;}r@{\;}r@{\;}r@{\;}r@{}}
    %%
    % \begin{tabularx}{\textwidth}{l r r r r r r r r r r r r r}
    % \begin{tabularx}{\textwidth}{l c c c c c c c c c c c c c}
    % \begin{tabularx}{\textwidth}{l *{13}{>{\centering\arraybackslash}X}}
        \toprule
        \multicolumn{2}{c}{} & \multicolumn{2}{c}{\textbf{CPLEX}} & \multicolumn{4}{c}{\textbf{ALNS}} & \multicolumn{5}{c}{\textbf{LKH (Config 1)}} \\
        \cmidrule(lr){3-4} \cmidrule(lr){5-8} \cmidrule(lr){9-13} 
        \raisebox{1.5ex}[0pt][0pt] {Instance} & \raisebox{1.5ex}[0pt][0pt] {$N_D$} & $z^*$ & $t^{\textit{MIP}}$(s) & $C^{\textit{ALNS}}$ & Dc & V & Time(s) & $\text{C}^1$ & Dc & V & Iter & Time(s) & \raisebox{1.5ex}[0pt][0pt] {$\%gap^{\textit{ALNS}}$} \\
        \midrule
6.5.1 & 5 & 1.09821 & 0.93 & 1.09821 & 3 & 1 & 0.014 & 1.09819 & 3 & 1 & 237 & 1.4 & 0 \\
6.5.2 & 6 & 0.84215 & 3.51 & 0.84215 & 3 & 1 & 0.001 & 0.84264 & 3 & 1 & 185 & 1.03 & 0 \\
6.5.3 & 5 & 1.21137 & 2.43 & 1.21137 & 3 & 1 & 0.001 & 1.21111 & 3 & 1 & 129 & 0.68 & 0 \\
6.5.4 & 5 & 0.94599 & 1.37 & 0.94599 & 3 & 1 & 0.003 & 0.94707 & 3 & 1 & 139 & 0.74 & 0 \\
6.10.1$\ast$ & 5 & 2.40611 & 7479 & 2.40611 & 4 & 2 & 0.004 & 2.52919 & 3 & 1 & 326 & 1.7 & \texttt{-5.12} \\
6.10.2 & 6 & 1.67927 & 6.8 & 1.67927 & 4 & 2 & 0.002 & 1.68018 & 4 & 2 & 34083 & 184.56 & 0 \\
6.10.3 & 6 & 1.32552 & 5.36 & 1.32552 & 4 & 2 & 0.003 & 1.32487 & 4 & 2 & 15818 & 104.7 & 0 \\
6.10.4 & 6 & 1.44307 & 5.23 & 1.44307 & 3 & 1 & 0.01 & 1.44363 & 3 & 1 & 105 & 0.48 & 0 \\
6.20.1 & 6 & 2.67759 & 4.73 & 2.67759 & 4 & 2 & 0.011 & 2.67777 & 4 & 2 & 2685 & 15.6 & 0 \\
6.20.2 & 5 & 4.31959 & 1.164 & 4.31959 & 3 & 1 & 0.054 & 4.32017 & 3 & 1 & 828 & 4.91 & 0 \\
6.20.3 & 6 & 3.82475 & 1.81 & 3.82475 & 4 & 2 & 0.002 & 3.82531 & 4 & 2 & 9294 & 62.81 & 0 \\
6.20.4 & 6 & 3.67872 & 2.17 & 3.67872 & 4 & 2 & 0.001 & 3.67917 & 4 & 2 & 527 & 3.01 & 0 \\
10.5.1 & 5 & 1.65563 & 13.37 & 1.65563 & 2 & 1 & 0.002 & 1.65493 & 2 & 1 & 175 & 1.56 & 0 \\
10.5.2 & 9 & 1.45185 & 433.15 & 1.45185 & 5 & 1 & 0.339 & 1.45223 & 5 & 1 & 2078 & 22.75 & 0 \\
10.5.3$\ast$ & 8 & 1.47357 & 236.11 & 1.47357 & 5 & 1 & 0.193 & 1.49945 & 5 & 1 & 14252 & 222.08 & \texttt{-1.76} \\
10.5.4 & 9 & 1.28489 & 345.22 & 1.28489 & 5 & 1 & 0.002 & 1.28497 & 5 & 1 & 859 & 10.47 & 0 \\
10.10.1 & 8 & 2.32647 & 369.92 & 2.32647 & 5 & 1 & 0.026 & 2.32692 & 5 & 1 & 2642 & 35.52 & 0 \\
10.10.2$\ast$ & 8 & 3.15856 & 121.28 & 3.15856 & 5 & 1 & 0.075 & 3.19025 & 5 & 1 & 9058 & 118.12 & \texttt{-1}\\
10.10.3$\ast$ & 7 & 2.55274 & 88.41 & 2.55274 & 6 & 2 & 0.427 & 2.62386 & 5 & 1 & 19011 & 247.81 & \texttt{-2.79} \\
10.10.4$\ast$ & 9 & 2.53931 & 246.43 & 2.53931 & 5 & 1 & 0.008 & 2.5414 & 5 & 1 & 17511 & 287.3 & \texttt{-0.08} \\
10.20.1 & 7 & 4.4524 & 6.65 & 4.4524 & 4 & 1 & 3.946 & 4.45219 & 4 & 1 & 5496 & 62.63 & 0 \\
10.20.2 & 8 & 6.16776 & 180.16 & 6.16776 & 4 & 1 & 0.011 & 6.1675 & 4 & 1 & 4489 & 62.06 & 0 \\
10.20.3 & 9 & 4.5463 & 251.14 & 4.5463 & 5 & 1 & 1.197 & 4.54611 & 5 & 1 & 2084 & 33.76 & 0 \\
10.20.4$\ast$ & 7 & 6.15355 & 275.36 & 6.15355 & 4 & 2 & 49.17 & 6.17069 & 4 & 1 & 1678 & 21.47 & \texttt{-0.28} \\
12.5.1 & 9 & 1.37381 & 1161.88 & 1.37381 & 6 & 1 & 31.444 & 1.37465 & 6 & 1 & 3754 & 70.38 & 0 \\
12.5.2$\ast$ & 12 & 1.05899 & 62131.17 & 1.05899 & 7 & 2 & 1.11 & 1.12875 & 7 & 2 & 10998 & 248.73 & \texttt{-6.59} \\
12.5.3 & 10 & 1.44765 & 433.9 & 1.44765 & 6 & 1 & 0.028 & 1.4466 & 6 & 1 & 3457 & 58.1 & 0 \\
12.5.4$\ast$ & 10 & 1.581 & 2259.66 & 1.581 & 6 & 1 & 0.1 & 1.61216 & 6 & 1 & 10224 & 161.28 & \texttt{-1.97} \\
12.10.1$\ast$ & 10 & 2.68103 & 811.26 & 2.68103 & 7 & 2 & 81.447 & 2.79981 & 6 & 1 & 1651 & 31.55 & \texttt{-4.43} \\
12.10.2$\ast$ & 10 & 2.6842 & 1004.35 & 2.6842 & 6 & 1 & 0.059 & 2.68753 & 6 & 1 & 3711 & 71.99 & \texttt{-0.12} \\
12.10.3$\ast$ & 9 & 2.88048 & 793.87 & 2.88048 & 6 & 1 & 0.03 & 2.93852 & 4 & 1 & 7404 & 149.62 & \texttt{-2.01} \\
12.10.4 & 10 & 2.31418 & 176.74 & 2.31418 & 6 & 1 & 0.011 & 2.31397 & 6 & 1 & 5458 & 92.62 & 0 \\
12.20.1 & 11 & 5.77759 & 3723.83 & 5.77759 & 7 & 2 & 0.272 & 5.77717 & 7 & 2 & 9856 & 200.92 & 0 \\
12.20.2 & 10 & 8.27254 & 1081.57 & 8.27254 & 4 & 1 & 0.004 & 8.27336 & 4 & 1 & 14550 & 209.82 & 0 \\
12.20.3 & 9 & 4.16693 & 24.52 & 4.16693 & 5 & 1 & 0.054 & 4.16692 & 5 & 1 & 2755 & 38.32 & 0 \\
12.20.4$\ast$ & 11 & 6.08859 & 1335.74 & 6.08859 & 7 & 2 & 0.21 & 6.40926 & 6 & 1 & 13456 & 188.9 & \texttt{-5.27} \\
20.5.1$\ast$ & 15 & & & 1.79347 & 9 & 1 & & 1.81974 & 8 & 1 & 2041 & 79.72 & \texttt{-1.46} \\
\textbf{20.5.2}$\ast$ & 14 & & & \textbf{1.95401} & 8 & 1 & & \textbf{1.82229} & 8 & 1 & 1828 & 54.18 & \textbf{6.74} \\
20.5.3$\ast$ & 19 & & & 1.48658 & 9 & 1 & & 1.8087 & 8 & 1 & 8370 & 300.02 & \texttt{-21.67} \\
20.5.4$\ast$ & 18 & & & 1.37893 & 10 & 1 & & 1.573 & 9 & 1 & 5362 & 171.34 & \texttt{-14.07} \\
20.10.1$\ast$ & 17 & & & 3.25253 & 10 & 1 & & 3.42054 & 8 & 1 & 5828 & 178.69 & \texttt{-5.17} \\
20.10.2$\ast$ & 19 & & & 3.08938 & 10 & 1 & & 3.47658 & 9 & 1 & 4580 & 216.47 & \texttt{-12.53} \\
20.10.3$\ast$ & 19 & & & 3.70226 & 9 & 1 & & 3.82191 & 8 & 1 & 6861 & 228.47 & \texttt{-3.23} \\
20.10.4$\ast$ & 15 & & & 3.3089 & 10 & 1 & & 3.63194 & 9 & 1 & 2289 & 75.03 & \texttt{-9.76} \\
20.20.1$\ast$ & 19 & & & 7.34453 & 10 & 1 & & 7.4771 & 9 & 1 & 7342 & 281.74 & \texttt{-1.81} \\
20.20.2$\ast$ & 16 & & & 7.54889 & 9 & 1 & & 7.66324 & 7 & 1 & 1728 & 63.32 & \texttt{-1.51} \\
20.20.3$\ast$ & 18 & & & 7.461 & 10 & 1 & & 7.94447 & 9 & 1 & 2583 & 117.93 & \texttt{-6.48} \\
20.20.4$\ast$ & 17 & & & 7.01331 & 9 & 1 & & 7.78857 & 8 & 1 & 5554 & 203.43 & \texttt{-11.05} \\
        \bottomrule
    \end{tabular*}
    \begin{tablenotes}\footnotesize
    \item 24$\ast$: Different visualization between ALNS and LKH.
    \end{tablenotes}
\end{table}

\begin{table}[H]
    \centering
    % \tiny
    % \scriptsize
    \footnotesize  % Slightly smaller to better fit
    % \small
    % \normalsize
    \caption{Comparison of ALNS and LKH (Config 1 24h)}
    \label{tab:vrpd_1LKH_long}
    \begin{tabular*}{\textwidth}{@{\extracolsep{\fill}}l@{\;}r@{\;}r@{\;}r@{\;}r@{\;}r@{\;}r@{\;}r@{\;}r@{\;}r@{\;}r@{\;}r@{\;}r@{\;}r@{\;}r@{\;}r@{\;}r@{}}
    %%
    % \begin{tabularx}{\textwidth}{l r r r r r r r r r r r r}
    % \begin{tabularx}{\textwidth}{l c c c c c c c c c c c c}
    % \begin{tabularx}{\textwidth}{l *{12}{>{\centering\arraybackslash}X}}
        \toprule
        \multicolumn{2}{c}{} & \multicolumn{3}{c}{\textbf{ALNS}} & \multicolumn{5}{c}{\textbf{LKH (Config 1 5min)}}  & \multicolumn{5}{c}{\textbf{LKH (Config 1 24h)}} \\
        \cmidrule(lr){3-5} \cmidrule(lr){6-10} \cmidrule(lr){11-15}
        \raisebox{1.5ex}[0pt][0pt] {Instance} & \raisebox{1.5ex}[0pt][0pt] {$N_D$} & $C^{\textit{ALNS}}$ & Dc & V & $\text{C}^1_{\text{5min}}$ & Dc & V & Iter & Time(s) & $\text{C}^1_{\text{24h}}$ & Dc & V & Iter & Time(s) & \raisebox{1.5ex}[0pt][0pt] {$\%gap^{\textit{ALNS}}$} & \raisebox{1.5ex}[0pt][0pt] {$\%gap^{\text{5min}}$} \\
        \midrule
6.10.1 & 5 & 2.40611 & 4 & 2 & 2.52919 & 3 & 1 & 326 & 1.7 & 2.4063 & 4 & 2 & 92629 & 351.22 & 0 & 4.86 \\
10.5.3 & 8 & 1.47357 & 5 & 1 & 1.49945 & 5 & 1 & 14252 & 222.08 & 1.47334 & 5 & 1 & 30989 & 507.22 & 0 & 1.74 \\
10.10.2 & 8 & 3.15856 & 5 & 1 & 3.19025 & 5 & 1 & 9058 & 118.12 & 3.15958 & 5 & 1 & 146290 & 20216.68 & 0 & 0.96 \\
10.10.3 & 7 & 2.55274 & 6 & 2 & 2.62386 & 5 & 1 & 19011 & 247.81 & 2.55349 & 6 & 2 & 719385 & 10482.62 & 0 & 2.68 \\
10.10.4 & 9 & 2.53931 & 5 & 1 & 2.5414 & 5 & 1 & 17511 & 287.3 & 2.53917 & 5 & 1 & 20170 & 326.35 & 0 & 0.09 \\
12.5.2 & 12 & 1.05899 & 7 & 2 & 1.12875 & 7 & 2 & 10998 & 248.73 & 1.05765 & 7 & 2 & 24420 & 813.07 & 0 & 6.3 \\
12.5.4 & 10 & 1.581 & 6 & 1 & 1.61216 & 6 & 1 & 10224 & 161.28 & 1.58106 & 6 & 1 & 451586 & 8501.58 & 0 & 1.93 \\
12.10.1 & 10 & 2.68103 & 7 & 2 & 2.79981 & 6 & 1 & 1651 & 31.55 & 2.68074 & 7 & 2 & 22371 & 524.13 & 0 & 4.25 \\
12.10.2 & 10 & 2.6842 & 6 & 1 & 2.68753 & 6 & 1 & 3711 & 71.99 & 2.68392 & 6 & 1 & 2515995 & 54229.45 & 0 & 0.13 \\
12.10.3 & 9 & 2.88048 & 6 & 1 & 2.93852 & 4 & 1 & 7404 & 149.62 & 2.88068 & 6 & 1 & 77985 & 1556.9 & 0 & 1.97 \\
12.20.4 & 11 & 6.08859 & 7 & 2 & 6.40926 & 6 & 1 & 13456 & 188.9 & 6.08876 & 7 & 2 & 568178 & 10189.36 & 0 & 5 \\
20.5.1 & 15 & 1.79347 & 9 & 1 & 1.81974 & 8 & 1 & 2041 & 79.72 & 1.7931 & 9 & 1 & 63187 & 3123.05 & 0 & 1.46 \\
\textbf{20.5.2}$\dagger$ & 14 & \textbf{1.95401} & 8 & 1 & \textbf{1.82229} & 8 & 1 & 1828 & 54.18 & \textbf{1.82229} & 8 & 1 & 1828 & 54.64 & \textbf{6.74} & 0 \\
20.5.3 & 19 & 1.48658 & 9 & 1 & 1.8087 & 8 & 1 & 8370 & 300.02 & 1.48661 & 9 & 1 & 107995 & 3963.18 & 0 & 17.8 \\
20.5.4 & 18 & 1.37893 & 10 & 1 & 1.573 & 9 & 1 & 5362 & 171.34 & 1.37868 & 10 & 1 & 23896 & 881.13 & 0 & 12.4 \\
20.10.1$\dagger$ & 17 & 3.25253 & 10 & 1 & 3.42054 & 8 & 1 & 5828 & 178.69 & 3.39656 & 8 & 1 & 363743 & 15728.05 & \texttt{-4.43} & 0.7 \\
20.10.2$\dagger$ & 19 & 3.08938 & 10 & 1 & 3.47658 & 9 & 1 & 4580 & 216.47 & 3.11575 & 10 & 1 & 219886 & 11673.55 & \texttt{-0.85} & 10.4 \\
20.10.3$\dagger$ & 19 & 3.70226 & 9 & 1 & 3.82191 & 8 & 1 & 6861 & 228.47 & 3.79814 & 9 & 1 & 11957 & 404.5 & \texttt{-2.59} & 0.62 \\
20.10.4$\dagger$ & 15 & 3.3089 & 10 & 1 & 3.63194 & 9 & 1 & 2289 & 75.03 & 3.37225 & 9 & 1 & 1393556 & 68910.21 & \texttt{-1.91} & 7.15 \\
\textbf{20.20.1}$\dagger$ & 19 & \textbf{7.34453} & 10 & 1 & \textbf{7.4771} & 9 & 1 & 7342 & 281.74 & \textbf{7.32289} & 10 & 2 & 125236 & 6285.33 & \textbf{0.29} & 2.06 \\
\textbf{20.20.2}$\dagger$ & 16 & \textbf{7.54889} & 9 & 1 & \textbf{7.66324} & 7 & 1 & 1728 & 63.32 & \textbf{7.53844} & 9 & 1 & 1102908 & 35954.28 & \textbf{0.14} & 1.63 \\
20.20.3 & 18 & 7.461 & 10 & 1 & 7.94447 & 9 & 1 & 2583 & 117.93 & 7.46033 & 10 & 1 & 26777 & 1386.68 & 0 & 6.09 \\
20.20.4 & 17 & 7.01331 & 9 & 1 & 7.78857 & 8 & 1 & 5554 & 203.43 & 7.01322 & 9 & 1 & 399602 & 19776.2 & 0 & 9.95 \\
        \bottomrule
    \end{tabular*}
    \begin{tablenotes}\footnotesize
    \item 7$\dagger$: Different visualization between ALNS and LKH (Config1 24h).
    \end{tablenotes}
\end{table}

Subsequently, we examine two configurations of LKH in Table~\ref{tab:vrpd_2LKH}. The bold instances indicate cases where Config2 revisits certain customer nodes. Note that 41 out of 48 $\%gap^{1,2} = 100 \cdot (\text{M}^1-\text{M}^2)/\text{M}^1$ values are nonnegative. This suggests that regardless of whether node revisit is implemented, LKH Config2 tends to derive superior solutions compared to Config1. This consistent improvement is not observed in Tables~\ref{tab:fstsp_2LKH_small}--~\ref{tab:fstsp_2LKH_medium} for FSTSP or Tables~\ref{tab:tspmd_2LKH_small}--\ref{tab:tspmd_2LKH_medium} for TSP-mD, possibly because the cost-minimization objective of VRP-D is more sensitive to the flexibility of revisiting nodes than makespan-minimization.

\begin{table}[H]
    \centering
    % \tiny
    % \scriptsize
    \footnotesize  % Slightly smaller to better fit
    % \small
    % \normalsize
    \caption{Comparison of LKH (Config 1) and LKH (Config 2) on small instances}
    \label{tab:vrpd_2LKH}
    \begin{tabular*}{\textwidth}{@{\extracolsep{\fill}}l@{\;}r@{\;}r@{\;}r@{\;}r@{\;}r@{\;}r@{\;}r@{\;}r@{\;}r@{\;}r@{\;}r@{\;}r@{}}
    %%
    % \begin{tabularx}{\textwidth}{l r r r r r r r r r r r r}
    % \begin{tabularx}{\textwidth}{l c c c c c c c c c c c c}
    % \begin{tabularx}{\textwidth}{l *{12}{>{\centering\arraybackslash}X}}
        \toprule
        \multicolumn{2}{c}{} & \multicolumn{5}{c}{\textbf{LKH (Config 1)}} & \multicolumn{5}{c}{\textbf{LKH (Config 2)}} \\
        \cmidrule(lr){3-7} \cmidrule(lr){8-12}
        \raisebox{1.5ex}[0pt][0pt] {Instance} & \raisebox{1.5ex}[0pt][0pt] {$N_D$} & $\text{C}^1$ & Dc & V & Iter & Time(s) & $\text{C}^2$ & Dc & V & Iter & Time(s) & \raisebox{1.5ex}[0pt][0pt] {$\%gap^{1,2}$} \\
        \midrule
6.5.1 & 5 & 1.09819 & 3 & 1 & 237 & 1.4 & 1.09819 & 3 & 1 & 115 & 0.64 & 0 \\
\textbf{6.5.2} & 6 & \textbf{0.84264} & 3 & 1 & 185 & 1.03 & \textbf{0.84253 revisit} & 4 & 1 & 941 & 4.1 & \textbf{0.01} \\
6.5.3 & 5 & 1.21111 & 3 & 1 & 129 & 0.68 & 1.21111 & 3 & 1 & 222 & 0.98 & 0 \\
\textbf{6.5.4} & 5 & \textbf{0.94707} & 3 & 1 & 139 & 0.74 & \textbf{0.91534 revisit} & 4 & 1 & 734 & 4.31 & \textbf{3.35} \\
\textbf{6.10.1} & 5 & \textbf{2.52919} & 3 & 1 & 326 & 1.7 & \textbf{1.97394 revisit} & 4 & 1 & 20 & 0.1 & \textbf{21.95} \\
\textbf{6.10.2} & 6 & \textbf{1.68018} & 4 & 2 & 34083 & 184.56 & \textbf{1.31437 revisit} & 4 & 1 & 43 & 0.25 & \textbf{21.77} \\
\textbf{6.10.3} & 6 & \textbf{1.32487} & 4 & 2 & 15818 & 104.7 & \textbf{1.09299 revisit} & 4 & 1 & 436 & 2.57 & \textbf{17.5} \\
\textbf{6.10.4} & 6 & \textbf{1.44363} & 3 & 1 & 105 & 0.48 & \textbf{1.30429 revisit} & 4 & 1 & 140 & 0.59 & \textbf{9.65} \\
\textbf{6.20.1} & 6 & \textbf{2.67777} & 4 & 2 & 2685 & 15.6 & \textbf{2.28309 revisit} & 4 & 1 & 86 & 0.46 & \textbf{14.74} \\
6.20.2 & 5 & 4.32017 & 3 & 1 & 828 & 4.91 & 4.32017 & 3 & 1 & 411 & 2.06 & 0 \\
6.20.3 & 6 & 3.82531 & 4 & 2 & 9294 & 62.81 & 3.82531 & 4 & 2 & 1417 & 9.24 & 0 \\
6.20.4 & 6 & 3.67917 & 4 & 2 & 527 & 3.01 & 3.67917 & 4 & 2 & 1165 & 8.07 & 0 \\
10.5.1 & 5 & 1.65493 & 2 & 1 & 175 & 1.56 & 1.65493 & 2 & 1 & 85 & 0.79 & 0 \\
10.5.2 & 9 & 1.45223 & 5 & 1 & 2078 & 22.75 & 1.45223 & 5 & 1 & 10179 & 155.59 & 0 \\
\textbf{10.5.3} & 8 & \textbf{1.49945} & 5 & 1 & 14252 & 222.08 & \textbf{1.41848 revisit} & 6 & 1 & 10662 & 204.6 & \textbf{5.4} \\
10.5.4 & 9 & 1.28497 & 5 & 1 & 859 & 10.47 & 1.28497 & 5 & 1 & 1177 & 16.67 & 0 \\
\textbf{10.10.1} & 8 & \textbf{2.32692} & 5 & 1 & 2642 & 35.52 & \textbf{2.11488 revisit} & 6 & 1 & 7921 & 103.35 & \textbf{9.11} \\
10.10.2 & 8 & 3.19025 & 5 & 1 & 9058 & 118.12 & 3.19025 & 5 & 1 & 2583 & 37.53 & 0 \\
\textbf{10.10.3} & 7 & \textbf{2.62386} & 5 & 1 & 19011 & 247.81 & \textbf{2.35047 revisit} & 6 & 1 & 22217 & 285.74 & \textbf{10.42} \\
10.10.4 & 9 & 2.5414 & 5 & 1 & 17511 & 287.3 & 2.53917 & 5 & 1 & 5675 & 91.68 & 0.09 \\
10.20.1 & 7 & 4.45219 & 4 & 1 & 5496 & 62.63 & 4.46195 & 4 & 1 & 4231 & 36.67 & -0.22 \\
10.20.2 & 8 & 6.1675 & 4 & 1 & 4489 & 62.06 & 6.1675 & 4 & 1 & 299 & 2.99 & 0 \\
10.20.3 & 9 & 4.54611 & 5 & 1 & 2084 & 33.76 & 4.5926 & 5 & 1 & 1073 & 17.46 & -1.02 \\
\textbf{10.20.4} & 7 & \textbf{6.17069} & 4 & 1 & 1678 & 21.47 & \textbf{5.09245 revisit} & 6 & 2 & 909 & 14.12 & \textbf{17.47} \\
12.5.1 & 9 & 1.37465 & 6 & 1 & 3754 & 70.38 & 1.37465 & 6 & 1 & 106 & 1.34 & 0 \\
\textbf{12.5.2} & 12 & \textbf{1.12875} & 7 & 2 & 10998 & 248.73 & \textbf{0.87851 revisit} & 10 & 4 & 6218 & 261.91 & \textbf{22.17} \\
12.5.3 & 10 & 1.4466 & 6 & 1 & 3457 & 58.1 & 1.4466 & 6 & 1 & 12430 & 174.61 & 0 \\
12.5.4 & 10 & 1.61216 & 6 & 1 & 10224 & 161.28 & 1.59804 & 6 & 1 & 391 & 6.3 & 0.88 \\
12.10.1 & 10 & 2.79981 & 6 & 1 & 1651 & 31.55 & 2.76468 & 6 & 1 & 6005 & 128.19 & 1.25 \\
12.10.2 & 10 & 2.68753 & 6 & 1 & 3711 & 71.99 & 2.69549 & 6 & 1 & 6443 & 146.06 & -0.3 \\
12.10.3 & 9 & 2.93852 & 4 & 1 & 7404 & 149.62 & 2.88068 & 6 & 1 & 6273 & 149.84 & 1.97 \\
12.10.4 & 10 & 2.31397 & 6 & 1 & 5458 & 92.62 & 2.31397 & 6 & 1 & 2170 & 44.89 & 0 \\
\textbf{12.20.1} & 11 & \textbf{5.77717} & 7 & 2 & 9856 & 200.92 & \textbf{5.79988 revisit} & 7 & 1 & 14592 & 268.59 & \textbf{-0.39} \\
12.20.2 & 10 & 8.27336 & 4 & 1 & 14550 & 209.82 & 8.27336 & 4 & 1 & 4933 & 81.7 & 0 \\
12.20.3 & 9 & 4.16692 & 5 & 1 & 2755 & 38.32 & 4.16692 & 5 & 1 & 3667 & 56.65 & 0 \\
\textbf{12.20.4} & 11 & \textbf{6.40926} & 6 & 1 & 13456 & 188.9 & \textbf{5.97297 revisit} & 7 & 2 & 8819 & 158.97 & \textbf{6.81} \\
20.5.1 & 15 & 1.81974 & 8 & 1 & 2041 & 79.72 & 1.79438 & 9 & 1 & 5323 & 279.33 & 1.39 \\
20.5.2 & 14 & 1.82229 & 8 & 1 & 1828 & 54.18 & 1.82229 & 8 & 1 & 3795 & 109.74 & 0 \\
20.5.3 & 19 & 1.8087 & 8 & 1 & 8370 & 300.02 & 1.48725 & 9 & 1 & 2977 & 120.01 & 17.77 \\
\textbf{20.5.4} & 18 & \textbf{1.573} & 9 & 1 & 5362 & 171.34 & \textbf{1.3097 revisit} & 12 & 1 & 5841 & 248.26 & \textbf{16.74} \\
20.10.1 & 17 & 3.42054 & 8 & 1 & 5828 & 178.69 & 3.45015 & 7 & 1 & 4940 & 198.35 & -0.87 \\
20.10.2 & 19 & 3.47658 & 9 & 1 & 4580 & 216.47 & 3.36939 & 10 & 1 & 6845 & 298.71 & 3.08 \\
20.10.3 & 19 & 3.82191 & 8 & 1 & 6861 & 228.47 & 3.83581 & 9 & 1 & 6200 & 292.42 & -0.36 \\
20.10.4 & 15 & 3.63194 & 9 & 1 & 2289 & 75.03 & 3.44654 & 9 & 1 & 5412 & 220.98 & 5.1 \\
\textbf{20.20.1} & 19 & \textbf{7.4771} & 9 & 1 & 7342 & 281.74 & \textbf{7.23237 revisit} & 10 & 1 & 3117 & 150.19 & \textbf{3.27} \\
20.20.2 & 16 & 7.66324 & 7 & 1 & 1728 & 63.32 & 7.58566 & 8 & 1 & 711 & 29.19 & 1.01 \\
20.20.3 & 18 & 7.94447 & 9 & 1 & 2583 & 117.93 & 8.1077 & 10 & 1 & 3142 & 149.46 & -2.05 \\
20.20.4 & 17 & 7.78857 & 8 & 1 & 5554 & 203.43 & 7.02277 & 9 & 1 & 4243 & 224.23 & 9.83 \\
        \bottomrule
    \end{tabular*}
\end{table}

In short, under 5min and extended runtimes in Tables~\ref{tab:vrpd_1LKH}--\ref{tab:vrpd_1LKH_long}, LKH Config1 can generate the optimal solutions for all instances up to 12 customers except one instance 10.20.4, and even improve the BKS results for three 20-customer instances as bold in Table~\ref{tab:vrpd_1LKH_long}.

\section{Discussion} \label{sec:discussion} 
% Leveraging this penalty-based solution framework opens a new direction for solving truck-and-drone routing problems.

This paper presents the first work to propose a unified, flexible, and malleable solution framework based on LKH-3 to systematically address a wide range of truck-and-drone routing problems, where each problem variant enriches the structural skeleton by incorporating diverse operational constraints. This framework includes three phases: setup in Matlab, execution in LKH-3, and closure in Matlab. Note that LKH-3 is a state-of-the-art heuristic solver for constrained VRPs that handles constraints through penalty functions. Regarding the penalty function embedded in our solution framework, it can not only help generate the best valid truck-and-drone tour under its specific coordination structure, but also address realistic features, such as drone flight endurance, drone operation time, customer service time, subset of drone-eligible customers and customer time windows, by their corresponding penalty terms. Notably, we can adapt the first phase of generating input files including ``\texttt{.tour}, \texttt{.drone}, \texttt{.par}" for LKH-3 and the second phase of designing the penalty function ``Penalty\_drone.c", to potentially handle all truck-and-drone variants. 

The effectiveness of this penalty-based framework is validated on three problem variants, that is, FSTSP, TSP-mD, and 1-1 VRP-D where each group comprises one truck and one dedicated drone. For each variant, we design a specific penalty function and compare its performance with a customized algorithm in the literature. According to the computational results in Section 5, we have several findings:

\begin{enumerate}
\item Our adjustable framework achieves optimal solutions for nearly all small Poikonen Set4 and Sacramento Set6 with limited 5min runtime (Tables~\ref{tab:fstsp_1LKH_small}, \ref{tab:tspmd_1LKH_small}, \ref{tab:vrpd_1LKH}) and outperforms the compared algorithms on medium instances with extended runtime (Tables~\ref{tab:fstsp_1LKH_longer}, \ref{tab:tspmd_1LKH_longer}, \ref{tab:vrpd_1LKH_long}). For Murray Set1, it can match the optimal solutions for 69 out of 72 10-customer instances (Tables~\ref{tab:Boccia10cusE20_1LKH}, \ref{tab:Boccia10cusE40_1LKH}), and even improve the BKS upper bound of two 20-customer instances (Tables~\ref{tab:Boccia20cusE20_1LKH}, \ref{tab:Boccia20cusE40_1LKH}), implementing exactly the same FSTSP settings within 1h time limit. Regrading Sacramento Set6, it improves the BKSes of three 20-customer instances for 1-1 VRP-D (Table~\ref{tab:vrpd_1LKH_long}) which involves realistic constraints.

\item For FSTSP and TSP-mD variants, the compared algorithms - HGA and ALNS - cannot generate valid tours for the small instances poi-7-23 and poi-10-23 across five independent runs (Tables \ref{tab:fstsp_1LKH_small}, \ref{tab:tspmd_1LKH_small}). Conversely, LKH consistently produces valid solutions.

\item Investigating the computational times averaged in Tables~\ref{tab:Boccia10cusE20_1LKH}, ~\ref{tab:Boccia10cusE40_1LKH}, ~\ref{tab:Boccia20cusE20_1LKH}, and ~\ref{tab:Boccia20cusE40_1LKH} to solve FSTSP in Murray Set1, we find that although LKH costs more time in small instances, it takes less time (369.81s) to solve the most complex case (20-customer instances with endurance 40), indicating the superior scalability of LKH.

% \begin{table}[H]
%     \centering
%     % \tiny
%     % \scriptsize
%     % \footnotesize  % Slightly smaller to better fit
%     % \small
%     % \normalsize  % (default size) 
%     % \large
%     \caption{Comparison of computational time of FSTSP on Set1}
%     \label{tab:set1_timecompare}
%     \begin{tabular*}{\textwidth}{@{\extracolsep{\fill}}l@{\;}c@{\;}c@{\;}c@{\;}c@{}}
%         \toprule
%         \textbf{Method} & \textbf{10-cus (20)} & \textbf{10-cus (40)} &
%         \textbf{20-cus (20)} & \textbf{20-cus (40)} \\
%         \midrule
%         \cite{boccia2023new} & 1.34 & 2.42 & 77.36 & 535.44 \\
%         LKH Config1 & 45.90 & 148.83 & 235.69 & 369.81 \\
%         \bottomrule
%     \end{tabular*}
% \end{table}

\item When solving small poi-7- instances optimally (Tables~\ref{tab:fstsp_1LKH_small}, \ref{tab:tspmd_1LKH_small}), CPLEX requires on average more time as problem complexity increases: from 2.99s for simple FSTSP to 9.12s for complex TSP-mD. In contrast, LKH Config1 requires less time, decreasing from 6.24s to 1.09s. As for the poi-10- instances, CPLEX takes 824.70s for FSTSP and cannot obtain optimal solutions within the 4h time limit for TSP-mD, whereas LKH reduces the average computational time from 53.83s (FSTSP) to 20.26s (TSP-mD). This improvement from FSTSP to TSP-mD is likely attributed to the design logic of the penalty function employed in LKH. 

\item The min-cost VRP-D shows distinct patterns compared to min-makespan FSTSP and TSP-mD. First, improved solutions given an extended runtime often coincide with fewer drone-served customers (Table~\ref{tab:vrpd_1LKH_long}), a pattern not observed in FSTSP (Table~\ref{tab:fstsp_1LKH_longer}) and TSP-mD (Table~\ref{tab:tspmd_1LKH_longer}). Second, removing the node revist restriction enables LKH Config2 to outperform Config1 in VRP-D (Table~\ref{tab:vrpd_2LKH}), while this improvement is not consistently observed for FSTSP (Tables~\ref{tab:fstsp_2LKH_small}, \ref{tab:fstsp_2LKH_medium}) and TSP-mD (Tables~\ref{tab:tspmd_2LKH_small}, \ref{tab:tspmd_2LKH_medium}). 

\end{enumerate}

Furthermore, some operational constraints such as visiting non-customer points can be readily incorporated, as discussed in the subsequent Extensions section.

\subsection{Extensions} \label{sec:extensions}

In this subsection, the flexibility of our framework is further demonstrated through two additional experiments, one modifies the input file \texttt{.drone} for LKH-3, and the other adapts the penalty function ``Penalty\_drone.c". Each experiment incorporates two extensions.

Firstly, for the first phase described in \hyperref[sec:methodology]{Methodology}, when generating input files in Matlab, we can adapt to TSP-D addressed in \cite{poikonen2019branch}, \cite{vasquez2021exact}, and \cite{roberti2021exact}. This generalized coordination of one truck and one drone eliminates the no-wait restriction inherent in FSTSP as mentioned in \hyperref[sec:introduction]{Introduction}, and may be more flexible by allowing nodes to be revisited multiple times in \cite{agatz2018optimization} and \cite{bouman2018dynamic}. In addition, we can easily introduce non-customer sites as candidate L/R locations for drones (\cite{carlsson2018coordinated}, \cite{madani2024hybrid}).

\begin{figure}[H]
    \centering
    \setlength{\tabcolsep}{1pt} % Space between columns
    \begin{tabular}{cc}        \includegraphics[width=0.5\textwidth]{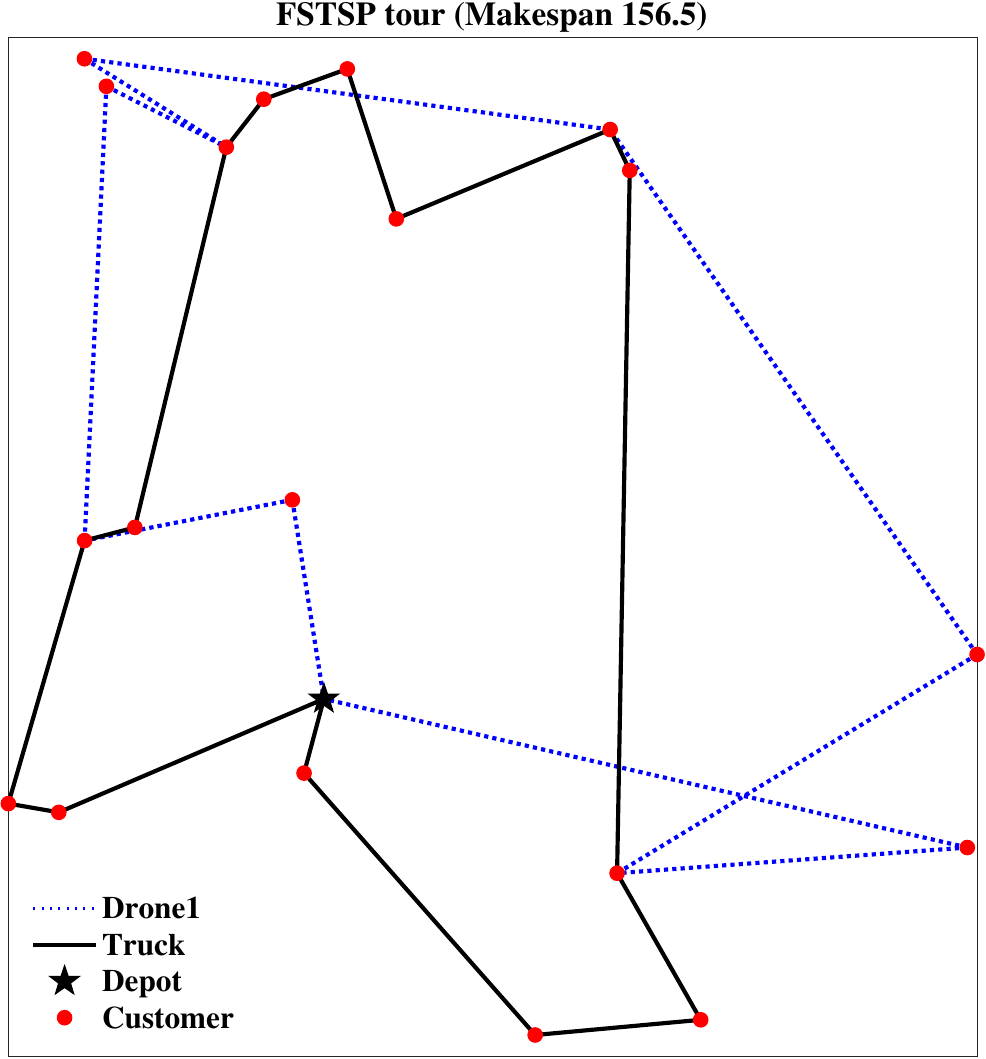} &        
    \includegraphics[width=0.5\textwidth]{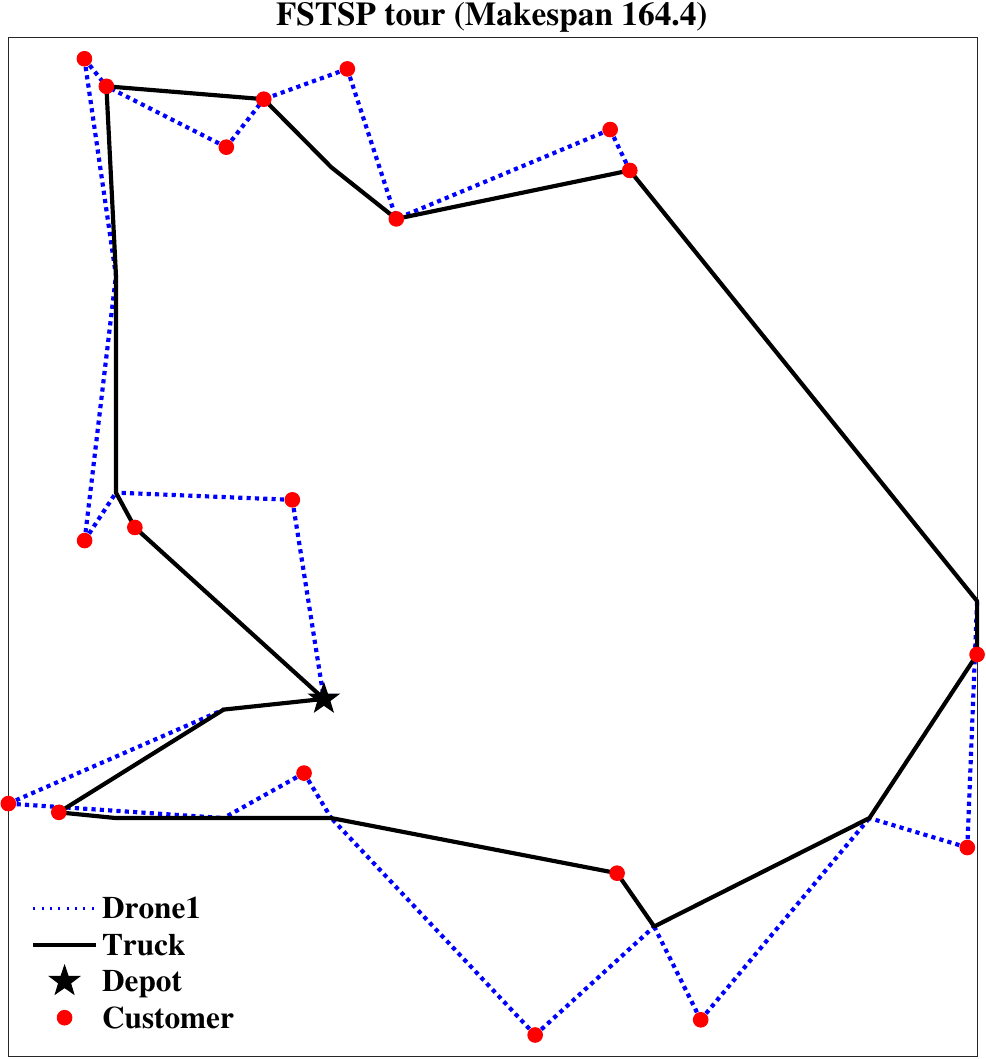} \\       
    \multicolumn{1}{c}{(a) Basic FSTSP.} & \multicolumn{1}{c}{(b) FSTSP extension.}
    \end{tabular}
    \caption{FSTSP tour and its extension on instance poi-20-1.}
    \label{fig:fstsp_loopNon}
\end{figure}

Thus, maintaining the same drone-to-truck speed ratio (50/35), we apply these two extensions to instance poi-20-1 and compare to the basic FSTSP in Fig.~\ref{fig:fstsp_loopNon}. The right solution is generated when allowing the truck to remain stationary at launch points for drone retrieval and enabling multiple visits to both customer and non-customer nodes. However, they do not necessarily improve the basic FSTSP within 5min runtime, which is indicated by the increased makespan of 164.4.

Secondly, according to the second phase described in \hyperref[sec:methodology]{Methodology}, we can accommodate alternative objectives in the penalty function. For example, the objective calculation of 1-1 VRPD can be adjusted from min-cost to min-makespan. Furthermore, the penalty function developed for 1-1 VRP-D can be easily modified to address 1-m VRP-D solved by \cite{gao2023scheduling}, wherein each truck is independently coordinated with multiple drones.

\begin{figure}[H]
    \centering
    \setlength{\tabcolsep}{1pt} % Space between columns
    \begin{tabular}{cc}        \includegraphics[width=0.5 \textwidth]{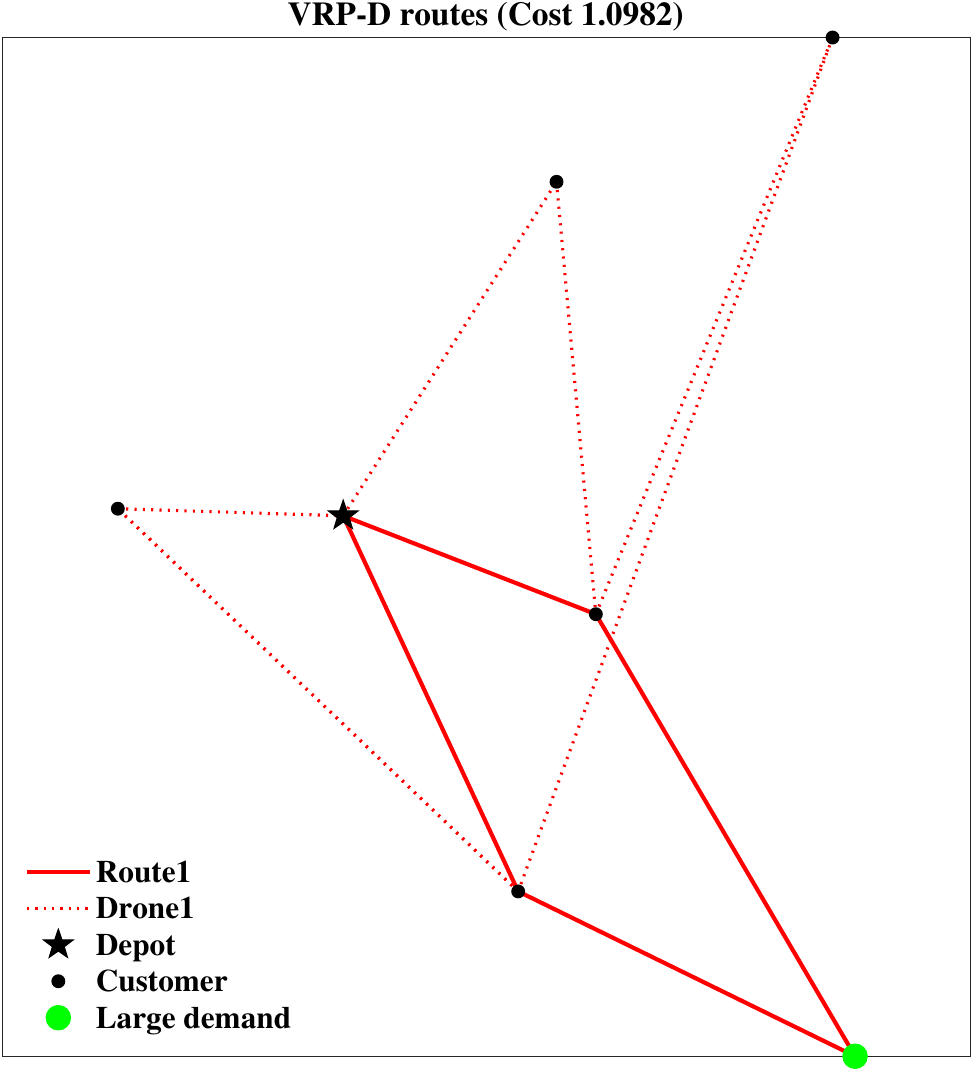} &        
    \includegraphics[width=0.5 \textwidth]{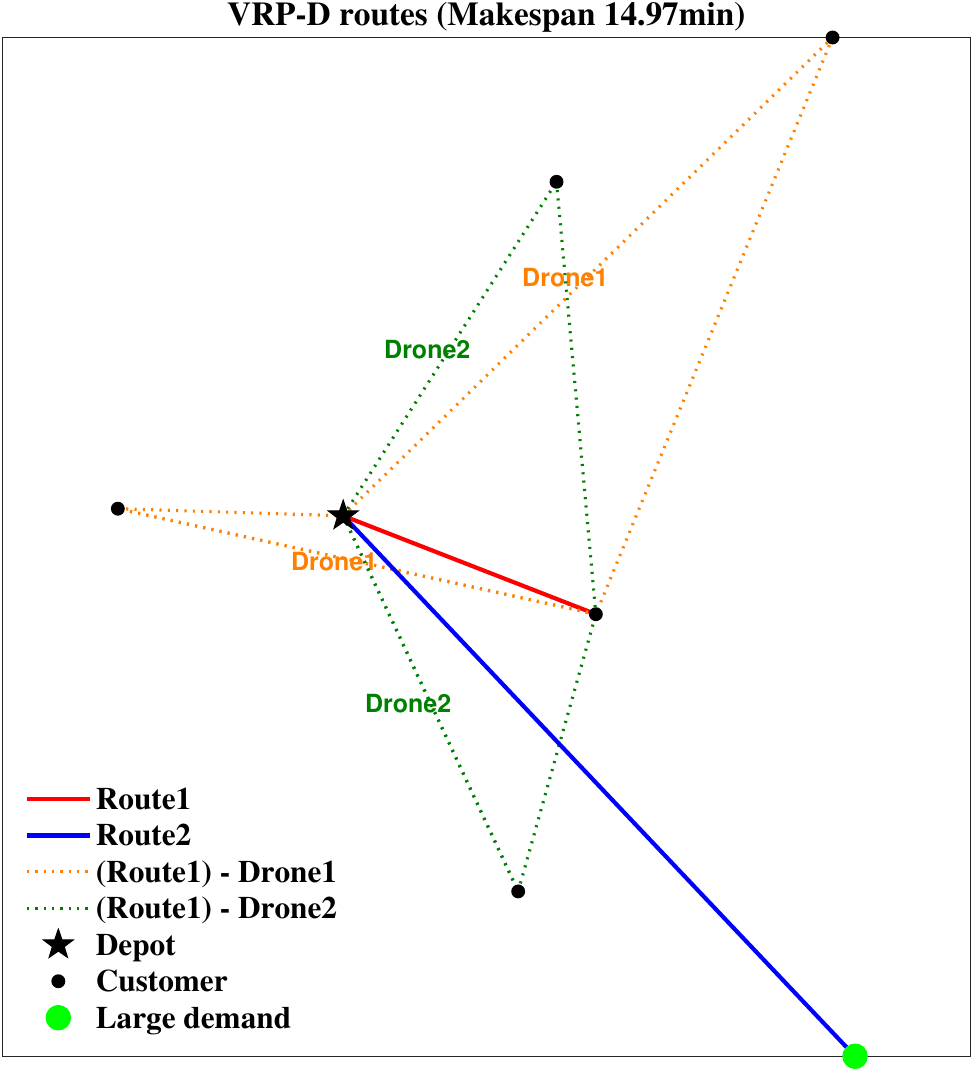} \\       
    \multicolumn{1}{c}{(a) 1-1 VRPD min-cost.} & \multicolumn{1}{c}{(b) 1-m VRPD min-makespan.}
    \end{tabular}
    \caption{VRPD tours on instance 6.5.1.}
    \label{fig:vrpd_costmakespan}
\end{figure}

Thus, we use instance 6.5.1 to solve 1-m VRP-D with min-makespan objective, as shown in Fig.~\ref{fig:vrpd_costmakespan}(b). Fig.~\ref{fig:vrpd_costmakespan}(a) minimizes the cost of 1-1 VRP-D to 1.0982€ with a corresponding makespan of 29.78min. In the meantime, the right tour costs higher at 1.388€ but with a smaller makespan of 14.97min, where Route 1 takes 14.97min using two drones and Route 2 takes 13.66min without drones.

In conclusion, our computational results and extensions presented here suggest that the penalty-based solution framework could be a basis to flexibly address all types of truck-and-drone routing problems. Incorporating some operational constraints may be easily handled by introducing additional penalty terms. But the core difficulty stems from certain complex problem settings that might change the coordination structure between trucks and drones, which demands tailoring the logic of the penalty function. This indicates several directions for future research, outlined as follows.

\subsection{Future work} \label{sec:future}

\begin{enumerate}
\item Truck-drone assignment flexibility (\cite{wang2019vehicle}, \cite{bakir2020optimizing}): Each drone is allowed to be retrieved by a truck different from the one originally assigned to it.
\item Multi-visit per drone flight (\cite{poikonen2020multi}, \cite{luo2021multi}): Drones can serve multiple customers sequentially before returning to the truck, increasing efficiency but complicating flight planning.
\item Cooperative pickup and delivery (\cite{karak2019hybrid}, \cite{gao2023multi}): Customers may also request package pickup services from their locations.
\item Framework optimization: Attempt to further make use of the high performance of LKH-3 and integrate Machine Learning techniques (\cite{bogyrbayeva2023deep}, \cite{sun2025physics}) to improve computational efficiency.

\end{enumerate}

\bibliographystyle{apalike}
\bibliography{reference}  % or whatever your .bib file is named without the extension
% \section{Appendix} \label{sec:appendix}

% \section*{Appendix}
% \addcontentsline{toc}{section}{Appendix}

\appendix
\renewcommand{\thesection}{Appendix \Alph{section}.}
\section{Input files}
\label{app:Inputfiles}

\begin{figure}[H]
    \centering
    \setlength{\tabcolsep}{1pt} % Space between columns
    \begin{tabular}{c}        \includegraphics[width=0.8\textwidth]{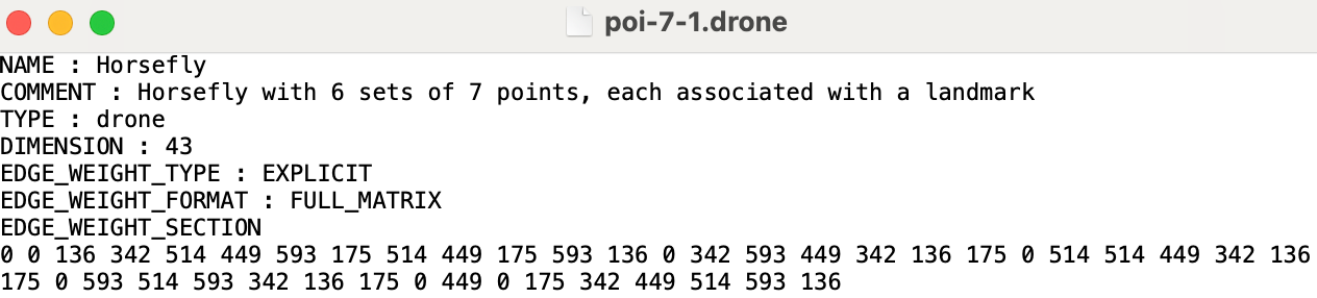} \\
        % \multicolumn{1}{c}{(a) ALNS.}
    \end{tabular}
    \caption{FSTSP input files for instance poi-7-1.}
    \label{fig:fstsp_exe1}
\end{figure}

For FSTSP, instance poi-7-1 includes 1 depot and $n=6$ customer locations. Based on the description in \hyperref[sec:methodology]{Methodology}, each customer $i$ generates a set of \texttt{new\_points} containing $k$ elements. For example, the input parameter $k$ is set to $7$. Then, these 7 points consist of the customer $i$ itself and its $k-1=6$ nearest neighboring nodes. Thus, the primary data, \texttt{flat\_points}, has a dimension of $n \cdot k+1=43$, which is specified as ``DIMENSION: 43" in Fig.~\ref{fig:fstsp_exe1}.

\begin{figure}[H]
    \centering
    \setlength{\tabcolsep}{1pt} % Space between columns
    \begin{tabular}{c}        \includegraphics[width=0.50\textwidth]{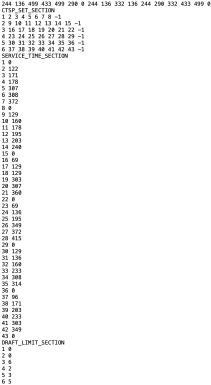} \\
        % \multicolumn{1}{c}{(a) ALNS.}
    \end{tabular}
    \caption{FSTSP input problem file for instance poi-7-1.}
    \label{fig:fstsp_exe2}
\end{figure}

% ``
% ----------
% For poi-7- instances, 6 sets of 7 points (6 nearest nodes + 1 itself). Dimension 43. 

% For poi-10- instances, 9 sets of 10 points (9 nearest nodes + 1 itself). Dimension 91.

% For poi-20- instances, 19 sets of 10 points (9 nearest nodes + 1 itself). Dimension 191.

% For poi-30- instances, 29 sets of 15 points (14 nearest nodes + 1 itself). Dimension 436.

% For poi-40- instances, 39 sets of 20 points (19 nearest nodes + 1 itself). Dimension 781."

\begin{figure}[H]
    \centering
    \setlength{\tabcolsep}{1pt} % Space between columns
    \begin{tabular}{c}        \includegraphics[width=0.8\textwidth]{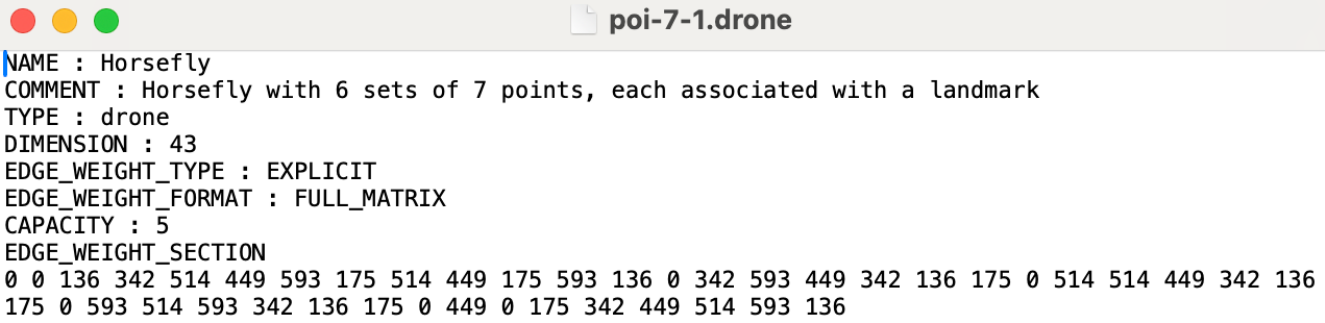} \\
        % \multicolumn{1}{c}{(a) ALNS.}
    \end{tabular}
    \caption{TSP-mD input files.}
    \label{fig:tspmd_exe1}
\end{figure}

The TSP-mD input files follow the same format as the FSTSP, with an additional line. For example, ``CAPACITY: 5" is added to specify the number of available drones, as shown in file \texttt{.drone} (Fig.~\ref{fig:tspmd_exe1}).

% ``----------
% For poi-7- instances, 6 sets of 7 points (6 nearest nodes + 1 itself). Dimension 43. 

% For poi-10- instances, 9 sets of 10 points (9 nearest nodes + 1 itself). Dimension 91.

% For poi-20- instances, 19 sets of 10 points (9 nearest nodes + 1 itself). Dimension 191.

% For poi-30- instances, 29 sets of 15 points (14 nearest nodes + 1 itself). Dimension 436.

% For poi-40- instances, 39 sets of 20 points (19 nearest nodes + 1 itself). Dimension 781."

\begin{figure}[H]
    \centering
    \setlength{\tabcolsep}{1pt} % Space between columns
    \begin{tabular}{c}        \includegraphics[width=0.8\textwidth]{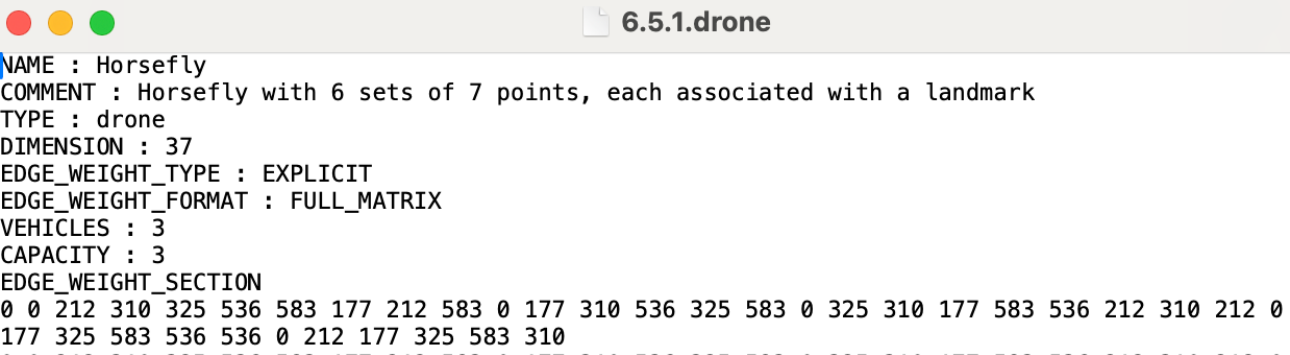} \\
        % \multicolumn{1}{c}{(a) ALNS.}
    \end{tabular}
    \caption{VRP-D input files.}
    \label{fig:vrpd_exe1}
\end{figure}

The VRP-D input files follow the same format as TSP-mD, with an additional line. As illustrated in Fig.~\ref{fig:vrpd_exe1}, ``VEHICLES: 3" is added to specify the number of available trucks. The existing line ``CAPACITY: 1" indicates that each truck carries a single drone.

For example, instance 6.5.1 includes 1 depot and 6 customer locations. Following \hyperref[sec:methodology]{Methodology}, each customer generates a set of points, i.e., \texttt{new\_points}, yielding 6 sets in total. The high-demand customer refers to a set containing only its own location, while each of the other 5 customers derives a set comprising the customer point plus its 6 nearest points for L/R operations. Therefore, the total number of \texttt{flat\_points} is calculated as 1 truck-only customer + 5 sets of 7 elements + 1 depot = 37, which matches the ``DIMENSION: 37" stated in the input file ``6.5.1.drone".

% ``----------
% For 6- instances, 6 sets of 7 points (6 nearest nodes + 1 itself). 6.5.1 has Dimension 37 = 5*7 + (6-5)truck-only customers + 1depot. VEHICLES 3.

% For 10- instances, 10 sets of 10 points (9 nearest nodes + 1 itself). 10.5.1 has Dimension 56 = 5*10 + (10-5)truck-only customers + 1depot. VEHICLES 10.

% For 12- instances, 12 sets of 10 points (9 nearest nodes + 1 itself). 12.5.1 has Dimension 94 = 9*10 + (12-9)truck-only customers + 1depot. VEHICLES 10.

% For 20- instances, 20 sets of 10 points (9 nearest nodes + 1 itself). 20.5.1 has Dimension 156 = 15*10 + (20-15)truck-only customers + 1depot. VEHICLES 10."

% \section{Penalty functions}
% \label{app:Penaltyfunctions}

\end{document}